\documentclass{amsart}

\usepackage{amsmath,amssymb,mathrsfs,amsthm}

\usepackage{enumerate}
\usepackage{mathtools}\mathtoolsset{showonlyrefs=true} 
\usepackage{docmute} 

\usepackage{tikz-cd}
\usetikzlibrary{arrows.meta, positioning, calc}

\usepackage{xspace} 

\usepackage[utf8]{inputenc} 
\usepackage{graphicx} 

\usepackage[expansion=false]{microtype}
\usepackage[mathcal]{euscript}
\usepackage{mathalfa}

\usepackage{bm} 
\usepackage[all]{xy}
\usepackage{comment}

\usepackage{color}

\numberwithin{equation}{section}
\usepackage{hyperref}

\newcommand{\nc}{\newcommand}
\newcommand{\oname}{\operatorname}

\nc{\myfavorite}{\mathbb} 

\newcommand{\C}{\myfavorite{C}}

\newcommand{\E}{\myfavorite{E}}
\newcommand{\F}{\myfavorite{F}}

\newcommand{\N}{\myfavorite{N}}
\newcommand{\Q}{\myfavorite{Q}}
\newcommand{\R}{\myfavorite{R}}
\newcommand{\T}{\myfavorite{T}}

\newcommand{\Z}{\myfavorite{Z}}
\nc{\RZ}{\R/\Z}
\nc{\QZ}{\Q/\Z}

\newcommand{\myb}{\mathbb}

\newcommand{\bC}{\myb C}

\newcommand{\bP}{\myb P}

\newcommand{\bR}{\myb R}

\newcommand{\bZ}{\myb Z}

\newcommand{\bbC}{\mathbb{C}}

\newcommand{\bbR}{\mathbb{R}}

\nc{\mcal}{\mathcal}
\nc{\calA}{\mcal A}
\nc{\calB}{\mcal B}
\nc{\calC}{\mcal C}
\nc{\calD}{\mcal D}
\nc{\calE}{\mcal E}
\nc{\calF}{\mcal F}
\nc{\calG}{\mcal G}
\nc{\calH}{\mcal H}
\nc{\calI}{\mcal I}
\nc{\calJ}{\mcal J}
\nc{\calK}{\mcal K}
\nc{\calL}{\mcal L}
\nc{\calM}{\mcal M}
\nc{\calN}{\mcal N}
\nc{\calO}{\mcal O}
\nc{\calP}{\mcal P}
\nc{\calQ}{\mcal Q}
\nc{\calR}{\mcal R}
\nc{\calS}{\mcal S}
\nc{\calT}{\mcal T}
\nc{\calY}{\mcal U}
\nc{\calV}{\mcal V}
\nc{\calW}{\mcal W}
\nc{\calX}{\mcal X}
\nc{\calZ}{\mcal Z}

\DeclareMathSymbol{A}{\mathalpha}{operators}{`A}
\DeclareMathSymbol{B}{\mathalpha}{operators}{`B}
\DeclareMathSymbol{C}{\mathalpha}{operators}{`C}
\DeclareMathSymbol{D}{\mathalpha}{operators}{`D}
\DeclareMathSymbol{E}{\mathalpha}{operators}{`E}
\DeclareMathSymbol{F}{\mathalpha}{operators}{`F}
\DeclareMathSymbol{G}{\mathalpha}{operators}{`G}
\DeclareMathSymbol{H}{\mathalpha}{operators}{`H}
\DeclareMathSymbol{I}{\mathalpha}{operators}{`I}
\DeclareMathSymbol{J}{\mathalpha}{operators}{`J}
\DeclareMathSymbol{K}{\mathalpha}{operators}{`K}
\DeclareMathSymbol{L}{\mathalpha}{operators}{`L}
\DeclareMathSymbol{M}{\mathalpha}{operators}{`M}
\DeclareMathSymbol{N}{\mathalpha}{operators}{`N}
\DeclareMathSymbol{O}{\mathalpha}{operators}{`O}
\DeclareMathSymbol{P}{\mathalpha}{operators}{`P}
\DeclareMathSymbol{Q}{\mathalpha}{operators}{`Q}
\DeclareMathSymbol{R}{\mathalpha}{operators}{`R}
\DeclareMathSymbol{S}{\mathalpha}{operators}{`S}
\DeclareMathSymbol{T}{\mathalpha}{operators}{`T}
\DeclareMathSymbol{U}{\mathalpha}{operators}{`U}
\DeclareMathSymbol{V}{\mathalpha}{operators}{`V}
\DeclareMathSymbol{W}{\mathalpha}{operators}{`W}
\DeclareMathSymbol{X}{\mathalpha}{operators}{`X}
\DeclareMathSymbol{Y}{\mathalpha}{operators}{`Y}
\DeclareMathSymbol{Z}{\mathalpha}{operators}{`Z}

\nc{\mrm}{\mathrm}

\nc{\ideal}{\mathfrak}
\newcommand{\ideala}{\ideal a}
    \newcommand{\ida}{\ideala}
\newcommand{\idealb}{\ideal b}
    \newcommand{\idb}{\idealb}
\newcommand{\idealc}{\ideal c}
    \newcommand{\idc}{\idealc}
\newcommand{\ideald}{\ideal d}
    \newcommand{\idd}{\ideald}
\newcommand{\idealn}{\ideal n}
    \newcommand{\idn}{\idealn}
\newcommand{\idealp}{\ideal p}
    \newcommand{\idp}{\ideal p}
\newcommand{\idealq}{\ideal q}
    \newcommand{\idq}{\idealq}
    \newcommand{\idealqSiegel}{\idealq_{\mrm{Siegel}}}
\newcommand{\idealm}{\ideal m}

\nc{\wtil}{\widetilde} 
\nc{\wti}{\widetilde} 
\nc{\super}[1]{^{(#1)}}

\nc{\ssec}{\subsection} 
\nc{\subsec}{\subsection}

\nc{\sssec}{\subsubsection}
\nc{\subsubsec}{\sssec}

\nc{\te}{\text}

\nc{\baci}{^{\times }} 

\nc{\numberOfTerms}[2]{\underset{#2}{\underbrace{#1}}}

\newcommand{\lnorm}[1]{\left\Vert #1 \right\Vert} 

\newcommand{\Gnorm}[2]{\left\Vert #1 \right\Vert _{U^{s+1}(#2)}} 

\newcommand{\mgn}[1]{\left| #1 \right|} 

\nc{\br}[1]{\left\{ #1 \right\} } 
\nc{\paren}[1]{\left( #1 \right)} 

\nc{\Ideles}{\mathbf I_K} 
\nc{\IdelesK}{\Ideles}
\nc{\idele}{id{\`e}le\spacex} 
\nc{\ideles}{id{\`e}les\spacex}
\nc{\adele}{ad{\`e}le\spacex}
\nc{\adeles}{ad{\`e}les\spacex}

\nc{\Primes}{\calP} 
\nc{\primes}{\Primes}

\nc{\Cl}{\oname{Cl}} 
\nc{\ClK}{\oname{Cl}(K)} 

\nc{\Nrm}{\mathbf{N}} 

\nc{\OK}{\mathcal O_K}
\nc{\Ok}{\OK}
\nc{\OKS}{\mathcal O_{K}[S^{-1}]} 
    \nc{\idaS}{\mathfrak a [S^{-1}]} 
    \nc{\idas}{\idaS} 
\nc{\OKN}{\mathcal O_{K,\le N}}
\nc{\OKp}{\mathcal O_{K,\mathfrak p}}

\nc{\Zn}{\Z ^n} 
\nc{\Zd}{\Z ^d} 
\nc{\Fp}{\F _p}

\nc{\totient}{\varphi _K}
\nc{\muK}{\mu _K} 
    \nc{\muk}{\mu _K} 

\nc{\setm}{\setminus}

\nc{\nonzero}{\setminus\{0\}}

\nc{\sgn}{\oname{sgn}}

\nc{\Kp}{K_{\idealp}}
\nc{\kp}{\Kp}

\nc{\Spec}{\oname{Spec}}
\nc{\SpecOK}{\oname{Spec}\OK}
\nc{\specOK}{\SpecOK}
\nc{\spOK}{\SpecOK}
\nc{\SpecK}{\SpecOK} 
\nc{\spK}{\SpecOK} 

\nc{\Ideals}{\oname{Ideals}}
\nc{\IdealsK}{\Ideals_K}
\nc{\IdealsKN}{\Ideals_{K,\le N}}

\nc{\Lambdaa}{\Lambda ^{\ideala }_K}
    \nc{\wtiLambdaa}{\widetilde{\Lambda} ^{\ideala }_K}
\nc{\LambdaCramer}{\Lambda _{\mrm{Cram{\acute e}r} ,Q}}
    \nc{\LambdaCramerIntro}{\Lambda _{\mrm{Cram{\acute e}r},Q }}
\nc{\LambdaaCramer}{\Lambda ^{\ideala }_{\mrm{Cram{\acute e}r} ,Q}}
    \nc{\LambdaaCramerw}{\Lambda ^{\ideala }_{\mrm{Cram{\acute e}r} ,w}}
\nc{\LambdaSiegel}{\Lambda _{\mrm{Siegel},Q}}
    \nc{\LambdaSiegelIntro}{\Lambda _{\mrm{Siegel},Q}}
\nc{\LambdaaSiegel}{\Lambda ^{\ideala}_{\mrm{Siegel},Q}}
    \nc{\wtiLambdaaSiegel}{\widetilde{\Lambda} ^{\ideala}_{\mrm{Siegel},Q}}

\nc{\psiSiegel}{\psi_{\mrm{Siegel}}}
\nc{\betaSiegel}{\beta_{\mrm{Siegel}}}
\nc{\idqSiegel}{\mathfrak{q}_{\mrm{Siegel}}}

\nc{\Pic}{\oname{Pic}}

\nc{\blog}{\mathbf{log}} 
\nc{\barg}{\mathbf{arg}} 

\nc{\vol}{\oname{Vol}}
    \nc{\Vol}{\vol}
\nc{\area}{\oname{Area}}
\DeclareMathOperator*{\residue}{res}
\nc{\residuezeta}{\residue\limits_{s=1}(\zeta _K(s))}

\nc{\Supp}{\oname{Supp}}

\nc{\vectorOneTwo}{(\numberOfTerms{1,\dots,1}{r_1}, 
\numberOfTerms{2,\dots ,2}{r_2} )}

\nc{\bmv}{\bm{v}} 
\nc{\bmw}{\bm{w}} 
\nc{\bmx}{\bm{x}} 
\nc{\bmxi}{\bm{\xi}} 
\nc{\ebmxi}{e_{\bm{\xi}}} 
\nc{\etheta}{e^{2\pi i \theta (-)}} 

\nc{\innpro}{\cdot} 
\nc{\innpr}{\innpro} 

\nc{\dlog}[1]{\frac{#1 '}{#1}} 
\nc{\ufrac}[1]{\frac{1}{#1}} 
\nc{\fracu}{\ufrac} 

\nc{\Aff}{\oname{Aff}}

\nc{\abspseudopolydecay}[1]{e^{-\sqrt{\log #1}}} 
\nc{\Kpseudopolydecay}[1]{e^{-\sqrt{\log #1}/O_K(1)}} 

\nc{\pseudopolygrowth}[1]{e^{\sqrt{\log #1}/O(1)}} 
\nc{\abspseudopolygrowth}[1]{e^{\sqrt{\log #1}}} 
\nc{\Kpseudopolygrowth}[1]{e^{\sqrt{\log #1}/O_K(1)}} 

\nc{\mydecayexponent}{\rho}
\nc{\mypseudopolygrowth}{\exp ((\log N)^{O_{k,n}(\mydecayexponent )})}
\nc{\mypseudopolydecay}{\exp (-(\log N)^{O_{k,n}(\mydecayexponent )})}
\nc{\pseudopoly}[1]{\exp ((\log N)^{#1})}
\nc{\pseudopolydecay}[1]{\exp (-(\log N)^{#1})}

\nc{\KR}{K\otimes _\Q \R } 
\nc{\Rpos}{\R\baci_{>0}} 
\nc{\OurRegion}{(\KR )_{<X_1,\dots ,X_{r_1+r_2}}}

\nc{\comp}{\oname{comp}} 
\nc{\step}{\oname{step}} 
\nc{\Cinfty}{{C}^\infty}
\nc{\nilcounter}{s}
\nc{\etacounter}{m}
\nc{\oneortwo}{b} 

\nc{\diag}{\oname{diag}} 
\nc{\bdiag}{\oname{\mathbf{diag}}} 

\nc{\id}{\oname{id}} 
\nc{\pr}{\oname{pr}} 

\nc{\can}{\mrm{can}} 

\nc{\tors}{\mrm{tors}} 
\nc{\add}{\mrm{add}} 
\nc{\mult}{\mrm{mult}} 

\nc{\directsum}{\bigoplus}
\nc{\dsum}{\directsum}

\nc{\disjointunion}{\bigsqcup} 
\nc{\dunion}{\disjointunion}

\nc{\coker}{\oname{coker}}
\nc{\rank}{\oname{rank}}

\nc{\inj}{\hookrightarrow} 
\nc{\injfrom}{\hookleftarrow} 

\nc{\surj}{\twoheadrightarrow} 
\nc{\surjfrom}{\twoheadleftarrow} 

\nc{\actson}{\curvearrowright} 
\nc{\acts}{\actson}

\nc{\actedby}{\curvearrowleft} 
\nc{\acted}{\actedby}

\nc{\xto}{\xrightarrow}
\nc{\xfrom}{\xleftarrow}

\nc{\isoto}{\xrightarrow\cong}
\nc{\isofrom}{\xrightarrow[\cong]{}}

\nc{\too}{\longrightarrow}

\nc{\inv}{^{-1}} 

\nc{\ol}{\overline}

\nc{\undl}{\underline}
\nc{\undli}{\underline i}
\nc{\undlj}{\underline j}
\nc{\undlx}{\underline x}

\nc{\noin}{\not\in}

\nc{\nearlyeq}{\fallingdotseq} 
\nc{\nearlyequal}{\nearlyeq}

\nc{\Podual}{^\wedge} 
\nc{\podual}{^\wedge} 

\nc{\Hom}{\oname{Hom}}
\nc{\Map}{\oname{Map}}
\nc{\GL}{\oname{GL}}

\nc{\cube}{\square} 

\nc{\divides}{\mid} 
\nc{\notdivide}{\nmid} 
\nc{\notdiv}{\notdivide} 
\nc{\doesnotdivide}{\notdivide} 

\nc{\mtx}[1]{\begin{pmatrix}#1\end{pmatrix} }

\nc{\hatT}{\widehat\T} 

\nc{\expexp}{\exp\exp}
\nc{\Expec}{\mathop{\mathbb{E}}} 

\nc{\spacex}{\xspace} 

\nc{\Cramer}{Cram{\'e}r\spacex}
    \nc{\mCramer}{\mrm{Cram{\acute e}r}}
\nc{\CauSch}{Cauchy--Schwarz\xspace}
\nc{\CauSchGow}{Cauchy--Schwarz--Gowers\xspace}
\nc{\Frohlich}{Fr{\"o}hlich\xspace}
\nc{\Froh}{\Frohlich}
\nc{\Fou}{Fourier\spacex}
\nc{\GTZ}{Green--Tao--Ziegler\spacex}
\nc{\Holder}{H{\"o}lder\xspace}
\nc{\Lie}{\oname{Lie}} 
\nc{\Lip}{Lipschitz\spacex}
    \nc{\mLip}{\oname{Lip}} 
\nc{\Malcev}{Mel'cev\spacex}
    \nc{\Mal}{\Malcev}
\nc{\Mit}{Mitsui\xspace}
\nc{\Mobius}{M{\"o}bius\spacex}
\nc{\mSiegel}{\mrm{Siegel}}
\nc{\Tera}{Ter{\"a}v{\"a}inen\spacex}
\nc{\vMang}{von Mangoldt\spacex}
\nc{\Yil}{Y{\i}ld{\i}r{\i}m\spacex}

\nc{\red}{\color{red}}
\nc{\brown}{\color{brown}}
\nc{\blue}{\color{blue}}

\nc{\ordinalth}{{\hspace{0pt}}^{\text{th}}} 
    \nc{\ordinalst}{{\hspace{0pt}}^{\text{st}}} 
    \nc{\ordinalnd}{{\hspace{0pt}}^{\text{nd}}} 
    \nc{\ordinalrd}{{\hspace{0pt}}^{\text{rd}}} 

\nc{\cont}{continuous\xspace}
\nc{\conti}{\cont}

\nc{\homo}{homomorphism\xspace}

\nc{\almbda}{\lambda} 

\nc{\loccit}{{\it loc.\,cit.}\xspace}
\nc{\midsep}{\,:\, } 

\nc{\disp}{\displaystyle}

\nc{\coeff}{coefficient\xspace}
\nc{\coeffs}{coefficients\xspace}

\nc{\poly}{\oname{poly}}
\nc{\incl}{\oname{incl}}

\theoremstyle{plain}

\newtheorem{theorem}{Theorem}[section]

\newtheorem{proposition}[theorem]{Proposition}
\newtheorem{prop-def}[theorem]{Proposition-Definition}

\newtheorem{lemma}[theorem]{Lemma}

\newtheorem{corollary}[theorem]{Corollary}

\newtheorem{claim}[theorem]{Claim}

\newtheorem{situation}[theorem]{Situation}

\newtheorem{strategy}[theorem]{Strategy}

\newtheorem*{theorem*}{Theorem}
\newtheorem*{thm*}{Theorem}
\newtheorem*{proposition*}{Proposition}
\newtheorem*{prop*}{Proposition}
\newtheorem*{lemma*}{Lemma}
\newtheorem*{lem*}{Lemma}
\newtheorem*{corollary*}{Corollary}
\newtheorem*{cor*}{Corollary}
\newtheorem*{sublemma*}{Sublemma}
\newtheorem*{sublem*}{Sublemma}
\newtheorem*{hypothesis*}{Hypothesis}
\newtheorem*{hyp*}{Hypothesis}
\newtheorem*{claim*}{Claim}
\newtheorem*{clm*}{Claim}

\theoremstyle{definition}
\newtheorem{definition}[theorem]{Definition}

\newtheorem*{definition*}{Definition}
\newtheorem*{defn*}{Definition}
\newtheorem*{setup*}{Set-up}

\theoremstyle{remark}
\newtheorem{remark}[theorem]{Remark}

\newtheorem*{remark*}{Remark}
\newtheorem*{rem*}{Remark}

\newtheorem*{disclaimer*}{Disclaimer}
\usepackage{Leitfaden_2026_2}

\date{\today}

\begin{document}

\title{Linear patterns of prime elements in number fields}
\author{Wataru Kai}
\address{Mathematical Institute, Tohoku University, Aramaki Aoba 6-3, 980-8578 Sendai, Japan}
\email{kaiw@tohoku.ac.jp}

\dedicatory{In memory of Prof.\ Noriyuki Suwa (1955--2022)} 

\begin{abstract}
    We prove a number field analogue of the Green--Tao--Ziegler theorem on simultaneous prime values of degree 1 polynomials whose linear parts are pairwise linearly independent.
        Applications of our results include 
        a Hasse principle of rational points for certain fibrations 
    $X\to \mathbb{P}^1$ over number fields $K$ which had only been available over $\Q $ by Harpaz--Skorobogatov--Wittenberg,
    and construction of elliptic curves having some specified ranks due to Koymans--Pagano and Zywina.
    This latter family of results led to a negative answer to a generalized Hilbert Tenth Problem.
    %
\end{abstract}

\maketitle

\tableofcontents

\section{Introduction}

Schinzel and Sierpinski \cite{Schinzel-Sierpinski} conjectured that if finitely many (distinct) irreducible $\Z $-coefficient polynomials $f_1,\dots ,f_t\in \Z [X]$ satisfy the condition that there is no prime $p$ that divides the product $f_1(n)\cdots f_t(n)$ for all $n\in \Z $ (the so-called {\em Bouniakowsky condition}), then there should be infinitely many $n\in \Z $ such that the values $f_1(n), \dots , f_t(n)$ are simultaneously prime (positive or negative).
Bateman and Horn 
\cite{Bateman-Horn} 
give the conjectural asymptotic frequency of this event: on the interval $[1,N]\subset \N $ it should happen 
\begin{align}
    \frac{C_{f_1,\dots ,f_t}}{\deg (f_1)\dots \deg (f_t)}
    \cdot
    \frac{N}{(\log N)^t}
    +
    o\left( 
    \frac{N}{(\log N)^t}
     \right)
\end{align}
times, where letting $0\le \omega (p)\le p$ be the number of elements $a\in \Fp$ such that $\prod _{i=1}^t f_i(a)=0$ in $\Fp$, we define 
\begin{align}
    C_{f_1,\dots ,f_t}
    := 
    \prod _{p} 
    \left( \frac{p}{p-1}
        \right) ^t 
        \left( 1-\frac{\omega (p)}{p} \right)
        .
\end{align}
Their conjectures extend the older conjectures by Dickson \cite{Dickson}
and by Hardy and Littlewood 
\cite{Hardy-Littlewood} 
about the case where $f_i$ all have degree $1$.
The twin prime conjecture is a special case of these conjectures, $f_1=X$, $f_2=X+2$.
The only known case of these conjectures is when $t=1$ and $f_1(X)=aX+b$ where $a\neq 0,b\in \Z $ are coprime.

It becomes easier to attain simultaneous prime values when more variables are allowed.
A milestone theorem of Green--Tao--Ziegler \cite{LinearEquations, MobiusNilsequences, GowersInverse} 
affirms the multi-variable version of the Hardy--Littlewood conjecture 
under a ``finite complexity'' assumption:
namely, they established an asymptotic formula for how frequently finitely many degree $1$ polynomials $\psi _1,\dots ,\psi _t \in \Z [X_1,\dots ,X_d]$ attain simultaneous prime values
under the additional condition that 
their linear parts 
$\dot \psi _i $ (obtained by discarding the constant terms) are pairwise linearly independent.
Note that this linear independence condition requires at least two variables and the twin prime conjecture is beyond the scope their result.

In this paper, we establish an analogue of the theorem of Green--Tao--Ziegler for number fields $K$.
To discuss the asymptotic aspect of the problem, fix a norm $\lnorm - \colon K\to \R _{\ge 0}$ in the $\R $-linear algebra sense: namely it satisfies the triangle inquality, non-degeneracy and $\lnorm{ax}=|a|\cdot \lnorm x$ for all $a\in \Q$ and $x\in K$, where $|a|$ is the absolute value as a real number (not any $p$-adic). 
For example one can consider the $l^\infty $-norm associated to a chosen $\Q $-linear basis $K\cong \Q^n$ (where $n:= [K:\Q ]$), 
or consider the canonical norm $\lnorm -_{\can }$
recalled in \S \ref{appendix:norm-length}.
By adding a subscript 
\begin{align}
    \OKN \subset \OK 
\end{align}
let us mean the set of elements $x$ satisfying $\lnorm x \le N$.
Recall that 
the von Mangoldt function $\Lambda _K $ for $K$ associates to an 
element
$x\in \OK $ the value $\log \Nrm (\idp )$ if the ideal $(x)$ is a power of a non-zero prime ideal $\idp $ with exponent $\ge 1$, and $0$ if not.
As usual, the contribution from prime powers of exponents $\ge 2$ is asymptotically negligible.
In particular, one can consider the usual \vMang function $\Lambda (\Nrm (x\OK ) ) $ of the norm of the principal ideal $x\OK $ if they prefer.

A thumbnail version of our main result is the following:

\begin{theorem}[Theorem \ref{thm:prime-values} and Remark \ref{rem:product-beta}]
    \label{thm:GTZ-number-field}
    Let $\psi _1,\dots ,\psi _t\in \OK [X_1,\dots ,X_d]$ be polynomials of degree $1$
    such that their linear parts $\dot\psi _i$ are pairwise linearly independent over $K$.
    Then we have an asymptotic formula as $N\to +\infty $:
    \begin{align}
        \sum _{\undl x\in \OKN ^d}
        \left( \prod _{i=1}^t 
        \Lambda _K (\psi _i (\undl x))
        \right) 
        = 
        C_{\psi _1,\dots ,\psi _t} \cdot |\OKN |^d
        +o(N^{nd})
        .
    \end{align}
    Here, $C_{\psi _1,\dots ,\psi _t }$ is a constant depending on $\psi _1,\dots ,\psi _t$,
    which is positive if (and only if) the following condition is satisfied:
    \begin{quotation}
        For every prime ideal $\idp \subset \OK $, 
        there is a vector $\undl x_{\idp }\in \OK ^d$ such that 
        $\psi _i (\undl x_{\idp }) \not\in \idp $ for all $1\le i\le t$.        
    \end{quotation}
    (Another equivalent condition is that for every prime number $p$, there is a vector $\undl x_p\in \OK ^d$ such that the norms $\Nrm (\psi _i (\undl x_p))$ are prime to $p$ for all $1\le i\le t$.)
\end{theorem}
In down-to-earth terms, Theorem \ref{thm:GTZ-number-field} says for example 
that 
if we are given a finite set of degree $1$ polynomials in two variables
\begin{align}\label{eq:example-two-variables}
    a_iX+b_iY+ 1 \in \OK [X,Y],\quad i=1,\dots ,t
\end{align}
such that any two vectors $(a_i,b_i)$, $(a_j,b_j)$ are linearly independent in $K^2$,
then there are $(x,y)\in \OK ^2$ such that the values $a_ix+b_iy+ 1$ are all prime elements,
and the number of such $(x,y)\in \OKN ^2$ is asymptotically positively proportional to $N^{2n}/(\log N)^t$.
By taking the norm $\OK \to \Z $ and noting that prime ideals of degree $\ge 2$ are negligibly few, Theorem \ref{thm:GTZ-number-field} gives a family of cases where the multivariate version 
of the Bateman--Horn conjecture holds
(in $dn$ variables ($d\ge 2$), degree $n$, number of polynomials $t\ge 1$).

The same form of asymptotic formula is 
valid when each variable $X_k$ is assumed to vary in a {\em chosen congruence class} modulo some non-zero integer $q _k$ if the definition of the constant $C_{\psi _1,\dots ,\psi _t}$ is changed accordingly. To see this, just consider the polynomials $\psi _i '(\undl X):=\psi _i (a_1+q_1X_1,\dots ,a_d+q_dX_d)$ and apply Theorem \ref{thm:GTZ-number-field} to them.
The criterion for the positivity 
of the leading coefficient 
will now say
that there is no fixed prime divisor for all the values of $\prod _{i=1}^t f(\undlx ) $
even under this congruence constraint.
The definitive form (Theorem \ref{thm:prime-values}) of Theorem \ref{thm:GTZ-number-field} also allows one to specify signs of $\psi _i(\undl x)$ at each real place of $K$ as long as they do not contradict each other.

In our earlier joint work \cite{KMMSY}, we considered polynomials \eqref{eq:example-two-variables} in the case where $a_i=1$ for all $i$.
This previous result is stronger in this particular case in that we can find prime-attaining pairs $(x,y)$ from a smaller set $\OK \times (\Z\nonzero )$ (though the asymptotic frequency has been identified only up to a bounded mutiplicative factor), 
but our current result allows
more general coefficients 
including the constant terms.

An outline of the proof of Theorem \ref{thm:GTZ-number-field} is given in \S \ref{sec:outline}.

\subsection{Applications}
We mention two applications of Theorem \ref{thm:GTZ-number-field}.
The first is about the Hasse principle for rational points on special types of varieties over number fields. 
This application is what motivated our attempt to prove the \GTZ theorem in the generality of number fields.

The second is 
the construction of elliptic curves of specified ranks,
due to Koymans--Pagano and Zywina.
Such constructions can be used to draw a negative answer to a generalized Hilbert's Tenth Problem. 

\subsubsection{Hasse principle for rational points on varieties over number fields}
Harpaz, Skorobogatov and Wittenberg \cite{Harpaz-Skorobogatov-Wittenberg} used the original \GTZ theorem 
for the case of $a_iX+b_iY+c_i \in \Z [X,Y]$
to study $\Q $-raional points on algebraic varieties. 
Their technique allows one to show that a $\Q$-variety $X$ fibered over the projective line $\pi \colon X\to \bP ^1_\Q $ 
satisfies the Hasse principle (in that the Brauer--Manin obstruction is the only obstruction) if generic fibers are known to satisfy the honest Hasse principle and one has control over finitely many bad fibers of $\pi $.
Their results include:

\begin{theorem}[{Harpaz--Skorobogatov--Wittenberg \cite[Corollary 3.4]{Harpaz-Skorobogatov-Wittenberg}}]
    Let $X$ be a smooth proper geometrically integral variety over $\Q $
    equipped with a morphism $\pi \colon X\to \bP ^1_\Q $.
    If the generic fiber of $\pi $ is a Severi-Brauer variety 
    and if all the fibers that are not geometrically integral are $\Q $-fibers,
    then the Brauer--Manin obstruction is the only obstruction 
    for the non-emptiness of $X(\Q )$.
\end{theorem}

This theorem is a somewhat immediate consequence of their main theorem 
\cite[Theorem 3.1]{Harpaz-Skorobogatov-Wittenberg}, which is also limited to (geometrically integral) varieties over $\Q$.
In \S \ref{sec:Hasse} we indicate how to extend it to the general number field $K$.
For this purpose we need a variant of Theorem \ref{thm:GTZ-number-field}
where we can verify under similar assumptions that a given finite set of 
degree $1$ polynomials with coefficients in the localized ring $\OKS $
attain simultaneous prime values in $\OKS $
(Theorem \ref{thm:prime-values-S}).

\subsubsection{Construction of elliptic curves of specified ranks and Hilbert's Tenth Problem}\label{sec:elliptic-curves-rank}

Using Theorem \ref{thm:GTZ-number-field} and its variants,
Koymans--Pagano \cite{KoymansPagano25} and Zywina \cite{Zywina2025a} independently showed the following remarkable result:
\begin{theorem}[{Koymans--Pagano \cite{KoymansPagano25}, Zywina \cite{Zywina2025a}}]\label{thm:rank1}
    Over every number field $K$, there are (infinitely many) elliptic curves $E/K$ such that $E(K)$ has rank $1$ as an abelian group.
\end{theorem}
The rank $0$ case had been known by Mazur--Rubin. 
Zywina \cite{Zywina26} more recently extended this statement to ranks $2,3,4$.

Hilbert's Tenth Problem asked for an algorithm to decide whether a given polynomial $f\in \Z [x_1,\dots ,x_d]$ has solutions in $\Z$.
In 1970, Matiyasevich showed that there is no such algorithm that works for all polynomials (the {\it MRDP Theorem}), building upon earlier work of Davis, Putnam and Robinson.




Soon people began to wonder about the case of other coefficient rings
and later it was shown that a key to extending the negative answer to rings such as $\OK $ was to construct, for quadratic extensions of number fields $L/K$, an abelian variety $A$ over $K$ satisfying the {\em rank stability} condition 
\[ 1\le \rank A(K) =\rank A(L) . \]
This can be achieved by techniques similar to Theorem \ref{thm:rank1}.
Thus:


\begin{theorem}[{Koymans--Pagano \cite{KoymansPagano24}, 
    Alp{\"o}ge--Bhargava--Ho--Shnidman \cite{ABHS}, 
    Zywina \cite{Zywina2025a}}]\label{thm:Hilbert10-number-rings}
    Hilbert's Tenth Problem has a negative answer over $\OK $ for every number field $K$.
\end{theorem}

By Eisentr\"{a}ger, it is known that Theorem \ref{thm:Hilbert10-number-rings} implies the negative answer to Hilbert's Tenth Problem over every commutative ring of finite type over $\Z $ that is an infinite set.
See \cite{KoymansPagano24} for more details and references.

\subsection*{Acknowledgments}
This research is a continuation of my joint work with Masato Mimura, Akihiro Munemasa, Shin-ichiro Seki and Kiyoto Yoshino \cite{KMMSY}, 
through which I learned a lot about additive combinatorics.

I got interested in this particular direction around New Year 2019/2020 when Federico Binda drew my attention to 
the work of Harpaz--Skorobogatov--Wittenberg 
and to the fact that it was limited to $\Q $ due to the lack of a number field analogue of the Green--Tao--Ziegler theorem.
Online discussions with Federico Binda, Hiroyasu Miyazaki and Rin Sugiyama in the subsequent lockdown period were very helpful---this activity was one of my few amusements during that time.

Joni \Tera and James Leng have kindly answered my questions about their respective works \cite{TaoTeravainen} (with Terence Tao) and \cite{EfficientEquidistribution}.

While I was working on this project, I was 
supported by Japan Society for the Promotion of Science (JSPS) through 
JSPS Grant-in-Aid for Young Scientists (JP18K13382, JP22K13886)
and JSPS Overseas Research Fellowship, 
which funded my long-term visit to University of Milan from March 2022 to March 2024.
My colleagues in Milan have always been immensely helpful both mathematically and non-mathematically.

\subsection{Notation and conventions}\label{sec:notation}

$K$ is a number field of degree $n:=[K:\Q ]$.
$N>1$ is the main integer parameter indicating the size of cubes $[-N,N]^n$.
The letter $s\in \N $ will be used to indicate
the following equivalent information: 
the complexity of the set of affine-linear forms $f_1,\dots ,f_t$,
the degree of the Gowers norm $U^{s+1}$,
and the nilpotence step of the nilmanifold $G/\Gamma $. 

\subsubsection{Asymptotic notation}

When $X$ is a positive quantity and $a,b,c,\dots $ are some parameters,
the expression $O_{a,b,c,\dots }(X)$ means a complex number quantity whose magnitude is bounded by $C_{a,b,c,\dots }\cdot X$ where $C_{a,b,c,\dots }>0$ is a positive number which can be written as a function in $a,b,c,\dots $.
The expressions $Y=O_{a,b,c,\dots ,}(X)$ and $Y\ll _{a,b,c,\dots }X$ mean that $Y$ is an $O_{a,b,c,\dots }(X)$.
An assertion of the form 
\begin{quote}
    ``If $Y\ll _{a,b,c,\dots }X$, then $Y'\ll _{a',b',c',\dots }X'$''
\end{quote}
means that there are constants $C_{a,b,c,\dots }>0$ (typically small) and $C_{a',b',c',\dots }>0$ (typically large)
which make the following true:
\begin{quote}
    ``If $|Y|< C_{a,b,c,\dots }X$, then $|Y'|< C_{a',b',c',\dots }X'$.''
\end{quote}

\subsubsection{The level of sieve}

Given $N>1$, let 
$Q$ be a quantity satisfying:
%
%
\begin{align}\label{eq:def-of-Q}
   \exp ((\log N)^{1/100}) < Q < \exp ( (\log N)^{1/3} ).
\end{align}
This will appear when we sieve out small primes $p<Q$.
The reader may fix it at $Q=\exp ((\log N)^{1/10})$, for example, for the most part of this paper.
Let
\begin{align}
    P(Q):= \prod _{0<p<Q} p 
\end{align}
be the product of all positive prime numbers smaller than $Q$.

As in this case, the letter $p$ is reserved for positive prime numbers. In particular, symbols $\sum _p , \prod _p $ will mean sums and products over positive prime numbers (often obeying some additional conditions) even when not explicitly so mentioned.

\subsubsection{Indicator functions of subsets and conditions}

When $A'$ is a subset of a set $A$, denote by $1_{A'}\colon A\to \{ 0,1 \} $
the indicator function of $A'$.

When $\mathbf P(-)$ is a condition for elements of $A$, the indicator function $1_{\mathbf P(-)}\colon A\to \{ 0,1\} $ is defined to take value $1$ when the condition holds and $0$ if not.

Sometimes we write indicator functions as $1[A']$ or $1[\mathbf P (-)]$ for legibility. 

\subsubsection{Cubes in Euclidean spaces}
For a finite-dimensional $\R$-vector space $V\cong \R^n$ with a chosen basis,
we write 
\begin{align}
    [a,b]^n_V ,\ (a,b]^n_V,\ \dots \subset V
\end{align}
for cubes with respect to the basis.
For a finite-rank free abelian group $\ida \cong \Zn$ with a chosen basis,
we write 
\begin{align}
    [a,b]^n_\ida := \ida\cap [a,b]^n_{\ida\otimes \R } 
\end{align}
and so on.
The subscripts will be omitted when there is not so much risk of confusion.
We use the shorthand symbol $[\pm N]$ to mean $[-N,N]$.

\subsubsection{Ideals}
We assume that the reader is familiar with the very basics of algebraic number theory as in \cite[Chapter I]{Neukirch}.
This includes objets like 
$\OK $, $\Cl (K)$ and the theory of fractional ideals.
We write 
\begin{align}
    \IdealsK := \{ \ida\subset \OK \mid \te{$\ida $ is a non-zero ideal} \}
\end{align}
for the multiplicative monoid of non-zero ideals of $\OK $.

We will make routine use of {\em norm-length compatible bases} of ideals, which are explaind in \S \ref{appendix:norm-length}.

Lower case German letters $\ida ,\idb, \idc, \idd ,\idn ,\idq $ are used to denote 
non-zero fractional ideals of $\OK $. 
The letter $\idp $ is reserved for non-zero {\em prime } ideals.
In particular, sums and products of the form $\sum _\idp $, $\prod _\idp $ are always those over prime ideals,
often with some additional conditions indicated.
The adjective ``non-zero'' might be omitted from time to time.

An element $x\in \ida $ of a fractional ideal and an ideal $\idq \subset \OK $ are said to be {\em coprime} in $\ida $
if we have $x\OK + \idq \ida =\ida $.
We say $x\in \ida $ and $\alpha \in \OK $ are coprime if $x$ and $(\alpha )$ are.

See \S \ref{sec:Hecke_chars}
for some conventions related to Hecke characters such as 
$\psi ^{\ida }$ (the restriction of a Hecke character to a fractional ideal) 
and $C_K(\idq )$ (the mod $\idq $ \idele class group).

\section{Outline of proof, Gowers norms}\label{sec:outline}

We broadly follow the overall scheme of Green--Tao \cite{LinearEquations}, 
incorporating later improvements by Tao--\Tera \cite{TaoTeravainen} and Leng \cite{EfficientEquidistribution}.
With $\OK $, additive combinatorics gets a little more complicated because the geometry of $\R^n$ is richer than $\R $.
Algebra gets more complex as well because unlike in $\Z $, elements and ideals are two very different notions---the multiplicative structure is cleaner in the framework of ideals, but we have to make sure that they get along with the additive combinatorics concepts. 

The central combinatorial tool is the {\em Gowers $U^{s+1}$-norm}, where $s\ge 0$ is an integer 
which will eventually be set to be the {\em complexity} (\S \ref{sec:vonNeumann}) of the set $\dot\psi _1,\dots ,\dot\psi _t$ of linear maps.

\begin{definition}[Gowers norm]
    \label{def:Gowers-norm}    
Let $A$ be a finite abelian group and $f\colon A\to \C $ be a function.
The {\em Gowers $U^{s+1}$-norm} $\lnorm{f}_{U^{s+1}(A)}\ge 0$ is defined by 
\begin{align}
    \left( 
        \lnorm{f}_{U^{s+1}(A)}
        \right)^{2^{s+1}}
        := 
        \Expec _{\substack{
            x\in A \\ (h_1,\dots ,h_{s+1})\in A^{s+1}}}
            \left[ 
                \prod _{ \substack{
                    (\omega _1,\dots ,\omega _{s+1}) \\ 
                    \in \{ 0,1 \}^{s+1}
                    }   
                    }
                    \calC ^{\omega _1+\dots +\omega _{s+1}}
                    f (x+ \sum _{i=1}^{s+1}\omega _ih_i ) 
                    \right]
                    ,
                \end{align}
                where $\calC $ is the complex conjugate operation.
                Henceforth, for $\undl h=(h_1,\dots ,h_{s+1})\in A^{s+1}$ and $\undl \omega =(\omega _1,\dots ,\omega _{s+1})\in \{ 0,1 \}^{s+1}$, it is convenient to write $|\undl \omega |:= \omega _1+\dots +\omega _{s+1}$ and 
                \begin{align}
                    \undl \omega \cdot \undl h := \sum _{i=1}^{s+1} \omega _ih_i .
                \end{align}
                
                Also, if $A$ is a finite subset of an abelian group $Z$
                and $f\colon A\to \C $ is a function (or a function whose domain of definition contains $A$), we define 
                $\lnorm f _{U^{s+1}(A)}\ge 0$ by 
                \begin{align}
                    \left( 
                        \lnorm{f}_{U^{s+1}(A)}
                        \right)^{2^{s+1}}
                        := 
                        \Expec _{\substack{
        x\in Z \\ 
        \undl h\in Z^{s+1} \\ 
        \te{such that} \\ 
        x+\undl \omega\cdot \undl h\in A \\ 
        \te{for all }\undl \omega \in \{ 0,1 \} ^{s+1} }
        }
        \left[ 
            \prod _{\undl \omega \in \{ 0,1 \}^{s+1}}
            \calC ^{|\undl \omega | }
            f (x+ \undl \omega\cdot \undl h ) 
            \right] 
            .
        \end{align}
        We will use this norm mainly for $A=\OKN \subset \OK $.
\end{definition}
See \cite[Appendix B]{LinearEquations} for a comprehensive introduction to the theory of Gowers norms.

\subsection{The von Neumann theorem}
The von Neumann theorem of Gowers norms (Appendix \ref{sec:vonNeumann}) says that the sum 
$\sum _{x\in \OKN ^d} \prod _{i=1}^t \Lambda _K (\psi _i (x))$
will change little if we replace the function $\Lambda _K $
by another function which is close in the 
$U^{s+1}(\OKN )$-norm.
Our choice of the replacement for 
$\Lambda _K$ 
is its {\em \Cramer model} $\LambdaCramerIntro $,
\begin{align}
    \LambdaCramer (\alpha ) = \begin{cases}
        \disp
        \frac{1}{\residuezeta }\frac{P(Q)^n}{\totient (P(Q))} & \te{if $\alpha $ is coprime to $P(Q)$ in $\OK $} \\
        0 & \te{if not,} 
    \end{cases}
\end{align}
for which the sum is computable (\S \ref{sec:Cramer-model}).
(Recall 
$\exp ( (\log N)^{1/100} ) <Q<\exp ((\log N )^{1/3})$.)
This reduces the problem to showing that $\lnorm{\Lambda_K - \LambdaCramerIntro }_{U^{s+1}(\OKN )} $ is small: 
\begin{theorem}[Corollary \ref{cor:vonMang-Cramer-close!}]
    \label{thm:Mangoldt-close-intro}
    We have 
    \begin{align}
        \lnorm{\Lambda_K - \LambdaCramerIntro }_{U^{s+1}(\OKN )}
        \ll _{s,K,\lnorm{-},A} (\log N)^{-A} 
    \end{align}
    for any positive number $A>1$.
\end{theorem}

As it turns out from basic properties of the Gowers norms (see \eqref{eq:unnormalized-norm-of-1}), we may assume that the chosen linear norm $\lnorm{-}\colon \OK \to \R _{\ge 0}$ is the $l^\infty $-norm associated to a $\Z $-basis $\OK \cong \Z ^n$, which we fix throughout Introduction.
Thus our target is now the quantity $\lnorm{\Lambda_K - \LambdaCramerIntro }_{U^{s+1}([-N,N]^n )}$.

\subsection{Siegel model}\label{subsubsec:Siegel}
When $K=\Q $, 
Green--Tao \cite{QuadraticUniformity}
directly 
compared $\Lambda $ and $\LambdaCramerIntro $.
Tao--\Tera \cite{TaoTeravainen} clarified that it is more efficient to introduce an intermediate model,
the {\em Siegel model} $\LambdaSiegelIntro $ (\S \ref{sec:Siegel-model}).
This modifies $\LambdaCramerIntro $ by incorporating the error caused by the potential Siegel zeros of $L$-functions of Dirichlet characters. 
This reduced them to the estimate of the two quantities
\begin{align}
\lnorm{\Lambda - \LambdaSiegelIntro }_{U^{s+1}([N])}
\quad 
\te{and} \quad
\lnorm{\LambdaSiegelIntro - \LambdaCramerIntro }_{U^{s+1}([N])}
.
\end{align}
Green--Tao's original computations can be used to bound the first difference, 
while the second difference can be bounded by utilizing the Weil bounds for character sums.
The faster decay achieved by this separation allows one to break the ``sound barrier'' so to speak, which makes the remaining procedure substantially simpler.

The Weil bound argument easily carries over to the number field situation (\S \ref{sec:Cramer-Siegel-close}). 
Thus we are left with the question of showing $\lnorm{\Lambda_K - \LambdaSiegelIntro }_{U^{s+1}([-N,N]^n )}$ is small.

\subsection{Gowers Inverse Theorem}
There is a combinatorial tool called {\em Gowers Inverse Theorem}, which gives a sufficient (and in a sense necessary) condition for the Gowers norm $\lnorm{f}_{U^{s+1}}$ of a given function $f$ to be small.
It produces an upper bound for 
the Gowers $U^{s+1}$-norm of a function $f\colon [-N,N]^n\to \C $ 
out of an upper bound for its correlations with {\em $s$-step nilsequences} $\theta (x)$:
\begin{align}\label{eq:intro-correlation}
    \left| \Expec _{x\in [-N,N]^n} f(x) \theta (x) \right| 
    < \epsilon 
    \quad\Rightarrow\quad  
    \lnorm{f}_{U^{s+1}([ - N,N]^n)} 
    = o_{\epsilon \to 0}(1).
\end{align}
Here {\em $s$-step nilsequences} are functions  
that arise as composites of the form 
\begin{align}
    \theta =F\circ g\colon\quad  \Z ^n \xto g G \surj G/\Gamma \xto F \C 
\end{align}
where $G$ is a connected, simply connected Lie group which is $s$-step nilpotent,
$g$ is a {\em polynomial map}, 
$\Gamma \subset G$ is a cocompact discrete subgroup and $F$ is a \Lip function 
(with respect to an appropriate metric)
whose values have absolute value $\le 1$. 
The precise definitions of these notions,
together with various complexity parameters associated with them, are recalled in \S \ref{sec:nilsequences}.

By convention, $1$-step nilpotence means $G$ is abelian, in which case $G/\Gamma $ is a real torus isomorphic to an $\R^m/\Z^m$. The function $F$ can be well approximated by waves $x\mapsto e^{2\pi i \cdot \xi (x)}$ of various frequencies $\xi \in \Hom _{\mathrm{cont}}(\R^m/\Z^m ,\RZ )\cong \Z^m $.
If moreover $g$ is just an additive map, then we are essentially in the realm of classical Fourier analysis.
The general situation, so-called {\em higher order Fourier analysis}, is a 
branch of analysis/combinatorics which 
has emerged since the late 90's. 
Some of the results that we use 
had been formulated for nilsequences on $\Z $ instead of $\Z^n$ in the literature, so in Appendices \ref{appendix:nilsequences}, \ref{sec:vonNeumann} we give some supplements on how to extend them to $\Z^n$.

The first version of the {\em Gowers Inverse Theorem} (\S \ref{sec:Gowers-inverse}) was 
proven by Green--Tao--Ziegler \cite{GowersInverse} in a qualitative form: that is, with an unspecified decay rate function $o_{\epsilon \to 0}(1)$ in formula \eqref{eq:intro-correlation}. 
Manners \cite{Manners} gave the first quantitative bound for the decay rate.
Leng--Sah--Sawhney \cite{QuasipolynomialInverseGowers} established a drastically stronger bound
which (coupled with Leng's strengthened equidistribution theory \cite{EfficientEquidistribution}) allows us to eliminate the need of the so-called 
{\em pseudorandom majorant} for (the $W$-tricked) $\Lambda _K$.

\subsection{Mitsui's Prime Number 
Theorem}

Now our goal is to establish a decay of the correlation $\left| \Expec _{x\in [-N,N]^n}f(x)\theta (x) \right| $ from \eqref{eq:intro-correlation} where $f:= \Lambda _K - \LambdaSiegel $.
Recall $\theta $ is the composite $\OK \cong \Z ^n\xto g G\surj G/\Gamma \xto F \C $ where $G$ is a nilpotent Lie group.
We shall show a subpolynomial decay for this correlation by induction on the nilpotence step $s$ of $G$.
For the induction it is more convenient to aim for a slightly more general statement:
\begin{equation}\label{eq:intro-slightly-general}
    \frac{1}{N^n} \left\vert \sum _{x \in \Omega \cap (a + P \Z^n)} 
    (\Lambda _K - \LambdaSiegel )(x) F (g (x)) \right\vert \overset{?}{<} \exp (-(\log N)^{c_{s,n}}) 
    , 
\end{equation}  
where there is extra flexibility to choose a convex set $\Omega \subset [-N,N]^n$ and a congruence class $a+P\Z^n $ modulo a sublattice, and the exponent 
$0< c_{s,n}< 1$ only depends on $s,n$.

When $s=0$, the group $G$ is trivial and so $\theta =F\circ g $ is a constant function of absolute value $\le 1$.
In this case
the problem is to approximate the sum 
\[ \sum _{
    x\in \Omega \cap (a +P\Z^n) 
     } 
\Lambda _K (x)\] 
by the more predictable sum for $\LambdaSiegel $.
This problem has been considered by Mitsui \cite{Mitsui} at least for $\LambdaCramer $, building on Landau's Prime Ideal Theorem. 
The version with the Siegel model 
was worked out in our preceding work (\cite{KaiMit}, see also Proposition \ref{prop:essentially-Mitsui}).


This argument also shows that we may assume $F$ has average $0$ on $G/\Gamma $ by subtracting the average from $F$.

\subsection{Vaughan decomposition}
To reduce the step-$s$ case to lower step cases,
we perform a {\em Vaughan decomposition} (\S \ref{sec:Vaughan})
of $\Lambda _K $ and $\LambdaSiegel $. 
This expresses them as the sum of functions behaving similarly to the divisor function $\tau (x):= \sum _{\ida \divides x} 1$ or the twofold divisor function $\sum _{\idb \divides \ida \divides x}1$.
Since the von Mangoldt function $\Lambda _K$ is theoretically well behaved only as a function on 
the monoid $\IdealsK $ of all non-zero ideals,
our Vaughan decomposition inevitably takes place in the context of ideals.

It follows that if we assume the left-hand side of \eqref{eq:intro-slightly-general} is larger than some positive number $\delta >0$, then by the pigeonhole principle we obtain an inequality that look similar to:
\begin{align}\label{eq:intro-goal-divisor}
    \frac{1}{N^n}
    \left| 
        \sum _{x\in \Omega \cap (a+P\OK )}
    \sum _{\ida \divides x}  F(g(x))
    \right|
    &
    > \delta (\log N)^{-O_{n}(1)} 
    \\
    \frac{1}{N^n}
    \left| 
        \sum _{x\in \Omega \cap (a+P\OK )}
    \sum _{\ida \divides x}
    \sum _{\idb \divides \ida } F(g(x))
    \right|
    &
    > \delta (\log N)^{-O_{n}(1)} ,
\end{align}
where one of the sources of the $\log N$ factor is the number of terms involved in the Vaughan decomposion.

\subsection{Equidistribution theory of nilsequences}\label{subsec:equidistribution_theory}

Now it is time to look more closely at the nilsequence $F(g(x))$.
We make use of the 
equidistribution theory of nilsequences, more specifically what is called Factorization Theorem. 
This is a machinery which, 
given a witness 
of the failure of equidistribution of an $s$-step polynomial map $g\colon \Z^n \to G\surj G/\Gamma $,
produces a factorization of $g$ into three polynomial maps
\begin{align}
    g = \varepsilon \cdot g' \cdot \gamma \qquad \te{ pointwise product as maps }\Z^n \to G,
\end{align}
where
\begin{itemize}
    \item $\varepsilon $ varies slowly on $[-N,N]^n$ (in technical terms ``{\em smooth}'');
    \item $g'$ takes values in a ``smaller'' (in the senses below) closed subgroup $G'\subset G$;
    \item $\gamma $ is periodic in $G/\Gamma $ with a not too large period.
\end{itemize}
The exact meaning of ``smaller'' for $G'$ was that $G'$ has strictly smaller {\em dimension} than $G$ in Green--Tao \cite{PolynomialOrbits}, while
Leng \cite{EfficientEquidistribution} improved it to the condition that $G'$ has strictly smaller {\em step of nilpotence} than $G$.

A ``witness of failure of equidistribution'' is given by function(s) on $G/\Gamma $
but its technical definition is 
tricky to explain here. 
To illustrate how the arguments proceed in practice, we suppose that the left-hand side of one of \eqref{eq:intro-goal-divisor} is larger than some positive number $\delta >0$ and perform an ideal-theoretic adaptation of the computations of Green--Tao \cite{QuadraticUniformity}, \cite{MobiusNilsequences}.
This produces lower bounds that look similar to 
\begin{align}
    \left\vert \Expec _{\alpha \in [-\frac N D ,\frac N D ]^n } F(g(\alpha b)) \right\vert ,\ 
    \left\vert \Expec _{\alpha \in [-\frac N D, \frac N D]^n } F(g(\alpha b)) \overline{F(g(\alpha b'))} \right\vert 
    > \exp (-(\log (1/\delta ))^{1/O_{s,n}(1)}) 
\end{align}  
for some value $D\in [1,N^n]$ and 
for a positive portion of
elements $b$ (and $b'$) of norm 
in the range $[D,2D]$.
Since by a preliminary Fourier-analytic argument we may assume $F$ has a non-trivial moderate period and in particular average $0$, 
these lower bounds say that the sample points $\{ g(\alpha b)\} _{\alpha \in [-\frac{N}{D}, \frac{N}{D}]^n} $ and $\{ (g(\alpha b),g(\alpha b')) \} _{\alpha \in [-\frac{N}{D}, \frac{N}{D}]^n}$ fail to capture the average of $F$ (resp.\, of the function $\pr _1^*F\cdot \overline{\pr _2^*F}\colon G\times G\to \C $) correctly. 
This essentially witnesses the failure of equidistribution and triggers Factorization Theorem that gives $g=\varepsilon \cdot g'\cdot \gamma $
with complexity parameters bounded by quasipolynomial quantities in $1/\delta $.

To finish off the induction step, we substitute this factorization  in the original formula \eqref{eq:intro-slightly-general}.
We can regard the polynomial maps $\varepsilon $ and $\gamma $ as virtually constant
by chopping the convex set $\Omega $ and congruence class $a+P\Z ^n$
into (quasipolynomially many) small enough ones.
We are left with the same type of sum as \eqref{eq:intro-slightly-general} where the polynomial map $g$ has been replaced by $g'$, which has strictly lower step $\le s-1$.
Then we can employ the induction hypothesis and deduce the step $s$ case with an exponent $c_{n,s}:= c_{n,s-1}/O_{s,n}(1)$.

\subsection{Parametrization of ideals}\label{subsec:parametrization_of_ideals}

Green--Tao's computations mentioned above are already highly intricate,
but in retrospect, two algebraic facts were keeping them from becoming even more so:
all ideals of $\Z $ are principal;
every non-zero ideal $\ida $ of $\Z $ has a canonical basis $\ida \cong \Z $
(namely the unique positive generator).

While it is true that every non-zero ideal $\ida $ of $\OK $ is isomorphic to 
$\Z ^n $ as an abelian group, usually there is no canonical choice of a basis $\ida \cong \Z^n$.
Since many of the additive combinatorics quantities in the abelian group $\Z^n$ are basis-dependent,
taking random bases can cause a problem 
when we have to deal with unboundedly many ideals.

In \S \ref{appendix:norm-length} we recall that Geometry of Numbers provides a natural way of furnishing all non-zero ideals with a basis depending only on a few initial choices.
We see that these bases enjoys a uniformity property which we call 
the {\em norm-length compatibility}.
As long as the choice of bases is concerned, this property makes each ideal $\ida $ look almost as good as the sublattice $\lfloor \Nrm (\ida )^{1/n} \rfloor \cdot \OK $.

\subsection{Conclusion}
As was already explained, 
from the fact that the correlation \eqref{eq:intro-slightly-general} is small,
{\em Gowers Inverse Theorem} implies that $\lnorm{\Lambda _K - \LambdaSiegelIntro }_{U^{s+1}}$ is small.
Combined with the smallness of $\lnorm{\LambdaSiegelIntro - \LambdaCramerIntro }_{U^{s+1}}$ proven separately by Weil bounds, 
we conclude that $\lnorm{\Lambda _K-\LambdaCramerIntro }_{U^{s+1}}$ is small,
and therefore by the {\em von Neumann theorem} the sum 
\begin{align}
    \sum _{x\in \OKN }\prod _{i=1}^t \Lambda _K(\psi _i(x))
\end{align}
is close to the one where $\Lambda _K$ is replaced by $\LambdaCramerIntro $, which is elementarily computable.
This is how we will prove Theorem \ref{thm:GTZ-number-field}.

\section{Norm-length compatible bases}\label{appendix:norm-length}

We sometimes want to choose bases $\ida \cong \Zn $ of unboundedly many ideals $\ida $ in some uniform manner, in order to get a handle on quantities which depend on the choice of bases.
The notion of norm-length compatible bases and their existence serve this purpose. 

Given a ring homomorphism $\sigma \colon K\inj \C $, use the same symbol to denote its $\R $-linear extension $\sigma \colon \KR \to \C $.
We define the {\em canonical norm} $\lnorm{-}_{\can }$ on $\KR $ 
by 
\begin{align}
    \lnorm{x}_{\can }:= \max\limits _{\sigma \colon K\inj \C } |\sigma (x)|.
\end{align}

It is indeed a {\em norm} in the $\R$-linear algebra sense (= ``{\it length}'') because $\sigma $'s give rise to an $\R$-linear 
(ring) isomorphism $\KR \cong \R ^{r_1}\times \C^{r_2}$ 
where $r_1$ is the number of ring homomorphisms $K\inj \R $ and $r_2$ 
is the number of complex conjugate pairs of homomorphisms $K\inj \C $ whose image is not contained in $\R $.

Let us recall that a subset $\calD \subset (\KR)\baci $ is said to be {\em norm-length compatible}
(terminology from \cite{KMMSY})
if there is $C>1$ such that for all $\sigma \colon K\inj \C $ and $x\in \calD $, we have 
\begin{align}
    \frac{1}{C} \Nrm (x)^{1/n}
    < |\sigma (x)|
    < C \Nrm (x)^{1/n}
\end{align}
(the existence of an upper bound of this form implies the existence of a lower bound by the formula $\prod _{\sigma }|\sigma (x)|=\Nrm (x)$ and vice versa).

It is known that there are $\OK\baci $-fundamental domains $\calD \subset (\KR )\baci $ which are norm-length compatible (see e.g.\ \cite[\S 4.3]{KMMSY}).

Now we want to introduce the notion of norm-length compatible {\it bases} of ideals.

\begin{proposition}[existence of norm-length compatible bases]
    \label{prop:NLC-bases-exist}
    There is a positive number $C_K>1$ such that the following holds:
    for any non-zero fractional ideal $\ideala \subset K$, there is a 
    $\Z$-basis $\bm x = \{ x_1,\dots ,x_n\} $ of $\ideala $ 
    such that if we write $\lnorm{-}_{\bm x,\infty }$ for the associated $l^\infty $ norm, then we have the following relation of norms on $\ideala $
    \begin{align}\label{eq:def-of-norm-length-compatible-basis}
        \frac{1}{C_K}\Nrm (\ideala )^{1/n}\lnorm{-}_{\bm x,\infty }
        \le 
        \lnorm{-}_{\can }
        \le 
        {C_K}\Nrm (\ideala )^{1/n}\lnorm{-}_{\bm x,\infty }
        .
    \end{align}
\end{proposition}
\begin{proof}
    Choose a representative $\ideala _\lambda $ for each ideal class $\lambda \in \Cl (K)$.
    Choose a $\Z$-basis $x^\lambda _1,\dots ,x^\lambda _n$ of $\ideala_\lambda $ and let $\lnorm{-}_{\lambda ,\infty }$ be the associated $l^\infty $ norm.
    Since any two norms on a finite dimensional real vector space are equivalent,
    there is a $C>1$ such that for all $\lambda $ and $x\in \ideala_\lambda $, we have 
    \begin{align}\label{eq:equivalent-norms-on-ideala-lambda}
        \frac{1}{C}\Nrm (\ideala_\lambda )^{1/n}\lnorm{x}_{\lambda ,\infty } 
        \le  \lnorm{x}_{\can }
        \le  
        {C}\Nrm (\ideala_\lambda )^{1/n}\lnorm{x}_{\lambda ,\infty }  
        .
    \end{align}
    
    Choose also a norm-length compatible $\OK\baci $-domain $\calD \subset (\KR)\baci $. 
    Then there is a $C'>1$ such that for all $\sigma \colon K\inj \C $ and $y\in \calD $, we have 
    \begin{align}\label{eq:by-def-of-norm-length-compatibility}
        \frac{1}{C'}\Nrm (y)^{1/n}
        \le 
        |\sigma (y)|
        \le 
        {C'}\Nrm (y)^{1/n} 
        .
    \end{align}

    Now let $\ideala $ be an arbitrary non-zero fractional ideal.
    There are unique $\lambda \in \Cl (K)$ and $y\in \calD \cap K\baci $ such that $\ideala = y\ideala_\lambda $.
    Note that $\Nrm (\ideala)= \Nrm (y)\Nrm (\ideala_\lambda )$.
    Consider the $\Z$-basis $yx^\lambda _1,\dots ,yx^\lambda _n$ of $\ideala $.
    The associated $l^\infty $ norm $\lnorm{-}_{\ideala ,\infty }$ on $\ideala $ is identical to the norm induced by $\lnorm{-}_{\lambda ,\infty }$ via the multiplication-by-$y$ isomorphism $\ideala_\lambda \isoto \ideala $.

    For any $x\in \ideala_\lambda $, by multiplying \eqref{eq:equivalent-norms-on-ideala-lambda}
    and \eqref{eq:by-def-of-norm-length-compatibility} we get
    \begin{align}
        \frac{1}{CC'}\Nrm (\ideala )^{1/n}\lnorm{yx}_{\ideala ,\infty }
        \le 
        \lnorm{yx}_{\can }
        \le 
        CC' \Nrm (\ideala )^{1/n} \lnorm{yx}_{\ideala, \infty } .
    \end{align}
    Since every element of $\ideala $ has the form $yx$ with $x\in \ideala_\lambda $, 
    the coefficient $C_K:= CC'$ 
    and the bases $yx^\lambda _i$ ($1\le i\le n$)
    do the job.
\end{proof}

\begin{definition}[norm-length compatible bases]
    \label{def:NLC-bases}
    Choose a $C_K>1$ as in Proposition \ref{prop:NLC-bases-exist}
    and fix it throughout the paper.
    
    A $\Z $-basis $x_1,\dots ,x_n$ of a non-zero fractional ideal 
    $\ideala \subset K$
    is said to be {\em norm-length compatible} ({\em with coefficient $C_K$} if one wants to be precise)
    if the condition \eqref{eq:def-of-norm-length-compatible-basis} in Proposition \ref{prop:NLC-bases-exist}
    is satisfied.

    By Proposition \ref{prop:NLC-bases-exist}, every non-zero fractional ideal admits norm-length compatible bases.
\end{definition}

\begin{proposition}[change of basis matrices]
    \label{prop:matrix-between-norm-length-compatibles}
    Let $\ideala _1\subset \ideala_2 $ be an inclusion of two non-zero fractional ideals.
    Let $x^1_i$, $x^2_i$ ($1\le i\le n$) be a norm-length compatible 
    basis 
    (with coefficient $C_K$) 
    on each.
    Then the matrix representing the inclusion map with respect to these bases 
    has entries with absolute value $\le C_K^2 \Nrm (\ideala _1\ideala _2\inv )^{1/n}$.

    In particular, supposing $C_K$ has already been chosen, any two norm-length compatible bases of a given fractional ideal 
    are connected by matrices with entries of absolute value $\le O_K(1)$.
\end{proposition}
\begin{proof}
    Let $\lnorm{-}_{j,\infty }$ be the $l^\infty $ norm associated with the basis $x^j_1,\dots ,x^j_n$ of $\ideala _j$.
    By the definition of norm-length compatibility, we have for each $i$
    \begin{align}
        \lnorm{x^1_i}_{2,\infty }
        &<
        C_K \Nrm (\ideala _2)^{-1/n}\lnorm{x^1_i}_{\can }
        \\ 
        &<
        C_K \Nrm (\ideala _2)^{-1/n}\cdot C_K\Nrm (\ideala_1)^{1/n}\lnorm{x^1_i}_{1,\infty }
        \\ 
        &=
        C_K^2 \Nrm (\ideala_1\ideala_2\inv )^{1/n},
    \end{align}
    as is claimed.
\end{proof}

\begin{proposition}[normalized norm function]
    \label{prop:norm-function-coefficients}
    Let $\ida $ be a non-zero fractional ideal and $\iota \colon \Zn \cong \ida $ be a norm-length compatible basis.
    Then the function 
    \begin{align}
        \Zn &\to \Z 
        \\ 
        x&\mapsto N_{K/\Q} (\iota (x)) / \Nrm (\ida )
    \end{align}
    is represented by a homogeneous polynomial of degree $n$
    with coefficients of absolute value $\le O_K(1)$ (a bound independent of $\ida $).
\end{proposition}
\begin{proof}
    Choose $\calD \subset (\KR )\baci $ and $(\ida _\lambda )_{\lambda\in \Cl (K)}$ as in the proof of Proposition \ref{prop:NLC-bases-exist}.
    Take the unique $y\in \calD \cap K\baci $ and $\lambda $ such that 
        $\ida = y\ida _\lambda $.
    Since the matrix representing the multiplication-by-$y$ isomorphism 
    $\ida _\lambda \isoto \ida $ has entries of size $\le O_K(1)$,
    and since the functions $N_{K/\Q }(-)/\Nrm (\ida )$ and $N_{K/\Q }(-)/\Nrm (\ida _\lambda )$ are compatible via this isomophism at least up to sign,
    the problem is reduced to the case $\ida = \ida _\lambda $ for some $\lambda $.
    But there are only finitely many $\lambda \in \Cl (K)$.
\end{proof}

Let us take this opportunity to mention a classical result about the count of ideals:

\begin{proposition}[{density of ideals, see e.g.\ \cite[p.~210]{Cassels-Frohlich}}]
    \label{prop:density-of-ideals}
    Let $\lambda\in \Cl (K)$ and $N\ge 3$.
    We have 
    \begin{align}
        \sum _{\substack{
            \ida\in\IdealsK \\ \Nrm (\ida )\le N \\ [\ida ]=\lambda \te{ in }\Cl (K) }}
        1 
        = 
        \frac{\residue _{s=1} (\zeta _K(s))}{\# \Cl (K)}
        N 
        + O_K(N^{1-\frac 1n})
        .
    \end{align}
\end{proposition}

\newcommand{\complexity}{M} 

\section{Recollection of nilsequences}\label{sec:nilsequences}

Nilsequences are sequences of complex numbers obtained in a certain geometric way.
We recall a minimum amount of definitions 
(and fix conventions thereof)
that are needed to make sense of the discussions that follow.
For details and references, see \cite[\S\S 1, 2]{PolynomialOrbits}.
Those who already know this concept should skip this section.
Advanced results that we will eventually need in the most technical part of this article are listed in Appendix \ref{appendix:nilsequences}.

We start out with the notions of nilmanifolds and polynomial maps into them.

\begin{definition}[filtered Lie groups and nilmanifolds]
    \label{def:filtered_groups}
    Let $k\ge 0$ be an integer.
    A {\em filtered group of degree $\le k$}
    is an abstract group $G$ equipped with a descending sequence of normal subgroups
    \begin{align}
        G= G_1\supset \dots \supset G_k\supset G_{k+1}=\{ 1 \}
    \end{align}
    such that $[G_i,G_j]\subset G_{i+j}$ for all $i,j \ge 0 $.

    A {\em filtered Lie group of degree $\le k$} is a connected, simply connected Lie group $G$ equipped with a descending sequence of normal subgroups as above, where $G_i$ are also required to be closed and connected.
    In this case $G$ is a nilpotent Lie group of {\em step} $\le k$: if define its {\em lower central series} $G_{(i)}$ ($i\ge 1$)
    by $G_{(1)}:= G$ and $G_{(i+1)}:= \overline{ [G,G_{(i)}]}$,
    then we have $G_{(k+1)}=\{ 1 \}$.
    This follows directly from the condition $[G_i,G_j]\subset G_{i+j}$.

    A {\em nilmanifold of degree $\le k$}
    is the data of a filtered Lie group $G$ of degree $\le k$
    and a cocompact discrete subgroup $\Gamma \subset G$.
    Usually we simply say ``a nilmanifold $G/\Gamma $ of degree $\le k$'' when it is not strictly necessary to spell out the filtered Lie group.

\end{definition}

\begin{definition}[\Malcev basis]
    A {\em \Malcev basis $X_1,\dots ,X_{m}$ adapted to} a nilmanifold $G/\Gamma $ of degree $\le k$
    is a basis of the $\R$-vector space $\Lie (G)$ satisfying:
    \begin{itemize}
        \item The map 
        \begin{align}
           \psi_{X_1,\dots ,X_m}
           \colon 
           \R ^m &\to G 
            \\ 
            (t_1,\dots ,t_m)&\mapsto \exp (t_1X_1)\dots \exp (t_mX_m)
        \end{align}
        is a bijection mapping $\Z^m$ exactly onto $\Gamma $.
        \item For each $0\le j\le m$, the $\R$-span $\mathfrak h_j$ of the vectors $X_{j+1},\dots ,X_m$ 
        is an ideal of the Lie algebra $\Lie (G)$. (Note that $\dim (\mathfrak h_j)=m-j $.)
        \item For each $0\le i\le k$, the Lie algebra of $G_i$ is $\mathfrak h_{m-\dim (G_i)}$.
    \end{itemize}
\end{definition}

\Mal bases are known to exist for every nilmanifold \cite[p.~478]{PolynomialOrbits}.

\begin{definition}[metric on $G/\Gamma $ associated with a \Mal basis]
    \label{def:metric}
    Let $G/\Gamma $ be a nilmanifold equipped with a \Mal basis.
    Give $\Lie (G)$ the positive definite inner product where the given basis is an orthonormal basis.
    Translate it to the entire $G$ by the action from the right, which makes $G/\Gamma $ a Riemannian manifold and in particular a metric space.
\end{definition}
    In our main references such as \cite[Definition 2.2, also the footnote]{PolynomialOrbits}, they use a somewhat different metric based on the $l^\infty $-norm.
Since we have the inequality $\lnorm x_\infty \le \lnorm x_2 \le \sqrt{m}\cdot \lnorm x_\infty $ for $x\in \R ^m$
and as we will pretty quickly have to deal with much larger multiplicative losses than $\dim (G)^{1/2}$,
the choice of the metric here turns out insignificant.

\begin{definition}[complexity of a \Mal basis and subgroups]
    Let $\complexity >0$ be a positive number.

    A \Mal basis $X_1,\dots ,X_m$ of a nilmanifold $G/\Gamma $ is said to have {\em complexity $\le \complexity $}
    (or to be {\em $\complexity $-rational})
    if the structural constants $c_{ijk}$ defined by 
    \begin{align}
        [X_i,X_j] = \sum _{\ell =1}^m c_{ij\ell } X_\ell 
    \end{align}
    are all rational numbers of height $\le \complexity $.
    (I.e., written as the quotient of integers of size $\le \complexity $.)
    If the choice of a \Mal basis is implicitly assumed, we simply say 
    that the nilmanifold $G/\Gamma $ has {\em complexity} $\le \complexity $.

    A closed subgroup $G'\subset G$ is said to have {\em complexity $\le \complexity $}
    (or to be {\em $\complexity $-rational})
    with respect to a \Mal basis $X_1,\dots ,X_m$
    if $\Lie (G')$ has a basis where each member can be written as a linear 
    combination of $X_1,\dots ,X_m$ with rational coefficients of height $\le \complexity $.
\end{definition}

\begin{definition}[horizontal characters]
    \label{def:horizontal-char}
    A {\em horizontal character} $\eta $ on a nilmanifold $G/\Gamma $
    is a continuous group homomorphism $\eta \colon G\to \RZ $
    which annihilates $\Gamma $.

    If we are also given a \Mal basis $X_1,\dots ,X_m$ adapted to $G/\Gamma $,
    there is a unique $\undl a=(a_1,\dots ,a_m)\in \Z^m$ such that for all $(t_1,\dots ,t_m)\in \R^m$:
    \begin{align}
        \eta (\psi _{X_1,\dots ,X_m} (t_1,\dots ,t_m))
        = a_1t_1+\dots +a_m t_m  \mod \Z.
    \end{align}
    The {\em size} 
    $|\eta |$ of $\eta $ is defined by 
    \begin{align}
        |\eta |:= \max _{1\le i\le m} |a_i| .
    \end{align}
\end{definition}

\begin{definition}[polynomial maps]
    Let $H$, $G$ be filtered groups in the sense of Definition \ref{def:filtered_groups}.
    For a map of sets $g\colon H\to G$ and $h\in H$, let $\partial _h g\colon H\to G$ 
    be the map defined by $(\partial _h g )(x)= g(xh)g(x)\inv $ for $x\in G$.
    
    We say $g$ is a {\em polynomial map} if for any $i_1,\dots ,i_m\ge 0$ and $h_1 \in H_{i_1}$, \dots, $h_m\in H_{i_m}$,
    the map $\partial _{h_1}\dots \partial _{h_m}g$ takes values in $G_{i_1+\dots +i_m}$.
\end{definition}
It is known that the composite of two polynomial maps is a polynomial map \cite[Theorems 1.6.9, 1.6.10]{TaoFourier}.
We usually give abelian groups $G$ the unique degree-$1$ filtration $G=G_1\supset G_2=\{ 0\}$.
A map $g\colon H\to G$ between abelian groups with this filtration is a polynomial map 
if (and actually only if) it is an affine-linear map (= a group homomorphism followed by a translate).

Consider the situation $H=\Z ^n$ equipped with the degree-$1$ filtration and $G$ being a filtered, conneced, simply connected Lie group where $\Lie (G)$ is equipped with a basis $X_1,\dots ,X_m$ satisfying the conditions for a \Mal basis except the condition $\psi _{X_1,\dots ,X_m}(\Z ^m)=\Gamma $ (because no $\Gamma $ has been given yet).
In this case
there is an explicit description of a polynomial map in terms of coordinates $\psi = \psi _{X_1,\dots ,X_m}\colon \R^m\cong G$.
To explain it, for a vector $\undl x=(x_1,\dots ,x_n)\in \Z^n$ and a multi-index $\undl i=(i_1,\dots ,i_n)\in \N ^n$,
let us use the power notation
\begin{align}
    \undl x ^{\undli }:= x_1^{i_1}\cdot \dots \cdot x_n^{i_n} \in \Z .
\end{align} 
We also use the symbol $|\undli |:= i_1+\dots +i_n$.
It is known \cite[Lemma 6.7]{PolynomialOrbits} that a map $g\colon \Z^n \to G$ is a polynomial map if and only if 
the composite $\psi\inv \circ g\colon \Z^n\xto g G \xto[\cong ]{\psi\inv }\R^m$
is written as: 
\begin{align}
    \psi\inv \circ g (\undl x)
    = 
    \sum _{\undli \in \N^n }a_{\undli }\undl x ^{\undli }
\end{align}
for some coefficients $a_{\undli }\in \Lie (G_{|\undli|})=\{ 0 \}\times \R^{\dim (G_{\undli })}$.
Simply put, polynomial maps are maps written as polynomials 
where the coefficients in higher degrees are 
required to be in smaller subgroups.

\begin{definition}[smoothness norm]\label{def:smoothness-norm}
    For a real number $\alpha \in \R $, write 
    \begin{align}
        \lnorm \alpha _{\RZ } := \min _{n\in \Z } |\alpha -n| .
    \end{align}

    Let $g\colon \Z^n \to \R$ be a real-valued polynomial map of degree $s$,
    written as 
    \begin{align}
        g(\undl x)= \sum _{\substack{
            \undli \\  |\undli |\le s}}
            a_{\undli }\undl x^{\undli }.
    \end{align}
    Let $N\ge 1$ be a natural number.
    Its {\em smoothness norm} $\lnorm g _{C^\infty [-N,N]^n}$ is defined by 
    \begin{align}\label{eq:def-of-C-infty-norm}
        \lnorm g_{C^\infty [-N,N]^n}
        :=
        \max _{\undli \neq 0} \left( N^{\undli }\lnorm{a_{\undli }}_{\RZ } \right) .
    \end{align}
\end{definition}
In fact, this should be written as $\lnorm{g}_{C_*^\infty [-N,N]^n}$ (with a star $*$) if we follow \cite[p.~4]{PolynomialOrbitsErratum}, which is different from some pre-existing $C^\infty [-N,N]^n$ norm.
However, as in \cite[Lemma 2.1]{PolynomialOrbitsErratum} there is some integer $q=O(s!)$ such that $\lnorm{qg}_{C^\infty [-N,N]^n} \le q \lnorm{g}_{C_*^\infty [-N,N]^n} $
and vice versa, and this makes all the results we use valid regardless of which norm we choose.
As \eqref{eq:def-of-C-infty-norm} is more convenient for the computation we will perform, we use it in this paper and simply write it as $\lnorm{g}_{C^\infty [-N,N]^n}$ without a star.

\begin{definition}[nilsequences]
    A {\em nilsequence of degree $\le k$} is a chain of maps 
    \begin{align}
        \Z ^n \xto g G \surj G/\Gamma \xto F \C ,
    \end{align}
    where $G/\Gamma $ is a nilmanifold of degree $\le k$,
    $g$ is a polynomial map, 
    $G\surj G/\Gamma $ is the canonical map 
    and $F$ is a \Lip function.

    It is said to be {\em $1$-bounded} if the values of $F$ always have absolute value $\le 1$,
    to have {\em \Lip constant} $\le L$ if $F$ does,
    and to have {\em complexity} $\le M$ if $G/\Gamma $ is equipped with such a \Mal basis.
    
\end{definition}

Nilsequences are used as test functions for the smallness of the Gowers norm. 
Namely, 
the {\em Gowers Inverse Theorem} \ref{thm:GowersInverse}
roughly says that 
when one knows that a 1-bounded function $f\colon [-N,N]^n\to \C $ has small correlation 
\begin{align}\label{eq:sec_nil:correlation}
    \left| \Expec _{x\in [-N,N]^n} f(x) F(g(x)) \right|
\end{align}
with every $1$-bounded $s$-step nilsequences $F(g(-))\colon \Z^n\to \C $ of controlled \Lip constant and complexity,
then it has small $U^{s+1}$-norm $\lnorm f _{U^{s+1}[-N,N]^n}$.
For this reason, we will be studying the correlation \eqref{eq:sec_nil:correlation} with $f=\Lambda _K - \LambdaSiegel $ in the most part of this paper.

Moreover, there is {\em (quantitative) equidistribution theory} of polynomial sequences (such as Theorem \ref{thm:Leibman}) to the effect that the image of $[-N,N]^n \xto g G \surj G/\Gamma $ is equidistributed on some submanifold within a certain error when $N$ is large.
This will be used as one of the main tools in Sections \ref{sec:Type-I}, \ref{sec:Type-II} for the analysis of the correlation \eqref{eq:sec_nil:correlation}. 

\section{The \Cramer model and some $\Z $-linear algebra}\label{sec:Cramer-model}

Recall from \S \ref{sec:notation} the quantity
$\exp ((\log N)^{1/100})<Q<\exp ( (\log N) ^{1/3})$ and 
symbol $P(Q):= \prod _{0<p<Q} p $.
Recall that the von Mangoldt function $\Lambda _K\colon \IdealsK \to \R_{\ge 0}$
is defined by 
\begin{align}
    \Lambda _K (\idc )  = 
    \begin{cases}
        \log \Nrm (\idp ) & \te{ if }\idc = \idp ^m \te{ ($\idp $ a non-zero prime ideal, $m\ge 1$) }
        \\ 
        0 &\te{ if not.}    
    \end{cases}
\end{align}
It is convenient to also define $\Lambda _K ((0)):= 0$.

For a non-zero fractional ideal $\ideala $ of $K$, 
we define $\Lambda ^{\ideala }_K\colon \ideala \to \R $ by 
\begin{align}
    \Lambda ^{\ideala }_K(x)= 
    \Lambda _K(x\ideala\inv )
    .
\end{align}

A primitive form of Mitsui's Prime Number Theorem (\cite{Mitsui}, recalled in the proof of Proposition \ref{prop:essentially-Mitsui})
says that for a convex body 
$\Omega \subset \KR $, we have 
\begin{align}\label{eq:Mitsui-primitive-form}
    \sum _{x\in \Omega \cap \ideala } \Lambda ^{\ideala }_K (x)
    &\fallingdotseq 
    \frac{1}{\residue _{s=1}\limits (\zeta _K(s))} \sum _{x\in \Omega\cap \ideala } 1 
    .
\end{align}
This is why we put the coefficient $\frac{1}{\residue _{s=1}\limits (\zeta _K(s))}$ in the next definition.

\begin{definition}
    We define the \Cramer model $\Lambda _{\mrm{Cram{\acute e}r},Q}$ of the von Mangold function $\Lambda _K$ to be the function $\IdealsK \cup \{ (0) \} \to \R $ defined by 
    \begin{align}
        \Lambda _{\mCramer ,Q} (\idc ) = \begin{cases}
            0 & \te{ if } \idealc + (P(Q)) \neq \OK \\
       \disp \frac{1}{\residue_{s=1}\limits\zeta _K(s)}\frac{P(Q)^n}{\totient (P(Q))} & \te{ if } \idealc + (P(Q)) = \OK .
    \end{cases}
\end{align}
For a non-zero fractional ideal $\ideala $, we define $\Lambda ^{\ideala }_{\mCramer ,Q}\colon \ideala \to \R $ by 
\begin{align}
    \Lambda ^{\ideala }_{\mCramer ,Q} (x)
    &= 
    \Lambda _{\mCramer ,Q} (x\ideala\inv ) 
    \\
    &= 
    \begin{cases}
        \displaystyle \frac{1}{\residue\limits (\zeta _K(s))}\frac{P(Q)^n}{\totient (P(Q))} & \te{ if $x$ coprime to $P(Q)$ in $\ida $}
        \\ 
        0 & \te{ if not}
         .
    \end{cases}
\end{align}
(Recall that $x\in \ida $ being coprime to $P(Q)$ in $\ida $ means $x\OK + P(Q)\ida = \ida $.)
\end{definition}

The usefulness of the \Cramer model is that its statistical behavior is quite predictable. 
In order to make the 
next proposition useful enough later on, 
we have to formulate it using {\em subsets of \Lip boundary} (Appendix \ref{appendix:lattice-points}), a more general notion than convex subsets.
This is mainly because the subset obtained by bounding the norm
$\Omega _{\le a}:= \{ x\in \KR \mid |\Nrm (x)|\le a \}$
is not convex unless $K$ is an imaginary quadratic field or $K=\Q $.  
Convex sets in $[-N,N]^d$ are known to have $\mLip (d,1,O_d(N))$ boundaries. 

\begin{proposition}[large scale behavior]
    \label{prop:sieve-theory}
    Let $\psi _1,\dots ,\psi _t \colon \Z ^d \to \ideala $ be $\Z $-affine-linear maps.
    Write $\dot\psi _i$ for their linear parts.
    Choose a norm-length compatible basis
    $\iota\colon \ideala \isoto \Z ^n$
    and assume:
    \begin{itemize}
        \item all $\dot\psi _i$ have finite cokernels,
        \item $\ker (\dot\psi _i)$ does not contain $\ker (\dot\psi _j)$ for all $i\neq j$,
        \item the coefficients of $\iota \circ \dot\psi _i\colon \Z ^d\to \Z ^n$ have absolute values bounded from above by some $1<L <\exp ((\log N)^{1/101})$.
    \end{itemize}
    
    Let $N>1$ be an integer and $\Omega\subset [\pm N]_\R ^d$ be a closed set with boundary of \Lip class $\mLip (d,M_1,M_2N)$ {\em (Definition \ref{def:Lip_class})} for some $M_1,M_2\ge 1$.
    Assume that $N$ is large enough depending on $M_1,M_2$.
\begin{enumerate}
    \item 
    We have 
    \begin{multline}\label{eq:prop-sieve}
        \sum _{x\in \Omega\cap \Z ^d}
        \prod _{i=1}^t
        \LambdaaCramer 
        (
            \psi _i (x)
            )
            =
            \frac{\vol (\Omega )}{\residuezeta ^t} \left(
                \prod _{0<p<Q} \beta _p (\psi _1,\dots ,\psi _t)
                \right)
                \\ 
                +
                O_{t,d,K}(N^{d}\exp (-\frac 15 (\log N)^{1/2})  ) ,
            \end{multline} 
            where $\beta _p(\psi _1,\dots ,\psi _t)$ is the following coefficient depending on the maps $\psi _i$ mod $p$:
            \begin{align}
                \beta _p (\psi _1,\dots ,\psi _t):= \left( \frac{p^n}{\totient (p)}\right) ^t 
                \Expec _{x\in (\Z /p\Z )^d }
                \left[
                    \prod _{i=1}^t \prod _{\idealp \divides p}
                    1[\psi _i(x)\neq 0 \te{ in }\ida /\idealp\ideala  ]
                    \right]
                    .
                \end{align} 
                \item Suppose $M=\Z^d $ carries an $\OK $-module structure and $\psi _i$ are affine-$\OK $-linear.  
                For each prime ideal $\idp $, define 
                \begin{align}
                    \beta _\idp (\psi _1,\dots ,\psi _t)
                    := \left( \frac{\Nrm (\idp )}{\totient (\idp )}\right) ^t 
                    \Expec _{x\in M /\idp M }
                    \left[
                        \prod _{i=1}^t 
                        1[\psi _i(x)\neq 0 \te{ in }\ida/\idealp\ideala  ]
                        \right]
                        .
                    \end{align}
                Then
                we have $\beta _p(\psi _1,\dots ,\psi _t)= \prod _{\idp\divides p}\beta _\idp (\psi _1,\dots ,\psi _t)$.
            \end{enumerate}
            \end{proposition}
\begin{proof}
    This proof is kind of routine, and the only purpose of writing it down is to explicate the dependence of the error term and the threshold for $N$ on the parameters.

    The decomposition (2) of $\beta _p(\psi _1,\dots ,\psi _t)$ is a consequence of Chinese Remainder Theorem for torsion $\OK $-module.
    To prove the main assertion (1), we shall first estimate the sum 
    \begin{align}\label{eq:sieve-lemma-can-be-applied}
        \sum _{x\in \Omega \cap \Z ^d}
        \prod _{i=1}^t
        \prod _{p<Q}
        \prod _{\substack{
             \idealp \divides p  
                        } 
                }
            1[\psi _i (x)\not\in \idealp\ideala ] 
        \end{align} 
    and then multiply it by the coefficient $\frac{1}{\residuezeta ^t}\prod _{0< p<Q} \left( \frac{p^n}{\totient (p)} \right) ^t$ to make it equal the left-hand side of \eqref{eq:prop-sieve}.

    For each $m\in\N $, set
    \begin{align}
        \Sigma _m&:= \left\{ x\in \Omega \cap \Z ^d \ \middle|\ \prod _{i=1}^t \Nrm (\psi _i(x)\ideala\inv ) =m  \right\} ,
        \\
         a_m&:= |\Sigma _m|.
    \end{align}
    Note that the condition $p\divides \prod _{i=1}^t \Nrm (\psi _i(x)\ideala\inv )$ for a prime number $p$ is equivalent to the condition that $\psi _i(x)\in \idealp \ideala $ for some prime ideal $\idealp $ over $p$ and some index $i$.
    Thus a point $x\in \Omega \cap \Z ^d$ contributes to the sum 
    \eqref{eq:sieve-lemma-can-be-applied} if and only if 
    $x\in \Sigma _m$ for some $m\in N$ that is coprime to all $p<Q$.
    
    It follows that the sum \eqref{eq:sieve-lemma-can-be-applied} equals 
    \begin{align}
        \sum _{\substack{
            m\in\N \\ (m,P(Q))=1
        } }
        a_m .
    \end{align}
    We want to evaluate it using Fundamental Lemma \ref{lem:sieve-theory} of Sieve Theory.

    For each prime number $p$, write 
    \begin{align}\label{eq:def-of-g(p)}
        g(p)&:= \Expec _{x\in (\Z /p\Z )^d } 1[{\psi _i (x)\in \idealp\ideala \te{ for some }i\te{ and }\idealp \te{ over } p }] 
        \\ 
        & = \Expec _{x\in (\Z /p\Z )^d }
        \left[
        1- \prod _{i=1}^t \prod _{\idealp \divides p} 
            1[\psi _i (x)\not\in \idealp \ideala ]
        \right]
        .
    \end{align}

    Suppose first that $g(p)=1$ for some $p<Q$.
    Then $\beta _p = \left( \frac{p^n}{\totient (p)}\right) ^t (1-g(p)) =0$ for this $p$.
    Meanwhile, the function $\prod _{i=1}^t \Lambda^{\ideala }_{\mCramer ,Q} (\psi _i(x))$ is identically zero under the current supposition.
    Therefore Proposition \ref{prop:sieve-theory} is true without the need to use sieve method.

    Next, let us assume $g(p)<1$ for all $p<Q$.
    To invoke Lemma \ref{lem:sieve-theory}, let us first verify the hypothesis \eqref{eq:sieve-hyp-bound} on the size of the product $\prod _{w\le p<Q}(1-g(p))\inv $.
    For an index $i$ and a prime number $p$, write
    \begin{align}
        \psi _{i,p},\ \dot\psi _{i,p} \colon (\Z /p\Z )^d \to \ideala /p\ideala 
    \end{align}
    for the 
    (affine-)$\F_p$-linear 
    maps induced by $\psi _i$.
    When 
    a prime ideal $\idealp $ over $p$ is moreover given, write 
    \begin{align}
        \psi _{i,\idealp},\ \dot\psi _{i,\idealp } \colon (\Z /p\Z )^d \to \ideala / \idealp\ideala 
    \end{align}
    for the induced (affine-)$\F_p$-linear maps.
    \begin{claim}\label{claim:kernels-dont-contain-each-other}
        There exists a natural number $A=O_{d,t,K }(L^{O_{n,t}(1)})$ such that 
        the following hold for all prime numbers $p$ coprime to $A$:
        \begin{enumerate}
            \item $\psi _{i,p}$ is surjective for all $i$;
            \item for two indices $i,j$ and two prime ideals $\idealp ,\idealq $ of residue characteristic $p$,
            whenever the pairs $(i,\idealp )$, $(j,\idealq )$ are distinct we have 
            $\ker (\dot\psi _{i,\idealp })\not\subset \ker (\dot\psi _{j,\idealq })$.
        \end{enumerate}
    \end{claim}
    \begin{proof}[Proof of Claim]
        By undergraduate linear algebra ($n\times n$ minors for the $t$ maps $\dot\psi _i$), there is a natural number 
        \begin{align}
            A_1=O_{n,t}(L^{O_{n,t}(1)})
        \end{align}
        such that if $p\notdivide A_1$ the maps $\dot\psi _i $ mod $p$:
        $(\Z / p\Z )^d \to \ideala /p\ideala \overset{\iota}{\cong }(\Z /p\Z )^n $ are surjective.
        (So any multiple of $A_1$ will do for the first assertion.)

        By the assumption that $\ker (\dot\psi_i)$'s do not contain each other,
        one can find elements $x_{ij}\in \ker (\dot\psi_i)\setminus \ker (\dot\psi _j)$
        which moreover can be taken to have 
        coefficients of sizes $O_{d,n}(L^{O_n(1)}) $.
        Their images $\dot\psi _j (x_{ij})\in \ideala \nonzero \overset{\iota }{\cong} \Z ^n \nonzero $ have entries of sizes 
        $O_{d,n}(L^{O_n(1)})$. 
        Hence their norms $\Nrm (\dot\psi _j (x_{ij})\ideala \inv )$ have sizes $O_{d,n,K }(L^{O_n(1)}) $.
        Let $A_2$ to be the least common multiple of these norms: 
        \begin{align}
            A_2:= \oname{lcm}\{ \Nrm (\dot\psi _j (x_{ij})\ideala \inv ) \mid 1\le i,j\le t \}
            = O_{t,d,K } (L^{O_{n,t}(1)}) .
        \end{align}

        Let $A$ be the least common multiple of $A_2$ and $A_1$:
        \begin{align}
            A:= \oname{lcm}\{ A_1,A_2 \} .
        \end{align}

        Now let $(i,\idealp )$ and $(j,\idealq )$ be two distinct pairs
        such that $\idealp ,\idealq $ have the same residue characteristic $p\notdivide A$.
        We have to show $\ker (\dot\psi_{i,\idp})\not\subset \ker (\dot\psi_{j,\idq })$.

        First assume $i\neq j$. 
        In this case the residue class of $x_{ij}$ in $(\Z /p\Z )^d$
        certainly belongs to $\ker (\dot\psi _{i,\idealp} )$.
        On the other hand, 
        suppose $x_{ij}\in \ker (\dot\psi _{j,\idq })$, namely $\dot\psi _j (x_{ij})\in \idealq\ideala $.
        This would imply $\Nrm (\dot\psi_j(x_{ij})) \in \Nrm (\idq\ida )\Z $
        and hence 
        \begin{align}
            \frac{A_2}{p} = \frac{A_2}{\Nrm (\dot\psi_j(x_{ij})\ida\inv )}\cdot \frac{\Nrm (\dot\psi_j(x_{ij})) }{\Nrm (\ida )\Nrm (\idq )}\cdot \frac{\Nrm (\idq )}{p}
            \in \Z ,
        \end{align}
        contradicting the condition $p\notdivide A_2$.
        Thus we have shown $\ker (\dot\psi _{i,\idealp })\not\subset \ker (\dot\psi _{j,\idq })$.

        Second, assume $i=j$ and hence $\idealp \neq \idealq $.
        In this case, by Chinese Remainder Theorem we have a natural surjection 
        $\ideala /p\ideala \surj \ideala /\idealp \ideala \oplus \ideala /\idealq \ideala $.
        Since $\dot\psi _i$ mod $p$ is a surjection $(\Z /p\Z )^d\surj \ideala /p\ideala $ when $p\notdivide A_1$, we conclude that the map induced by $\dot\psi _i$
        \begin{align}
            (\Z /p\Z )^d \surj \ideala /\idealp \ideala \oplus \ideala /\idealq \ideala 
        \end{align}
        is surjective.
        In particular there is an element which maps to zero in $\ideala /\idealp\ideala $ and to a non-zero element in $\ideala /\idealq\ideala $.
        Therefore we conclude $\ker (\dot\psi _{i,\idealp })\not\subset \ker (\dot\psi _{i,\idealq })$.

        This completes the proof of Claim \ref{claim:kernels-dont-contain-each-other}.
    \end{proof}
    Resuming the proof of Proposition \ref{prop:sieve-theory},
    now we can estimate $g(p)$ for $p\notdivide A$ as follows. Note that by the inclusion-exclusion principle we have 
    \begin{align}
        g(p) &= \# \left(
            \bigcup _{\substack{
                (i,\idealp ) \\ 1\le i\le t \\ \idealp \divides p
            } }
            \psi _{i,\idealp }\inv (0) 
        \right)
        / p^d
        \\ 
        &=
        \left( 
        \left(\sum _{\substack{
                (i,\idealp ) \\ 1\le i\le t \\ \idealp \divides p
            } }
            \# \psi _{i,\idealp }\inv (0) 
            \right)
        - \left( \sum _{(i,\idealp )\neq (j,\idealq )}
        \#  (\psi_{i,\idealp }\inv (0)\cap \psi _{j,\idealq }\inv (0)) 
        \right)
        + \cdots  \right)
        /p^d ,
    \end{align}
    where the alternating sum has $O_{t,n}(1)$ members.
    Since the linear part $\dot\psi _{i,\idealp }$ is surjective,
    we know that the inverse image 
    $\psi _{i,\idealp }\inv (0)$ is a translate of $\ker (\dot\psi _{i,\idealp })$. 
    By Claim \ref{claim:kernels-dont-contain-each-other}, 
    the intersection $\ker (\dot\psi _{i,\idealp }) \cap \ker (\dot\psi _{j,\idq })$ or the intersection of more members has codimension $\ge 2$ in the $\F _p$-vector space $(\Z / p\Z )^d$.
    We conclude that the above expression can be written as  
    \begin{align}
        &
        \left(\sum _{\substack{
                (i,\idealp ) \\ 1\le i\le t \\ \idealp \divides p
            } }
            \# \ker (\dot\psi _{i,\idealp })
            \right)
        /p^d
        + O_{t,n}\left( \frac 1{p^2} \right)
        \\ 
        =
        &
        \left(
            \sum _{\substack{
                (i,\idealp ) \\ 1\le i\le t \\ \idealp \divides p
            } }
         \frac{1}{\Nrm (\idealp )}
        \right) 
        +O_{t,n}\left( \frac 1{p^2}\right)
        \\
        =
        &
        \left(
            \sum _{\substack{
                \idealp \divides p
            } }
         \frac t{\Nrm (\idealp )}
        \right) 
        +O_{t,n}\left( \frac 1{p^2}\right)
        .
    \end{align}
    It follows that for $p\notdivide A$ we have 
    \begin{align}\label{eq:1-g(p)_for_large_p}
        (1-g(p))\inv = (1+O_{t,n}(1/p^2))\prod _{\idealp\divides p} \left( 1-\frac{1}{\Nrm (\idealp )} \right) ^{-t} .
    \end{align}
    For $p\divides A$ we will simply use the crude inequality $g(p)\le 1-\frac{1}{p^{d}}$,
    which holds because the expectation $\Expec $ defining $g(p)$ is over the set $(\Z /p\Z )^d$ of cardinality $p^{d}$ and we are assuming $g(p)<1$ for all $p<Q$;
    equivalently 
    \begin{align}\label{eq:sieve-crude-inequality}
        (1-g(p))\inv \le p^{d}.
    \end{align}

    Rosen's Mertens' theorem for number fields (Theorem \ref{thm:Mertens-Rosen}) states 
    \begin{align}\label{eq:Mertens-Rosen}
        \prod _{\substack{ \idealp \\ \Nrm (\idealp )<x}}
        \left( 1-\frac{1}{\Nrm (\idealp )} \right) ^{-1}
        = e^\gamma \cdot \residue _{s=1}(\zeta _K(s)) \cdot \log x +O_K(1). 
    \end{align} 
    Formulas \eqref{eq:1-g(p)_for_large_p} \eqref{eq:sieve-crude-inequality} \eqref{eq:Mertens-Rosen},
    together with the fact that prime ideals $\idealp $ with norm $\ge p^2$ are negligibly few,\footnote{in the sense that the product $\prod _{\idealp } \left( 1-\frac{1}{\Nrm (\idealp )} \right)\inv $ over such $\idealp $'s converges, say.} 
    imply that 
    \begin{align}
        \prod _{w\le p<Q } (1-g(p)) \inv 
        &= 
        \prod _{w\le p<Q,\ p\divides A} (1-g(p)) \inv 
        \prod _{w\le p<Q,\ p\notdivide A} (1-g(p)) \inv 
        \\ 
        &\ll _{K,t} 
        \prod _{p\divides A }p^{d}
        \cdot \left( \frac{\log Q}{\log w} \right) ^{t } 
        \ll _{K,t}
        A^{d}  \left( \frac{\log Q}{\log w} \right) ^{t }   ,
    \end{align}
    verifying the hypothesis \eqref{eq:sieve-hyp-bound} with 
    \begin{align}\label{eq:M-and-kappa}
        M=O_{K,t}(A^{d})\te{ and }\kappa :=t. 
    \end{align}

    Now we want to establish a relation of type \eqref{eq:coefficient-X} in Lemma \ref{lem:sieve-theory}.
    Set $X:= \vol (\Omega )$.
    For a squre free natural number $f$,
    a point $x\in \Omega \cap \Z ^d$ is in $\dunion _{m\in\N } \Sigma _{fm} $ if and only if 
    $\prod _{i=1}^t \Nrm (\psi _i (x)\ideala \inv )$ is a multiple of $p$ for all prime factors $p\divides f$.
    This latter condition is equivalent to the condition that for all $p\divides f$, there is a prime ideal $\idealp \divides p $ and an index $i$ such that $\psi _i (x)\in \idealp \ideala $. Thus:
    \begin{align}
        \dunion _{m\in\N } \Sigma _{fm}
        = 
        \bigcap _{p\divides f} \{ x\in \Omega \cap \Z ^d\mid \psi _i (x)\in \idealp \ideala \te{ for some }i \te{ and }\idealp \divides p \} .
    \end{align} 
    Note that the condition in the bracket only depends on the mod $f$ class of $x\in \Z ^d$.

    Consider the complete set of representatives $[0,f-1]^d\subset \Z ^d$
    for the quotient $(\Z /f\Z )^d $. 
    Out of the $f^{d}$ residue classes in $(\Z /f\Z )^d$, there are exactly $f^{d} g(f)$ classes (any representatives of) which belong to $\dunion _{m\in\N } \Sigma _{fm}$; 
    this is true by the very definition of $g(p)$ in \eqref{eq:def-of-g(p)} and by Chinese Remainder Theorem for various $p\divides f$.
    Fill the interior of $\Omega \cap \Z ^d$ by as many disjoint translates of $[0,f-1]^{d}$ as possible.
    By Proposition \ref{cor:Widmer-rescaled} there are $\ll _{d,M_1,M_2}fN^{d-1}$ points of $\Omega \cap \Z^d$ not contained in any of these translates.
    Note that there are $\vol (\Omega ) + O_{d,M_1,M_2}(N^{d-1})$ points in $\Omega \cap \Z ^d$.
    We conclude:
    \begin{align}
        \sum _{m\in\N} a_{fm}
        = 
        \vol (\Omega ) g(f)
        + 
        O_{d,M_1,M_2}(fN^{d-1}),
    \end{align}
    which implies that 
    when we take $X=\vol (\Omega ) $
    the remainder $r_f$ in \eqref{eq:coefficient-X}
    is an $O_{d,M_1,M_2}(fN^{d-1})$.

    Now we can apply Fundamental Lemma \ref{lem:sieve-theory} of Sieve Theory. This gives
    \begin{multline}\label{eq:after-sieve-before-beta}
        \sum _{\substack{
            m\in\N \\ (m,P(Q))=1
        } }
        a_m 
        =
        \vol (\Omega )
        \left(
            \prod _{p<Q} (1-g(p))
            \right)
        (1+O(M^{10}e^{9t-\frac{\log D}{\log Q}}))
        \\ 
        +
        O_{d,M_1,M_2}\left( \sum _{
            \substack{1\le f\le D \\ f\divides P(Q)}
        } fN^{d-1} \right) .
    \end{multline}
    We choose $D:= N^{0.4} $, which is larger than $ \exp ((9t+1)(\log N)^{1/2}) >Q ^{9t+1}$ for $N\gg _t 1$.
    For this choice the second $O(-)$ term is an
    \begin{align}
        O_{d,M_1,M_2} (D^2 N^{d-1})
        = O_{d,M_1,M_2} (N^{d-0.2})    
        .
    \end{align}   
    Meanwhile we have $\frac{\log D}{\log Q}\ge  0.4 (\log N)^{1/2} $.
    From \eqref{eq:M-and-kappa} and Claim \ref{claim:kernels-dont-contain-each-other}, we know 
    \begin{align}
        M=O_{t,d,K }(L^{O_{t,d,n}(1)}) 
        &=O_{t,d,K } (\exp ((\log N)^{1/2})  ) 
    \end{align} 
    the first $O(-)$ term in \eqref{eq:after-sieve-before-beta} is an
    \begin{align}
        O_{t,d,K }( \exp (-0.3(\log N)^{1/2}) ).
    \end{align}
    Combined with the trivial bounds $\vol (\Omega )\ll _{d} N^{d}$ and $1-g(p)\le 1$,
    the right-hand side of \eqref{eq:after-sieve-before-beta} looks as follows if $N$ is large enough depending on $M_1,M_2$:
    \begin{align}
        \vol (\Omega )
        \left(
            \prod _{p<Q} (1-g(p))
            \right)
            +
        O_{t,d,K }(N^{d} \exp (-0.3 (\log N)^{1/2}) )
        .
    \end{align}
    
    Lastly we multiply \eqref{eq:after-sieve-before-beta} by $\frac 1{\residuezeta ^t}\prod _{p<Q}\left( \frac{p^n}{\totient (p)} \right) ^t $,
    which is an $O_K((\log Q)^t )$ by Rosen's Mertens' Theorem \ref{thm:Mertens-Rosen}.\footnote{
        An asymptotically weaker bound $\prod_{p<Q}\prod _{\idealp \divides p}\left( \frac{\Nrm (\idealp )}{\Nrm (\idealp )-1} \right) ^t = O((\log Q)^{nt})$ suffice here, which follows from the more naive bound $\prod _{\idealp \divides p} \left( \frac{\Nrm 
        (\idealp )}{\Nrm (\idealp )-1} \right) \le \left( \frac{p}{p-1} \right) 
        ^n$ and Mertens' theorem for $\Z $.}
    The right-hand side becomes (compromising the exponent in the error term a little to absorb this $O_K((\log Q)^t )$)
    \begin{align}
        \frac{1}{\residuezeta ^t}\vol (\Omega )
        \left(
            \prod _{p<Q} \beta _p
            \right)
            +
        O_{t,d,K }(N^{d} \exp (-0.2 (\log N)^{1/2}) )
    \end{align}
    and the left-hand side becomes 
    $\sum _{x\in \Omega \cap \Z ^d} \prod _{i=1}^t \Lambda^{\ideala }_{\mCramer ,Q} (\psi _i(x))$
    as was noted at the beginning of this proof.
    This completes the proof of Proposition \ref{prop:sieve-theory}.
\end{proof}

\section{The Siegel model and Mitsui's theorem}\label{sec:Siegel-model}

Now we define the Siegel model $\LambdaSiegel  $, which is a closer approximation to $\Lambda_K$ 
than the \Cramer model is.
To define it we have to recall the notion of Siegel zeros.

\begin{theorem}[Siegel zeros]
    \label{thm:Siegel-zero}
    There exists a positive constant $0< c_K <1$ such that for every $N\ge 1$ and corresponding $Q$: 
    \begin{enumerate}
        \item 
        there is at most one
        Hecke character $\psi $ whose conductor $\idq $ satisfies  
        $\Nrm (\idealq )<Q $
        and such that 
        the $L$-function $L(\psi ,s)$ has zero in the region (where $s=\sigma +it $)
        \begin{align}
            \sigma > 1 - \frac{c_K}{\log Q + \log  (|t|+4)  } .
        \end{align}
        When it exists, let us denote by $\psi _{\mrm{Siegel}}$ the unique character (called the {\em $Q$-Siegel character})
        and its conductor $\idq _{\mrm{Siegel}}$.
        
        \item The Siegel character $\psi _{\mrm{Siegel}}$ is real, namely its values are in $\{\pm 1\}$.
        
        \item 
        The $L$-function $L(\psi_{\mrm{Siegel}},s)$ can have at most one zero 
        $\beta _{\mrm{Sielge}}$ (the {\em $Q$-Siegel zero})
        in the said region. It is necessarily real and 
        for any $0<\epsilon <1$ there is a positive constant $0<c_K(\epsilon )<1$ satisfying:
        \begin{align}
            1 - \frac{c_K}{\log (4Q)}
            < \beta _{\mrm{Siegel}} < 
            1-\frac{c_K(\epsilon )}{\Nrm (\idqSiegel )^{\epsilon }} .
        \end{align}
        Consequently, for any $A>1$ we have $\Nrm (\idqSiegel )\gg _{K,A,s} (\log N)^A$.
    \end{enumerate}
\end{theorem}
\begin{proof}
    See e.g.\ \cite[Theorem 1.9 (p.~277)]{Weiss83}, and \cite[Theorem 5.28 (p.~122)]{Iwaniec-Kowalski} or \cite[Theorem 11.7 (p.~367), Corollary 11.15 (p.~372)]{Montgomery-Vaughan}.
\end{proof}

Of course, the existence of a $Q$-Siegel zero (for any single $Q$) would contradict the Generalized Riemann Hypothesis.
Nonetheless, being unable to disprove its existence so far, 
one has to 
take its potential influence into account
to achieve better approximations.

We always adopt the convention that we ignore the term related to $\psiSiegel $ when it does not exist.

\begin{definition}
    The Siegel model $\Lambda _{\mrm{Siegel},Q} \colon \IdealsK \cup \{ (0)\} \to \R $ is defined by:
    \begin{multline}\label{eq:def-of-LambdaSiegel}
        \LambdaSiegel (\idealc )
        := 
        \LambdaCramer (\idealc ) \cdot \left( 1- \psiSiegel ([\idealc ])\Nrm (\idealc )^{\beta -1} \right) 
        \\ 
        =
        \begin{cases}
            \displaystyle \frac{1}{\residue\limits (\zeta _K(s))}\frac{P(Q)^n}{\totient (P(Q))}
            \left( 1- \frac{\psiSiegel ([\idealc ])}{\Nrm (\idealc )^{1-\betaSiegel }} \right)  
            & \te{ if }\idealc +(P(Q)) = \OK             
            \\ 
            0 & \te{ if not.}
        \end{cases}
    \end{multline}
\end{definition}
Note that the value $\psiSiegel ([\ideala ])\in \{\pm 1\}$ is well defined 
in \eqref{eq:def-of-LambdaSiegel}
because the ideal $\ideala $ is coprime to $\idealq_{\mrm{Siegel}} $ 
under the assumption $\ideala +(P(Q)) = \OK $;
to see this, it suffices to verify that every prime factor $\idealp $ of $\idqSiegel $ is also a prime factor of $(P(Q))$ in $\OK $.
Suppose $\idealp$ has residue characteristic $p$.
We have $p\le \Nrm (\idealp )\le \Nrm (\idqSiegel )<Q$.
This implies that $p$ divides $(P(Q))$, and hence so does $\idealp $.

For a non-zero fractional ideal $\ideala $ we define $\LambdaaSiegel \colon \ideala \to \R $ by 
    \begin{align}
        &\LambdaaSiegel (x)= \LambdaSiegel (x\ideala\inv )
        \\
        &=
        \begin{cases}
            \displaystyle \frac{1}{\residue\limits (\zeta _K(s))}\frac{P(Q)^n}{\totient (P(Q))}
            \left( 1- \frac{\psiSiegel ^{\ida } (x )}{ \Nrm (x\ideala\inv )^{1-\betaSiegel }} \right)  
            & 
            \te{ if $x$ coprime to $P(Q)$ in $\ida $,}
            \\ 
            0 &\te{if not}
            .
        \end{cases}
    \end{align}
    Here, the function $\psiSiegel ^{\ida } \colon \ida \to \{ \pm 1 \}$ is defined by: 
    if $x\in \ida $ is coprime to $\idealqSiegel $, then $\psiSiegel ^{\ida }(x)=\psiSiegel ([x\ida\inv ])$;
    if not, $0$.
    

    By convention, when the Siegel character does not exist, we mean $\LambdaSiegel := \LambdaCramer $ and $\LambdaaSiegel := \LambdaaCramer $.

\begin{proposition}[pointwise bound]
    \label{prop:pointwise-bound}
    Let $N\ge 3$. 
    For all non-zero ideal $\idealc $ of norm $<N$,
    we have the bound 
    \begin{align}
        \Lambda _K(\idealc ),\ \LambdaCramer (\idealc ),\ \LambdaSiegel (\idealc ) \ll _K \log N
        .
    \end{align}
    
\end{proposition}
\begin{proof}
    For $\Lambda _K$, we have the obvious bound 
    $|\Lambda _K(\idealc )| \le \log \Nrm (\idealc ) < \log N $.

    Since $|\LambdaSiegel (\idealc )| \le 2 |\LambdaCramer (\idealc )|$, it remains to verify the assertion for $\LambdaCramer $, which is equivalent to showing 
    $    \frac{P(Q)^n}{\totient (P(Q))} \ll _K \log N $. 
    (Recall $Q$ is a quantity which is much smaller than $N$.)
    By the multiplicativity of $\totient $, the left hand side equals 
    \begin{align}
        \prod _{0<p<Q} \prod _{\idealp \divides p} 
        \left( 1-\frac{1}{\Nrm (\idealp )} \right) \inv .
    \end{align}
    Rosen's Mertens' Theorem \ref{thm:Mertens-Rosen} 
    implies that the previous quantity is bounded as 
    \begin{align}
        \ll _K \log Q <\log N ,
    \end{align}
    which completes the proof.
\end{proof}

\begin{proposition}[Mitsui's Prime Number Theorem]
    \label{prop:essentially-Mitsui}
    Let $\ideala $ be a non-zero fractional ideal,
    $N>2$ a positive number 
    and
    $\Omega \subset (\KR )_{\le N\Nrm (\ida )^{1/n}}$ a convex body.
    Let $\idq $ be a non-zero ideal with 
    $\Nrm (\idq )< \exp ((\log N)^{1/200}) $, 
    and $a\in \ideala /\idq \ideala $ be a residue class.
    Then we have 
    \begin{align}\label{eq:Mitsui-rephrased}
        \sum _{\substack{
            x\in \Omega \cap \ideala \\ x\equiv a \te{ in }\ideala / \idq \ideala  
            }}
            \Lambdaa (x) - \LambdaaSiegel (x)
        \ll _K 
        N^n \exp (-(\log N)^{1/201} )  .
    \end{align}
\end{proposition}
\begin{proof}
    This is essentially Mitsui's Prime Number Theorem with potential Siegel zero incorporated \cite[Theorem 5.1]{KaiMit}. 

    If $\Nrm (\idq )\ge \exp ((\log N)^{1/200})$, 
    then the number of $x\in \Omega \cap \ideala $ is $\ll _K N^n \exp (-(\log N )^{1/200} )$. As the summand is $\ll _K \log N$ by Proposition \ref{prop:pointwise-bound}, the left hand side of \eqref{eq:Mitsui-rephrased} is $\ll _K N^n \exp (-(\log N)^{1/201} )$.
        So we may assume 
        \[ \Nrm (\idq )< \exp ((\log N)^{1/200}). \]
    Also, suppose $a\in \ida /\idq\ida $ generates only a proper $\OK $-submodule of $\ida /\idq\ida $, say contained in $\idp\ida /\idq\ida $ for a prime ideal $\idealp \divides \idq $.
    Then every $x\in \ida $ congruent to $a$ in $\ideala /\idq \ideala $ satisfies $\LambdaaSiegel (x) =0 $
    because $x\ida\inv $ is divisible by $\idealp $ whose residue characteristic is smaller than $\Nrm (\idq )<Q$.
        On the other hand, every such $x$ is in $\idealp\ida $ and hence $\Lambdaa (x)=0$ unless $x$ is a generator of $\idp\ideala $.
        The contribution of such $x$'s (which is only $O_K((\log N)^n)$) is negligible compared to the right-hand side of \eqref{eq:Mitsui-rephrased}.
    It follows that we may assume $a\in \ida /\idq\ida $ is a generator.
    
    \ \ \  
    
    From now on, we assume that $\Nrm (\idq )< \exp ((\log N)^{1/200})$ and $a\in \ida /\idq\ida $ is a generator as an $\OK $-module.
    
    \ \ \ 
     
    First we consider the case $\idq _{\mSiegel }\divides \idq $.
    In this case $\psi _{\mSiegel }$ can be considered as a mod $\idq $ character.
    Mitsui's Prime Number Theorem \cite[Theorem 5.1]{KaiMit} reads 
    \begin{multline}\label{eq:Mitsui-reads}
        \left( \sum _{\substack{
            x\in \Omega\cap \ideala \\ x\equiv a \te{ in }\ideala /\idq \ideala }}
            \Lambdaa  (x) \right) - \\  
            \frac{1}{\Nrm (\ideala )\totient (\idq )}
            \frac{2^{r_2}}{\residue\limits _{s=1}\zeta _K(s) \sqrt{|D_K|}}
            \int _\Omega 1-\psiSiegel (-;a) \Nrm ((-)\ideala\inv )^{\beta _{\mSiegel }-1} d\mu _{\add}
            \\ 
            \ll _K 
            N^n \exp (-(\log N)^{1/2}/O_K(1) )  .
    \end{multline}
    Here, the symbol $\psiSiegel (-;a)$ means the composite 
    \begin{align}
        (\KR )\baci \xto{(\id , a )} (\KR )\baci \times (\ida /\idq\ida )\baci 
        \xto{\incl + \incl } C_K(\idq ) \xto{\psiSiegel} S^1 .
    \end{align}
    For $x\in a+\idq\ida $, via the inclusion $a+\idq\ida \subset (\KR)\baci $ we have 
    $\psiSiegel (x;a)= \psiSiegel ^{\ida }(x)$; see Appendix \ref{sec:Hecke_chars} for some detail.

    The integral $\frac{2^{r_2}}{\sqrt{|D_K|}}\int _\Omega$ is approximated by the sum 
    \begin{align}
        \Nrm (\idq \ideala ) \sum _{\substack{
            x\in \Omega\cap \ideala \\ x\equiv a \te{ in }\ideala /\idq \ideala }} 
            1- \psiSiegel ^{\ida }(x) \Nrm (x\ideala\inv )^{\betaSiegel -1}
        \end{align}
    within an error $O_K(N^{n-\frac 12})$ (see \eqref{eq:theory-of-integration} below for some detail), which is negligible here.
    Up to this error, the left hand side of \eqref{eq:Mitsui-reads} becomes 
    \begin{align}
        \sum _{\substack{
            x\in \Omega\cap \ideala \\ x\equiv a \te{ in }\ideala /\idq \ideala }} 
            \left( \Lambdaa  (x) - 
            \frac{1}{\residue\limits _{s=1}\zeta _K(s)}
            \frac{\Nrm (\idq )}{\totient (\idq )}
            (1- \psiSiegel ^{\ida }(x) \Nrm (x\ideala\inv )^{\betaSiegel -1})
            \right)
            .
    \end{align}
    On the other hand, the left hand side of \eqref{eq:Mitsui-rephrased} is by definition 
    \begin{align}
        \sum _{\substack{
            x\in \Omega\cap \ideala \\ x\equiv a \te{ in }\ideala /\idq \ideala }} 
            \Bigl( \Lambdaa  (x) -  
            \LambdaaCramer (x)
            (1-\psiSiegel ^{\ida }(x)\Nrm (x\ideala\inv )^{\betaSiegel -1} )
            \Bigr) .
    \end{align}
    We have to show that the above two quantities are within the claimed error.
    So we are reduced to showing the following two bounds:
    \begin{multline}\label{eq:we-are-reduced-to-two-1}
        \sum _{\substack{
            x\in \Omega\cap \ideala \\ x\equiv a \te{ in }\ideala /\idq \ideala }} 
            \left(
            \LambdaaCramer (x)
            -
            \frac{1}{\residue\limits _{s=1}\zeta _K(s)}
            \frac{\Nrm (\idq )}{\totient (\idq )}
            \right)
            \\ 
            =
            O_K\left(
                N^n \exp (-(\log N)^{1/201}  )
                 \right) ,
        \end{multline}
        and
    \begin{multline}\label{eq:we-are-reduced-to-two-2}
        \sum _{\substack{
            x\in \Omega\cap \ideala \\ x\equiv a \te{ in }\ideala /\idq \ideala }} 
            \left( 
                \LambdaaCramer (x) - \frac{1}{\residue\limits _{s=1}\zeta _K(s)}
            \frac{\Nrm (\idq )}{\totient (\idq )}
            \right)
            \cdot\psiSiegel ^{\ida }(x)\Nrm (x\ideala\inv )^{\betaSiegel -1}
            \\ 
            =
            O_K\left(
                N^n  \exp (-(\log N)^{1/201}  )
                \right) .
        \end{multline}

        To show \eqref{eq:we-are-reduced-to-two-1}
        (which turns out to not require the condition $\idqSiegel \divides \idq $), we apply Proposition \ref{prop:sieve-theory} to the (single) affine-linear map 
        $\idq \ideala \to \ideala $; $y\mapsto y+a $.
        By 
        \S \ref{appendix:norm-length},
        if we equip the source and target with norm-length compatible bases,
        the inclusion map 
        $\Z ^n \cong \idq \ideala \inj \ideala \cong \Z ^n$ is represented by a matrix with entries of sizes $\le \Nrm (\idq )^{1/n} < \exp ((\log N)^{1/200})$.
        We conclude 
        \begin{align}\label{eq:sieve-applied}
            &\sum _{\substack{
            x\in \Omega\cap \ideala \\ x\equiv a \te{ in }\ideala /\idq \ideala }} 
                \LambdaaCramer (x)
            =
            \sum _{\substack{
            y\in \idq \ideala \\ y+ a\in \Omega  
                    }} 
                \LambdaaCramer (y+a) 
            \\ 
            &=
            \frac{1}{\residuezeta } 
            \left(
                \sum _{\substack{
            x\in \Omega\cap \ideala \\ x\equiv a \te{ in }\ideala /\idq \ideala }}  
             1 
                \right)
                \prod _{\idealp } \beta _\idealp 
                + 
                O_{K}\left(
                    N^n \exp ({-\frac 15 (\log N )^{1/2}})
                \right)
                .
        \end{align}
        Let us evaluate $\beta _\idealp $'s in this case.
        If $\idealp\divides \idealq $, for every $y\in \idq\ideala $ we know that 
        $y+a \not\in \idealp\ideala $ because $a$ generates $ \ideala /\idq\ideala $.
        This implies that $\beta _\idealp = \frac{\Nrm (\idealp )}{\Nrm (\idealp )-1} \cdot 1 = \frac{\Nrm (\idealp )}{\Nrm (\idealp )-1}$.
        Next suppose $\idealp \notdivide \idq $.
        Then since $\idq\ideala $ and $\ideala $ localized at $\mathcal O _{K,\idealp }$ are equal, the map $\idq \ideala /\idealp (\idq\ideala ) 
        \to \ideala / \idealp \ideala $ induced by the inclusion is an isomorphism.
        Hence $\Expec\limits _{y\in \ideala/\idealp(\idealq\ideala )}1[y+a\not\in \idealp\ideala ] = 1-\frac 1{\Nrm (\idealp )}$
        so that $\beta _{\idealp }= \frac{\Nrm (\idealp )}{\Nrm (\idealp )-1}\cdot (1-\frac 1{\Nrm (\idealp )}) =1$.
        It follows that 
        \begin{align}
            \prod _{\idealp }\beta _\idealp = 
            \prod _{\idealp \divides \idq }\frac{\Nrm (\idealp )}{\Nrm (\idealp )-1}
            =
            \frac{\Nrm (\idq )}{\totient (\idq )}.
        \end{align}
        This together with \eqref{eq:sieve-applied} implies \eqref{eq:we-are-reduced-to-two-1}.

        Next we turn to \eqref{eq:we-are-reduced-to-two-2}.
        We will deduce this from \eqref{eq:we-are-reduced-to-two-1}
        by an Abel summation-like argument.

        Under the assumption $\idqSiegel \divides \idq $,
        the function 
        $x\mapsto \psiSiegel ^{\ida }(x)$ on $\{ x\in \ideala \mid x\equiv a \te{ in }\ideala /\idq\ideala \}$ can be seen as the restriction 
        to $\ida $
        of a continuous function $\psiSiegel (-;a)$ on $(\KR )\baci $. 
        Since $\psiSiegel (-;a)$ is a locally constant function into $\{\pm 1\}$, 
        by splitting the convex body $\Omega$ into $O_n(1)$ convex pieces, we may assume that the $\psiSiegel ^{\ida }(x)$ part is constant.
        Thus what we have to show is:
            \begin{multline}\label{eq:we-are-reduced-to-two-2-simplified}
                \sum _{\substack{
                    x\in \Omega\cap \ideala \\ x\equiv a \te{ in }\ideala /\idq \ideala }} 
                    \left( 
                        \LambdaaCramer (x) - \frac{1}{\residue\limits _{s=1}\zeta _K(s)}
                    \frac{\Nrm (\idq )}{\totient (\idq )}
                    \right)
                    \cdot\Nrm (x\ideala\inv )^{\betaSiegel -1}
                    \\ 
                    =
                    O_K\left(
                        N^n \left(  \exp(-(\log N )^{1/201} )\right)
                        \right) .
                \end{multline}
    Let us write 
    \begin{align}
        \Nrm ^{\ideala }(x):= \Nrm (x\ideala\inv ) 
        \quad\te{and}\quad 
        Y:=  \exp ((\log N)^{1/200}) 
        . 
    \end{align}

    We cover the cube 
    $[\pm O_K(1)N]^n_{\ida }\subset \ida $ 
    by disjoint translates $P_i$ ($i\in I$)
    of 
    $[1,N/Y]^n$.
    Thus $\# I=O_K(Y^n)$.
    We have 
    \begin{align}\label{eq:we-are-reduced-to-two-2-simplified-I}
        \te{L.H.S.\ of }\eqref{eq:we-are-reduced-to-two-2-simplified}
        = 
        \sum _{i\in I}
        \sum _{\substack{
                    x\in \Omega\cap P_i \cap \ideala \\ x\equiv a \te{ in }\ideala /\idq \ideala }} 
                    \left( 
                        \LambdaaCramer (x) - \frac{1}{\residue\limits _{s=1}\zeta _K(s)}
                    \frac{\Nrm (\idq )}{\totient (\idq )}
                    \right)
                    \cdot\Nrm ^{\ideala } (x )^{\betaSiegel -1}
                    .
    \end{align}
    
    Let $I_{\mathrm{negl}} \subset I $ be the set of those $i$'s
    such that the cube $P_i$ contains a point $x$ with 
    $|\sigma (x)|< N Y^{-1/2}$ for some $\sigma \colon K\inj \C $.
    It is easy to see that 
    $\# \bigcup _{i\in I_{\mathrm{negl}}} P_i 
    \ll _K N^n Y^{-1/2}
    $. 
    Since $|\LambdaaCramer (x) - \frac{1}{\residue\limits _{s=1}\zeta _K(s)}
    \frac{\Nrm (\idq )}{\totient (\idq )}| \ll _K \log N$ and $\Nrm ^{\ideala }(x)^{\beta -1}\le 1$,
    the contribution of $I_{\mathrm{negl}}$ in \eqref{eq:we-are-reduced-to-two-2-simplified-I} is bounded as:
    \begin{multline}\label{eq:contribution-of-negligible-i's}
        \mgn{\sum _{i\in I_{\mathrm{negl}}}
        \sum _{\substack{
                    x\in \Omega\cap P_i \cap \ideala \\ x\equiv a \te{ in }\ideala /\idq \ideala }} 
                    \left( 
                        \LambdaaCramer (x) - \frac{1}{\residue\limits _{s=1}\zeta _K(s)}
                    \frac{\Nrm (\idq )}{\totient (\idq )}
                    \right)
                    \cdot\Nrm ^{\ideala } (x )^{\betaSiegel -1}
                }        
            \\ 
            \ll _K 
            (\log N)\cdot N^n Y^{-1/2}
            < N^n \exp (-(\log N)^{1/201} )
            .
    \end{multline} 

    Next, for $i\in I\setminus I_{\mathrm{negl}}$, if we write 
    $M_i:= \max _{x\in P_i} \Nrm ^{\ideala }(x) $,
    we have 
    \begin{align}
        \left(
            1-O_K\left( Y^{-1/2} \right)
        \right)
        M_i
        \le 
        \Nrm ^{\ideala }(x)
        \le M_i
        \quad
        \te{ for all } x\in P_i
    \end{align}
    and hence
    \begin{align}\label{eq:norm-does-not-vary-a-lot}
        M_i^{\beta -1}
        \le 
        \Nrm ^{\ideala }(x)^{\beta -1}
        \le 
        M_i^{\beta -1} 
        \cdot \left(
            1+{O_K(Y^{-1/2})}
        \right)
        .
    \end{align}
    Also noting $0<M_i^{\beta -1}\le 1$, we establish the following bound:
    \begin{align}
        &\sum _{i\in I\setminus I_{\te{negl}}}
        \mgn{
            \sum _{ \substack{
                x\in \Omega\cap P_i\cap \ideala \\ x\equiv a \te{ in }\ideala /\idq \ideala 
                }
                }
            \left(
            \LambdaaCramer (x)-
            \frac{1}{\residue\limits _{s=1}\zeta _K(s)}
                    \frac{\Nrm (\idq )}{\totient (\idq )}
            \right)
            \cdot\Nrm ^{\ideala } (x )^{\betaSiegel -1}
        }
        \\ 
        \le 
        &\sum _{i\in I\setminus I_{\te{negl}}}
        \Bigl(
        \mgn{
            \sum _{ \substack{
                x\in \Omega\cap P_i\cap \ideala \\ x\equiv a \te{ in }\ideala /\idq \ideala 
                }
                }
            \left(
            \LambdaaCramer (x)-
            \frac{1}{\residue\limits _{s=1}\zeta _K(s)}
                    \frac{\Nrm (\idq )}{\totient (\idq )}
            \right)
            M_i^{\betaSiegel -1}
        }
        \\
        &+
        \mgn{
            \sum _{ \substack{
                x\in \Omega\cap P_i\cap \ideala \\ x\equiv a \te{ in }\ideala /\idq \ideala 
                }
                }
            O_K\left(
            \log N
            \right)
            \cdot O_K\left( {O_K(Y^{-1/2})} \right)
        }
        \Bigr)
        .
    \end{align}
    Since we can write $\Omega\cap P_i\cap \ideala =\Omega_i\cap \ideala $ for another convex open set $\Omega_i$, we can apply \eqref{eq:we-are-reduced-to-two-1}
    to bound the first term.
    For the second term we use $\# P_i =O_n(N^n /Y^n )$. 
    Let us also recall $\# I =O_K(Y^n)$.
    This implies that the previous quantity is:
    \begin{multline}\label{eq:contribution-of-not-negligible-i's}
        \sum _{i\in I\setminus I_{\te{negl}}}
        \left(
        O_K\left(
                N^{ n}
                \exp ({-({\log N})^{1/2} /O_K (1)})
                \right)
            +
            O_K\left((N^{n}/Y^{n})(\log N ) Y^{-1/2}  
            \right)
        \right)
        \\ 
        = 
        N^{ n}
        O_K\Bigl( 
            Y^n \exp ({-({\log N})^{1/2} /O_K (1)}) 
        + \exp (-(\log N)^{1/201}) \Bigr) ,
    \end{multline}
    which is an $N^{ n}
        O_K( \exp (-(\log N)^{1/201}) ) $.

    Combining this with \eqref{eq:contribution-of-negligible-i's}, 
    we conclude 
    that the estimate \eqref{eq:we-are-reduced-to-two-2}
    holds,
    which completes the proof of Proposition \ref{prop:essentially-Mitsui}
    for the case $\idqSiegel \divides \idq $.

    \ \ \ 
    
    Now consider the (supposedly) easier case $\idqSiegel \notdivide \idq $.
    In this case \eqref{eq:Mitsui-reads} does not have the Siegel term.
    Hence what we have to show is that the correction term for the Siegel zero is negligible in this case:
    \begin{align}\label{eq:what-we-have-to-show-when-not-divide}
        \sum _{\substack{
            x\in \Omega\cap \ideala \\ x\equiv a \te{ in }\ideala /\idq\ideala 
        }}
        \LambdaaCramer (x)\psiSiegel ^{\ida }(x)\Nrm ^{\ideala }(x)^{\beta -1}
        =
        O_K(N^{n} \exp (-(\log N)^{1/201}) ) 
        .
    \end{align}
    Let $\idq_{\mrm{lcm}}= \idqSiegel \cap \idq $ be the least common multiple of $\idqSiegel $ and $\idq $ in the monoid of ideals.
    Write 
    \[ (\ida /\idq\ida )\baci := \{ x\in \ida /\idq\ida \midsep x\te{ generates }\ida /\idq\ida \te{ as an $\OK $-module}\} . \]
    We have natural surjections with kernels of sizes $<Q$
    \begin{align}
        (\ideala /\idqSiegel \ideala )\baci 
        \xfrom{\pr _{\mSiegel}} 
        (\ideala /\idq_{\mrm{lcm}}\ideala )\baci 
        \xto{\pr _{\idq }} 
        (\ideala /\idq\ideala )\baci . 
    \end{align}
    Any sum of the form $\sum _{\substack{
        x\in \Omega\cap \ideala \\ x\equiv a \te{ in }\ideala /\idq\ideala 
    }}$
    can be partitioned as 
    $\sum _{b\in \pr_{\idq}\inv (a)}\sum _{\substack{
        x\in \Omega\cap \ideala \\ x\equiv b\te{ in }\ideala /\idq_{\mrm{lcm}}\ideala 
    }}$.
    By splitting $\Omega$ into its intersections with the connected components of $(\KR )\baci $, we may again assume that $x\mapsto \psiSiegel ^{\ida }(x)$
    is constant on each $(b+\idqSiegel \ida )\cap \Omega $, whose value we will write as $\psiSiegel ^{\ida } (\Omega , b )$ (assuming now $\Omega$ is contained in a connected component of $(\KR )\baci $).
    Hence the left hand side of \eqref{eq:what-we-have-to-show-when-not-divide} becomes:
    \begin{align}
        \sum _{b\in \pr_{\idq}\inv (a)}
            \psiSiegel ^{\ida }(\Omega ,b)
        \sum _{\substack{
        x\in \Omega\cap \ideala \\ x\equiv b\te{ in }\ideala /\idq_{\mrm{lcm}}\ideala 
    }}
        \LambdaaCramer (x)\Nrm ^{\ideala }(x)^{\beta -1}
        .
    \end{align}
    Since $\idqSiegel \divides \idq_{\mrm{lcm}}$, we can apply 
    \eqref{eq:we-are-reduced-to-two-2-simplified} to $\idq_\mrm{lcm}$.
    The previous value now looks like
    \begin{multline}\label{eq:results-applied-to-q-lcm}
        \frac{1}{\residuezeta }\frac{\Nrm (\idq )}{\totient (\idq )}
        \left( \sum _{b\in \pr_{\idq}\inv (a)}
            \psiSiegel ^{\ida }(\Omega ,b)
            \sum _{\substack{
                x\in \Omega\cap \ideala \\ x\equiv b\te{ in }\ideala /\idq_{\mrm{lcm}}\ideala 
            }}
            \Nrm ^{\ideala} (x)^{\beta -1}
        \right)
        \\ 
        +
        O_K(N \exp (-(\log N)^{1/201}) )
        .
    \end{multline}
    
    Choose a norm-length compatible basis $\idq_{\mrm{lcm}}\ida \cong \Zn $ and 
    and let $\Gamma \subset \ideala $ be the corresponding fundamental paralelogram.
    The convex body $\Omega\cap \ideala $ is covered by 
    $O_K(N^{ n}/\Nrm (\idq_{\mrm{lcl}})) $ 
    disjoint translates of $\Gamma $.
    The contribution from the points $x$ with $|\sigma (x)|< N^{1/2} $ for some $\sigma \colon K\inj \C $
    is $O_K((N)^{n-1/3})$, and is negligible as before.
    Outside this region (where the value of $\Nrm ^{\ideala }$ is always $\gg _K N^{n-1/2}$), the value of $\Nrm ^{\ideala}(-)^{\beta -1}$
    varies by at most $\pm O_K(\Nrm (\idq_{\mrm{lcm}}) N^{-1/2} )$ within any translate of $\Gamma $.
    It follows that 
    for any $b\in \pr_{\idq_{\mrm{lcm}}}\inv (a)$:
    \begin{multline}\label{eq:theory-of-integration}
        \sum _{\substack{
                x\in \Omega\cap \ideala \\ x\equiv b\te{ in }\ideala /\idq_{\mrm{lcm}}\ideala 
            }}
            \Nrm ^{\ideala} (x)^{\beta -1}
            =
            \frac{1}{\Nrm (\idq_{\mrm{lcm}}\ideala )} \frac{2^{r_2}}{\sqrt{|D_K|}}\int _\Omega \Nrm ^{\ideala} (x)^{\beta -1} dx
            \\ 
            +
            O_K\left(\frac{N^{ n}}{\Nrm (\idq_{\mrm{lcm}})} \cdot \frac{\Nrm (\idq_{\mrm{lcm}})}{N^{1/2}} \right)
            .
    \end{multline}
    (Yes, we are more or less doing the theory of integration of \Lip functions in a re-scaled environment.)
    If we substitute this into \eqref{eq:results-applied-to-q-lcm}, the problem boilds down to evaluating 
    \begin{align}
        \sum _{b\in \pr_{\idq}\inv (a)}
            \psiSiegel ^{\ida }(\Omega ,b),
    \end{align}
    which is zero because $\psiSiegel ^{\ida }$ restricted to $\pr _{\idq }\inv (a)$
    is not constant by the assumption $\idqSiegel \notdivide \idq $.

    This settles the case of $\idqSiegel \notdivide \idq $ and completes the proof of 
    Proposition \ref{prop:essentially-Mitsui}.
\end{proof}

\section{The two models are close}
\label{sec:Cramer-Siegel-close}

In this section, we use the Weil bound for character sums to prove that 
$\LambdaaCramer $ and $\LambdaaSiegel $ are close in the Gowers norms.

\begin{lemma}[the Weil bound]
    \label{lem:Weil-bound}
    Let $\psi \colon C_K(\idq )\to \{ \pm 1 \}$ 
    be a real primitive Hecke character.
    By the natural inclusion to the non-archimedean part $(\OK /\idq)\baci  \inj C_K(\idq )$
    and zero extension, we regard it as a function on $\OK /\idq $.
    
    Then for any positive number $0<\epsilon <1$, we have 
    \begin{align}
        \lnorm{\psi }_{U^{s+1}(\OK /\idq )}
        \ll _{\epsilon ,s,K}
        \Nrm (\idq )^{-\frac{1}{2^{s+2}}+\epsilon }
        .
    \end{align}
    More concretely, 
    \begin{align}
        \left| \sum _{ \substack{
            x\in \OK /\idq \\  \undl h \in (\OK /\idq )^{s+1} }} 
        \prod _{\undl\omega \in \{ 0,1\}^{s+1}}
        \psi (x+\undl\omega\cdot\undl h)
        \right| 
        \ll _{\epsilon ,s,K}
        \Nrm (\idq )^{(s+2)-\frac{1}{2}+\epsilon }
    \end{align}
    for all $0<\epsilon <1$.
\end{lemma}
\begin{proof}
    Write $\idq = \prod _{\idp }\idp ^{m_\idp }$. Note that then $\Nrm (\idq )=\prod _{\idp }\Nrm (\idp )^{m_\idp }$.
    By Chinese Remainder Theorem, the left-hand side of the assertion is 
    $\prod _{\idp } \lnorm{\psi }_{U^{s+1}(\OK /\idp ^{m_\idp })}$.
    Let us consider these factors separately.

    When the residue characteristic of $\idp $ is odd, 
    we have $m_\idp =1$ because the kernel of 
    $(\OK /\idp ^{m_\idp })\baci \surj (\OK/\idp )\baci =: \F (\idp )\baci $
    has odd order and hence no non-trivial homomorphism into $\{ \pm 1 \}$.

    We have by the definition of $\lnorm{-}_{U^{s+1}}$ and multiplicativity of $\psi $:
    \begin{align}
        \lnorm{\psi }_{U^{s+1}(\F (\idp ))}^{2^{s+1}}
        &=
        \frac{1}{\Nrm (\idp )^{s+2}}
        \sum _{\substack{
                x\in \F (\idp ) \\ \undl{h}\in \F (\idp )^{s+1}
                }}
        \prod _{\undl\omega \in \{ 0,1\}^{s+1} }
        \psi (x+\undl\omega \cdot \undl h)
        \\ 
        &=
        \frac{1}{\Nrm (\idp )^{s+2}}
        \sum _{\substack{
                x\in \F (\idp ) \\ \undl{h}\in \F (\idp )^{s+1}
                }}
                \psi \left(
            \prod _{\undl\omega \in \{ 0,1\}^{s+1} }
            (x+\undl\omega \cdot \undl h)
            \right)
            .
    \end{align}
    The $2^{s+1}$ values $\undl\omega \cdot \undl h$ (for various $\undl\omega \in \{0,1\}^{s+1} $)
    are all different
    except for $O(2^{s+1})\Nrm (\idp )^{s}$ values of $\undl h\in \F(\idp)^{s+1}$.
    In particular the polynomial $\prod _{\undl\omega \in \{ 0,1\}^{s+1} }
    (X+\undl\omega \cdot \undl h)\in \F (\idp )[X]$
    is not a square. Thus for such an $\undl h$, the Weil bound \cite[Theorem 11.23 with n=1]{Iwaniec-Kowalski} applies to give 
    \begin{align}
        \left| 
        \sum _{\substack{
                x\in \F (\idp ) 
                }}
                \psi \left(
            \prod _{\undl\omega \in \{ 0,1\}^{s+1} }
            (x+\undl\omega \cdot \undl h)
            \right)
        \right|
        \le s\Nrm (\idp )^{1/2}
        .
    \end{align}
    For those $O(2^{s+1})\Nrm (\idp )^s$ excluded values of $\undl h$,
    we bound the same sum by the trivial bound $\le \Nrm (\idp )$.
    It follows that for every $\idp $ with odd residue characteristic, we have 
    \begin{align}\label{eq:bound-odd-primes}
        \lnorm{\psi}_{U^{s+1}(\F (\idp))}^{2^{s+1}}
        < s\Nrm (\idp )^{-1/2} + O(2^{s+1})\Nrm (\idp )^{-1}
        ,
    \end{align}
    which is $\ll _s \Nrm (\idp )^{-1/2}$
    and even $< \Nrm (\idp )^{-\frac 12 +\epsilon }$
    for all but finitely many $\idp $ failing to satisfy $\Nrm (\idp )^{\epsilon }> O_s(1)$.

    When the residue characteristic of $\idp $ is $2$ (there are $\le n$ of them),
    it is no longer automatically true that the $\idp $-part of the conductor of $\psi $ satisfies $m_\idp =1$,
    but nonetheless there is an upper bound on $m_\idp $.
    Namely, by the structure of the completion $\calO_{K,\idp}\baci $
    \cite[Proposition 5.7]{Neukirch} there is an exponent $e_\idp $ such that 
    $1+\idp ^{e_\idp }\calO _{K,\idp }$ is contained in the set of squares 
    $(\calO _{K,\idp }\baci )^2$.
    Since any real character $\psi $ is trivial on squares, it follows that 
    $m_\idp \le e_\idp $.
    We conclude 
    \begin{align}\label{eq:bound-even-primes}
        \lnorm{\psi}_{U^{s+1}(\F (\idp ))}^{2^{s+1}}
        \le 1 =\Nrm (\idp ^{m_\idp } )^{1/2}\cdot \Nrm (\idp ^{m_\idp })^{-1/2}
        \ll _K \Nrm (\idp ^{m_\idp })^{-1/2}.
    \end{align}

    Taking the product of \eqref{eq:bound-odd-primes} \eqref{eq:bound-even-primes} over all $\idp $ and recalling 
    $ \prod _{\idp }\Nrm (\idp )^{m_\idp } = \Nrm (\idq )$, 
    we conclude that our assertion holds.
\end{proof}

\begin{proposition}\label{prop:Cramer-Siegel-close}
    Let $\ida \cong \Zn $ be a norm-length compatible basis
    and $\Omega \subset [\pm N]^n_{\ida \otimes \R}$ be a convex set.
    Then  
    \begin{align}
        \lnorm{(\LambdaaSiegel - \LambdaaCramer )\cdot 1_\Omega }_{U^{s+1}[\pm N]^n }
        \ll _{s,K,\epsilon }
        \Nrm(\idqSiegel )^{-\frac{1}{2^{s+2}} +\epsilon }
        \ll _{s,K,A}
        (\log N)^{-A}
    \end{align}
    for any $0<\epsilon < 1/2^{s+4} $ and $A>1$.
\end{proposition}
\begin{proof}
    The latter inequality follows from the last line of Theorem \ref{thm:Siegel-zero}.

    Given $0\le a \le O_K(N^n)$, consider the following subset 
    \begin{align}
        \Omega(a):= \{ x\in \Omega  \mid \Nrm (x\ida\inv )\le a \} ,
    \end{align}
    which we know to have boundary of \Lip class $\mLip (n,O_n(1),O_K(N))$
    (see Lemma \ref{lem:bounded-by-norm-Lipschitz}).
    
    Suppose for the time being that for some $0<c<1$ we have the following inequality 
    \begin{align}\label{eq:inequality-uniformly-for-all-Lip}
        \lnorm{ \LambdaaCramer (-)\psiSiegel ^{\ida }(- ) 
         1_{\Omega_1} (-) }_{U^{s+1}[\pm N]^n}
        \ll _{s,K}
        \Nrm (\idqSiegel )^{-c}
    \end{align}
    uniformly for all subsets $\Omega_1\subset [\pm N]^n_{\ida\otimes \R }$ with boundary of such \Lip class.
    If we write $ M =O_K(N^n)$ for the maximum value of norms on $[\pm N]^n_\ida $,
    we have the following identity of functions on $[\pm N]^n_\ida $:
    \begin{align}
        \Nrm ((-)\ida\inv )^{\betaSiegel -1}1_{\Omega }(-)
        =
         M  ^{\betaSiegel -1}1_{\Omega }(-)
        + 
        (1-\betaSiegel )
        \int _1^{ M }
        y^{\betaSiegel -2} 1_{\Omega(y)}(-) dy .
    \end{align}
    Then we deduce the inequalities below, where the ``$\le $'' follows from the triangle inequality for the Gowers norm and its integral version (noting that $\Vert - \Vert _{U^{s+1}[\pm N]^n}\colon \C ^{[\pm N]^n}\to \R $ is a continuous function with respect to the unique topology on the source) 
    and ``$\ll $'' follows from the assumption \eqref{eq:inequality-uniformly-for-all-Lip}:
    \begin{align}
        & \lnorm{(\LambdaaSiegel - \LambdaaCramer )\cdot 1_\Omega }_{U^{s+1}[\pm N]^n }
        \\ 
        &=
        \lnorm{ \LambdaaCramer (-)\psiSiegel ((-)\ida\inv )
        \Nrm ((-)\ida\inv )^{\betaSiegel -1} 1_{\Omega } (-) }_{U^{s+1}[\pm N]^n}
        \\ 
        &=
        \left\Vert \LambdaaCramer (-)\psiSiegel ((-)\ida\inv ) \cdot \phantom{\int_1^M}\right.
        \\ 
        &\qquad \left. \left[ 
            M  ^{\betaSiegel -1} 1_\Omega (-)
            + 
            (1-\betaSiegel )
            \int _1^{ M }
            y^{\betaSiegel -2} 1_{\Omega(y)}(-) dy 
        \right] \right\Vert _{U^{s+1}[\pm N]^n}
        \\ 
        & \le
        \lnorm{\LambdaaCramer \psiSiegel 1_\Omega}_{U^{s+1}[\pm N]^n}
        \cdot M  ^{\betaSiegel -1}
        \\ 
        &\qquad 
        +
        (1-\betaSiegel )
        \int _1^{ M } 
        y^{\betaSiegel -2} 
        \lnorm{\LambdaaCramer \psiSiegel  1_{\Omega(y)}}_{U^{s+1}[\pm N]^n}
        dy
        \\ 
        &\ll _{s,K}
        \Nrm (\idqSiegel )^{-c}
        \left(
             M  ^{\betaSiegel -1}
            +
        (1-\betaSiegel )
        \int _1^{ M } 
        y^{\betaSiegel -2} 
        dy 
        \right)
        \\ 
        &= \Nrm (\idqSiegel )^{-c}.
    \end{align}
    Thus we are reduced to showing an inequality of the form 
    \eqref{eq:inequality-uniformly-for-all-Lip} with $c=2^{-(s+2)}-\epsilon $.

    To show this, let $\wti\Omega  \subset \left( [\pm 2N]^{n}_{\ida\otimes \R} \right) ^{s+2}$ be the following bounded set:
    \begin{align}
        \wti\Omega  := \{ (x,\undl h)\in \ida\otimes \R \mid x+\undl\omega \cdot \undl h\in \Omega _1\te{ for all }\undl\omega\in \{0,1 \} ^{s+1} 
            \}
        .
    \end{align}
    It has boundary of \Lip class $\mLip (n(s+2),O_{n,s}(1),O_K(N))$.
    The $2^{s+1}$st power of the left-hand side of 
    \eqref{eq:inequality-uniformly-for-all-Lip}, multiplied by the normalization factor $\asymp _s N^{n(s+2)}$, is 
    \begin{align}\label{eq:power-of-left-hand-side}
        \sum _{(x,\undl h)\in \wti\Omega\cap \ida ^{s+2}}
        \prod _{\undl\omega \in \{ 0,1\}^{s+1} }
        \LambdaaCramer (x+\undl\omega\cdot\undl h)\psiSiegel ((x+\undl\omega\cdot\undl h)\ida\inv )
        .
    \end{align}
    Let $\calD \subset \ida $ be a complete set of representatives 
    for $\ida /\idqSiegel \ida $.
    The above sum can be written as a sum 
    $\sum _{(a,\undl b)\in \calD ^{s+2}} \sum _{(x,\undl h)\in \wti\Omega \cap (a,\undl b)+(\idqSiegel \ida )^{s+2}} (-)$.
    Since $\psiSiegel $ is a mod $\idqSiegel $ character, it is legitimate to write 
    $\psiSiegel ((x+\undl\omega\cdot\undl h )\ida\inv )= 
    \psiSiegel ((a+\undl\omega\cdot\undl b )\ida\inv )$
    for all $(x,\undl h)\in (a,\undl b)+(\idqSiegel \ida )^{\oplus s+2}$.
    The sum \eqref{eq:power-of-left-hand-side} thus equals
    \begin{multline}\label{eq:sum-in-power-of-left-hand}
        \sum _{(a,\undl b)\in \calD ^{s+2}} 
        \left[ 
        \left(
            \prod _{\undl\omega\in \{0,1\}^{s+1} }
            \psiSiegel ((a+\undl\omega b)\ida\inv )
        \right)
        \right.
        \\ 
        \left.
        \cdot 
        \sum _{(x,\undl h)\in \wti\Omega \cap (a,\undl b)+(\idqSiegel \ida )^{s+2}}
        \left( \prod _{\undl\omega\in \{0,1\}^{s+1} }
        \LambdaaCramer (x+\undl\omega\cdot\undl h)
        \right)
        \right] 
        .
    \end{multline}
    Note that the first product $\prod _{\undl\omega }$ in this expression is nonzero only when $a,\undl b$ are such that 
    $a+\undl\omega\cdot\undl b $ is coprime to $\idqSiegel $ in $\ida $
    for all $\undl\omega $. So we concentrate on such $a,\undl b$.
    
    For each $a,\undl b$, the latter factor of the expression has already been computed 
    in Proposition \ref{prop:sieve-theory}.
    The volume part from the proposition is now $\vol (\wti\Omega)/ \Nrm (\idqSiegel )^{s+2}$
    ($\ll _s N^{n(s+2)} / \Nrm (\idqSiegel )^{s+2} $) reflecting the fact that we are on the ideal $\idqSiegel \ida $.
    The coefficients $\beta _\idp $ from the proposition {\it a priori} seem to depend on $a,\undl b$,
    but it is easy to see they do not.
    For when $\idp\divides \idqSiegel $, all the elements of $a+\undl\omega\cdot\undl b +\idqSiegel \ida $ are coprime to $\idp $ in $\ida $ for every $\undl\omega $
    by the assumption that $a+\undl\omega\cdot\undl b$ is.
    It follows that $\beta _\idp = (\Nrm (\idp )/\totient (\idp ))^{2^{s+1}}$ regardless of $(a,\undl b)$.
    When $\idp\notdivide \idqSiegel $,
    by the affine-$\OK $-linear isomorphisms induced by the inclusion $\idqSiegel \ida \inj \ida $ and the translations by $a+\undl\omega\cdot\undl b$:
    \begin{align}
        \prod _{\undl\omega } \idqSiegel \ida /\idp \idqSiegel \ida 
        \isoto 
        \prod _{\undl\omega } \ida /\idp\ida 
        \isoto
        \prod _{\undl\omega } \ida /\idp\ida ,
    \end{align}
    it follows that 
    \begin{align}
        \beta _\idp 
        &= (\Nrm (\idp )/\totient (\idp ))^{2^{s+1}}
        \Expec _{(x,\undl h)\in \F (\idp )^{s+2}} 1_{x+\undl\omega\cdot\undl h\neq 0 \te{ in $\F(\idp )$ for all $\undl\omega $} }
        \\ &= 1+O_s\left( \frac{1}{\Nrm (\idp )^2}\right) 
        ,
    \end{align}
    in which $(a,\undl b)$ does not appear.
    It also follows that 
    \begin{align}
        \prod _{ \substack{
            p<Q \\ \idp \divides p } }
        \beta _\idp 
        \ll _K 
        \prod _{\idp \divides \idqSiegel } \left( 1-\frac{1}{\Nrm (\idp )} \right) ^{-2^{s+1}}
        \ll _K 
        \log \Nrm (\idqSiegel ) < \log N
    \end{align}
    by Mertens' theorem \ref{thm:Mertens-Rosen} say.
    We conclude that the sum $\sum _{(x,\undl h)}$ in 
    \eqref{eq:sum-in-power-of-left-hand}
    has the following form regardless of $(a,\undl b)$:
    \begin{align}\label{eq:regardless-of-(a,b)}
        \frac{\vol (\wti\Omega)/\Nrm (\idqSiegel )^{s+2}}{\residue\limits _{s=1}(\zeta _K(s))}\prod _{p<Q}\beta _p
         +
         O_{s,K}(N^{n(s+2)} \exp (-\frac{1}5 (\log N)^{1/2}) ).
    \end{align}
    The contribution of first term of this expression to 
    \eqref{eq:sum-in-power-of-left-hand}
    is, by Lemma \ref{lem:Weil-bound},
    \begin{align}
        \ll _{\epsilon ,s,K} 
        \Nrm (\idqSiegel )^{(s+2)-\frac 12 +\epsilon }
        \cdot \frac{N^{n(s+2)}}{\Nrm (\idqSiegel )^{s+2}}
        \cdot \log N .
    \end{align}
    By Theorem \ref{thm:Siegel-zero}, the $\log N$ factor can be absorbed in $\Nrm (\idqSiegel )^{\epsilon }$ by modifying $\epsilon $ a little, and so we end up with a bound
    \begin{align}\label{eq:first-term-1/2+epsilon}
        \ll _{\epsilon , s,K} N^{n(s+2)}\Nrm (\idqSiegel )^{-\frac 12 +\epsilon }. 
    \end{align}
    The contribution of the error term in \eqref{eq:regardless-of-(a,b)} 
    to \eqref{eq:sum-in-power-of-left-hand}
    is 
    \begin{align}\label{eq:second-term-N(q)^-5}
        &\ll _{s,K} 
        N^{n(s+2)}\Nrm (\idqSiegel )^{s+2} \exp (-\frac{1}5 (\log N)^{1/2}) 
        \\ 
        &< N^{n(s+2)} \exp (-(\log N)^{1/3}) .
    \end{align}

    The inequalities \eqref{eq:first-term-1/2+epsilon} \eqref{eq:second-term-N(q)^-5}
    give the assertion raised to the $2^{s+1}$st power.
    This completes the proof of Proposition \ref{prop:Cramer-Siegel-close}.
\end{proof}

\section{Vaughan's decomposition (a.k.a.\ the Type I/II sums)}
\label{sec:Vaughan}

Having established that $\LambdaaCramer $ and $\LambdaaSiegel $ are close in the $U^{s+1}$ norm (Proposition \ref{prop:Cramer-Siegel-close}),
we move on to the task of showing that $\lnorm{\Lambdaa  - \LambdaaSiegel 
 }_{U^{s+1}[\pm N]^n}$ is small.
For this it will be convenient to decompose $\Lambdaa $ and $\LambdaSiegel $ into manageable chunks in the style of Vaughan.

Recall the definition of the divisor function 
$\tau \colon \IdealsK \to \N $
for the number field:
\begin{align}
    \tau (\idealc ) = \sum _{\ideald \divides \idealc } 1 ,
\end{align}
where $\idd $ runs through the non-zero ideals of $\OK $ dividing $\idc $.

\begin{lemma}\label{lem:divisor-moment}
    For any integers $N\ge 3$ and $m\ge 1$ 
    we have 
    \begin{align}
        \sum _{\substack{\idealc \ \mrm{ such that}\\ \Nrm (\idealc )\le N}}
        \tau (\idc )^m
        &\ll _{K,m} N (\log N)^{2^m-1}
        \\ 
        \sum _{\substack{\idealc \ \mrm{ such that}\\ \Nrm (\idealc )\le N}} \frac{\tau (\idealc )^m}{\Nrm (\idealc )} 
        &\ll _{K,m} ( \log N)^{2^m}.
    \end{align}
\end{lemma}
\begin{proof}
    See e.g.\ \cite[Appendix C]{QuadraticUniformity} and \cite[p.~288]{Ruzsa} for the case of $\Z $.
    The case of $\OK $ can be treated in the same way, noting 
    $\# \{ \idealc \subset \OK \ :\  \Nrm (\idealc )\le N  \} \ll_K N $
    (Proposition \ref{prop:density-of-ideals}).
\end{proof}

Now we introduce the notion of 
(twisted) Type I/II sums for $K$.
We will give the minimum definitions that cover the objects of our interest.

\begin{definition}\label{def:divisor-bounded-functions}
    Let $N\ge 3$ be an integer and $A\colon \IdealsKN \to \C $ be a function.
    It is said to be {\it log-bounded} if we have a bound of the form 
    \begin{align}
        |A(\idealc )| \le \log  N 
    \end{align}
    for all $\idealc \in \IdealsKN $.
\end{definition}

\begin{definition}\label{def:type-I/II}
    A function $S\colon \IdealsKN \to \C $ for some $N\ge 3$ is said to be:
    \begin{itemize}
        \item a {\it Type I sum} if it is of the form 
            \begin{align}
                S(\idealc ) = 1_{\Nrm (\idealc )\le N'} 
                \sum _{\substack{\ideald \te{ such that } \\ 
                                   \Nrm (\ideald )\le N^{2/3} 
                                   }}
                A (\ideald )1_{\ideald\divides \idealc } ,    
            \end{align}
            where $A \colon \IdealsKN \to \C $ is log-bounded and $N'\le N$.
        \item a {\it twisted Type I sum} if it is of the form
            \begin{align}
                S(\idealc ) = \psiSiegel ([\idealc ])
                1_{\Nrm (\idealc )\le N'} 
                \sum _{\substack{\ideald \te{ such that } \\ 
                                   \Nrm (\ideald )\le N^{2/3} 
                                   }}
                A(\ideald )1_{\ideald\divides \idealc }  ,     
            \end{align}
            where $A\colon \IdealsKN \to \C$ is log-bounded and $N'\le N$;
            we set $\psiSiegel  ([\idc ]):=0$ if the ideal $\idc $ is not coprime to $\idqSiegel $.
            
        \item a {\it Type II sum} if it is of the form
            \begin{align}
                S(\idealc )= \sum _{\substack{\ideald _1 ,\ideald _2 \te{ such that} \\ 
                                        \Nrm (\ideald _1),\Nrm (\ideald _2) >N^{1/3} \\ 
                                        \te{ and } \ideald _1 \ideald _2= \idealc 
                                                }}
                            L(\ideald _1 )T(\ideald _2)
            \end{align}
            where $L,T \colon \IdealsKN \to \C $ obey the bounds 
            \begin{align}
                |L(\ideald _1 )| &\le \log N,
                \\ 
                |T(\ideald _2 )| & \le \tau (\ideald _2 )
                .
            \end{align}
        \item {\it negligible} if it is log-bounded and satisfies
            \begin{align}
                \sum _{\substack{\idealc \te{ such that} \\ \Nrm (\idealc )\le N}}
                |S (\idealc )|
                < N\exp (-\frac{1}{10}(\log N)^{1/2}  ) .
            \end{align}
    \end{itemize}
\end{definition}

\begin{proposition}[pointwise bound of (twisted) Type I, II sums]
    \label{prop:pointwise-bound-of-type-I/II-sums}
    For any (twisted) Type I, II sum 
    $S\colon \Ideals _{K,\le N}\to \C $,
    we have 
    \begin{align}
        |S(\idc )|\le (\log N )\tau (\idc )^2
    \end{align}    
    for all $\idc \in \Ideals _{K,\le N}$.
\end{proposition}
\begin{proof}
    First suppose $S$ is a (twisted) Type I sum.
    In either case, by $|A(-)|\le \log N$ and $|\psiSiegel (-)|\le 1$
    we have 
    \begin{align}
        |S(\idc )| 
        \le \log N 
        \sum _{\substack{
            \idd \\ \Nrm (\idd )\le N^{2/3}
       }}
       1_{\idd \divides \idc }
       \le 
       (\log N )\tau (\idc ).
    \end{align}
    
    Next suppose $S$ is a Type II sum. We have 
    \begin{align}
        |S(\idc )|
        =
        \left| 
            \sum _{\substack{\ideald _1,\ideald _2 \te{ such that} \\ 
                                        \Nrm (\ideald _1 ),\Nrm (\ideald _2 ) >N^{1/3} \\ 
                                        \te{ and } \ideald_1 \ideald _2 = \idealc 
                                                }}
                            L(\ideald _1 )T(\ideald _2)
        \right|
        \le 
        \log N 
        \left|  
            \sum _{\idd _2 \divides \idc } |T(\idd _2 )|
        \right|
        &\le 
        \log N \sum _{\idd _2 \divides \idc } \tau (\idd _2 )
        \\ 
        &\le 
        (\log N)\tau (\idc )^2 .
    \end{align}
    This completes the proof.
\end{proof}

\begin{proposition}[Vaughan's decomposition]\label{prop:Vaughan-decomposition}
    For any $N\ge 3$ sufficiently large depending on $K$, the functions $\Lambda _K$, $\LambdaCramer $, $\LambdaSiegel $ restricted to $\IdealsKN $ can be written as  
    $\C $-linear combinations of the form
        \begin{align}\label{eq:Vaughan-decomposition}
            \sum _{i\in I} c_i S_i ,
        \end{align}
    where 
    \begin{itemize}
        \item the coefficients $c_i\in \C $ satisfy $\sum _{i\in I} |c_i|\ll _K \log N$;
        \item each $S_i\colon \IdealsKN \to \C$ is either a Type I sum, a twisted Type I sum, a Type II sum or negligible.
    \end{itemize}
    Twisted Type I sums only appear in the case of $\LambdaSiegel $.
\end{proposition}
\begin{proof}
    First, we consider $\Lambda _K$. By the \Mobius inversion formula 
    \begin{align}
        \Lambda _K = \Lambda _K*\mu * 1 ,
    \end{align}
    for all ideal $\idealc $ we have 
    \begin{align}
        \Lambda _K(\idealc )
        = \sum _{\substack{\ideald _1,\ideald _2, \ideald _3 \\ \te{such that }\ideald_1\ideald_2\ideald _3 =\idealc }}
        \Lambda _K(\ideald_1 )\mu (\ideald_2 ) 1
        =
        \sum _{\substack{\ideald_1, \ideald_2 \\ \te{such that }\ideald_1\ideald_2 \divides \idealc }}
        \Lambda _K(\ideald_1 )\mu(\ideald_2 ) .
    \end{align}
    We partition the sum into four partial sums according to whether $\Nrm (\ideald_1 ), \Nrm (\ideald_2 )\le N^{1/3}$.
    Out of these four partial sums, let us treat the following three first.
    \begin{claim}\label{claim:partial-sums}
        \begin{enumerate}
            \item The sum 
                \begin{align}
                    \sum _{\Nrm (\ideald_1 ),\Nrm (\ideald_2 )\le N^{1/3}}
                    \Lambda _K(\ideald_1 )\mu(\ideald_2 )
                \end{align}
                is a Type I sum.
            \item The sum 
            \begin{align}
                \left( \sum _{\Nrm (\ideald_1 ),\Nrm (\ideald_2 )\le N^{1/3}}
                + 
                \sum _{\Nrm (\ideald_1 )>N^{1/3},\Nrm (\ideald_2 )\le N^{1/3}}
                \right)
                \Lambda _K(\ideald_1 )\mu(\ideald_2 )
            \end{align}
            can be written as a linear combination $\sum _{i\in I}c_iS_i$ where $\sum _{i\in I}|c_i|\le 2 \log N$ and all $S_i$ are Type I sums.
            \item The sum 
            \begin{align}
                \left( \sum _{\Nrm (\ideald_1 ),\Nrm (\ideald_2 )\le N^{1/3}}
                + 
                \sum _{\Nrm (\ideald_1 )\le N^{1/3},\Nrm (\ideald_2 )> N^{1/3}}
                \right)
                \Lambda _K(\ideald_1 )\mu(\ideald_2 )
            \end{align}
            is negligible.
        \end{enumerate}
        \label{claim:type-I/II-three-terms}
    \end{claim}
    \begin{proof}[Proof of Claim]
        For (1), we rearrange the sum as follows by setting $\ideald := \ideald_1\ideald_2 $ 
        \begin{align}
            \sum _{\substack{\ideald \te{ such that} \\ \Nrm (\ideald )\le N^{2/3}}}
            1_{\ideald \divides \idealc }
            \left(
                \sum _{\substack{\ideald_1,\ideald_2 \te{ such that} \\ \ideald_1\ideald_2 =\ideald \\ \Nrm (\ideald_1 ),\Nrm (\ideald_2 )\le N^{1/3}} }
                \Lambda _K(\ideald_1 )\mu (\ideald_2 ) .
            \right)
        \end{align}
        By setting $N'=N^{2/3}$ in Definition \ref{def:type-I/II}, it suffices to show that the sum $\sum _{\substack{\ideald_1,\ideald_2 }}$ 
        in the parenthesis above is log-bounded in the sense of Definition \ref{def:divisor-bounded-functions}. Name it $A$:
        \begin{align}
            A(\ideald )= \sum _{\substack{\ideald_1,\ideald_2 \te{ such that} \\ \ideald_1\ideald_2 =\ideald \\ \Nrm (\ideald_1 ),\Nrm (\ideald_2 )\le N^{1/3}} }
            \Lambda _K(\ideald_1 )\mu (\ideald_2 ) 
        \end{align} 
        By the formula $\Lambda _K*1=\log \Nrm (-)$ we have 
        \begin{align}
            |A(\ideald )| \le \sum _{\substack{\ideald_1,\ideald_2 \te{ such that} \\ \ideald_1\ideald_2 =\ideald } }
            \Lambda _K(\ideald_1 )\cdot 1
            = \log \Nrm (\ideald )\le \log N .
        \end{align}
        This shows Claim \ref{claim:type-I/II-three-terms}(1).

        We turn to (2) of Claim \ref{claim:partial-sums}. First, note that the sum equals
        \begin{align}
            \sum _{\substack{\ideald_1 ,\ideald_2 \te{ such that}\\ \ideald_1\ideald_2 =\idealc \\ \Nrm (\ideald_2 )\le N^{1/3}} }
            \Lambda _K(\ideald_1 )\mu (\ideald_2 )
            = 
            \sum _{\substack{\ideald_2 \divides \idealc \\ \Nrm (\ideald_2 )\le N^{1/3}}} 
            \mu (\ideald_2 )
            \left(
                \sum _{\ideald_1 \divides \idealc \ideald_2\inv }
                \Lambda (\ideald_1 )            
                \right)
            \end{align}
            By the formula $\Lambda _K*1=\log \Nrm (-)$ again, we can further compute as 
            \begin{align}
                = \sum _{\substack{\ideald_2 \te{ such that}\\ \Nrm (\ideald_2 )\le N^{1/3} }}
                \mu (\ideald_2 ) 1_{\ideald_2\divides \idealc } \log \Nrm (\idealc \ideald_2\inv )
                = \sum _{\substack{\ideald_2 \te{ such that}\\ \Nrm (\ideald_2 )\le N^{1/3} }}
                \mu (\ideald_2 ) 1_{\ideald_2\divides \idealc } 
                \left( \log \Nrm (\idealc ) -\log \Nrm(\ideald_2 ) \right) .
            \end{align}
            The sum arising from $\log \Nrm (\ideald_2 )$'s is a Type I sum 
            by setting $N'=N$ and 
            $A(\ideald_2 )= 1[{\Nrm (\ideald_2 )\le N^{1/3}}]\mu (\ideald_2 )\log \Nrm (\ideald_2 )$.
            To treat the sum arising from $\log (\Nrm (\idealc ))$'s, we use the identity 
            \begin{align}
                \log (\Nrm (\idealc ) )
                =
                \log N + \sum _{t=1}^{N-1} (\log t - \log (t+1)) \cdot 1_{\Nrm (\idealc )\le t}
            \end{align}
            to write the sum in question as 
            \begin{align}
                \left( 
                    \sum _{\substack{\ideald_2 \te{ such that}\\ \Nrm (\ideald_2 )\le N^{1/3}}}
                    \mu (\ideald_2 )1_{\ideald_2\divides \idealc }
                \right)\log N
                + \sum _{t=1}^{N-1} (\log t - \log (t-1)) 
                \left( 
                    1_{\Nrm (\idealc )\le t}
                    \sum _{\substack{\ideald_2 \te{ such that}\\ \Nrm (\ideald_2 )\le N^{1/3}}}
                    \mu (\ideald_2 )1_{\ideald_2\divides \idealc } 
                \right)
                .
            \end{align}
            As the function $A(\ideald_2):= 1[\Nrm (\ideald_2 )\le N^{1/3}]\mu (\ideald_2 )\in \{ -1,0,1 \}$ is log-bounded in the sense of Definition \ref{def:divisor-bounded-functions}, each sum $\sum _{\ideald_2 }$ is a Type I sum.
            The coefficients $\log N$, $\log t - \log (t+1)$ satisfy the property that the sum of their absolute values is $2\log N$.
            This proves Claim \ref{claim:type-I/II-three-terms}(2).

            Lastly we show (3) of Claim \ref{claim:partial-sums}.
            The sum in question equals 
            \begin{align}
                \sum _{\substack{\ideald_1 ,\ideald_2 \te{ such that}\\ \ideald_1\ideald_2\divides \idealc \\ \Nrm (\ideald_1)\le N^{1/3} }}
                \Lambda _K(\ideald_1)\mu(\ideald_2)
                =
                \sum _{\substack{\ideald_1 \te{ such that}\\ \ideald_1\divides \idealc \\ \Nrm (\ideald_1 )\le N^{1/3}}}
                \Lambda _K(\ideald_1 ) 
                \sum _{\substack{\ideald_2 \te{ such that}\\ \ideald_2\divides \idealc\ideald_1 \inv }} 
                \mu (\ideald_2 ) .
            \end{align}
            By the identity $\sum _{\substack{\ideald_2 \te{ such that}\\ \ideald_2\divides \idealc\ideald_1\inv }} 
            \mu (\ideald_2 ) = 1[\ideald_1 = \idealc ]$ (the \Mobius inversion formula), this equals 
            \begin{align}
                1[\Nrm (\idealc )\le N^{1/3}] \Lambda _K(\idealc ) ,
            \end{align}
            which is negligible in the sense of Definition \ref{def:type-I/II} because 
            \begin{align}
                \sum _{
                    \substack{\idealc \te{ such that}\\ \Nrm (\idealc )\le N }
                }
                |1[\Nrm (\idealc )\le N^{1/3}] \Lambda _K(\idealc ) |
                \le 
                \sum _{ \substack{
                    \idealc \te{ such that}\\ \Nrm (\idealc )\le N^{1/3}
                }}
                \log N
                \ll _K N^{1/3} \log N .
            \end{align}
            This completes the proof of Claim \ref{claim:type-I/II-three-terms}.
    \end{proof}

To finish the $\Lambda _K$ part of Proposition \ref{prop:Vaughan-decomposition},
it remains to treat the forth partial sum 
\begin{align}
    \sum _{\substack{ \ideald_1,\ideald_2 \te{ such that} \\ \ideald_1\ideald_2 \divides \idealc \\ \Nrm (\ideald_1 ),\Nrm (\ideald_2 )>N^{1/3} } }
    \Lambda _K(\ideald_1 )\mu (\ideald_2 ) .
\end{align}
We show that this is a Type II sum.
By introducing the dummy variable $\ideald := \idealc \ideald_1 \inv $, rewrite the sum as 
\begin{align}
    \sum _{\substack{
        \ideald_1 ,\ideald \te{ such that} \\ \ideald_1\ideald = \idealc \\ \Nrm (\ideald_1 )> N^{1/3}
    }}
    \Lambda _K(\ideald_1 ) 
    \sum _{\substack{
        \ideald_2 \te{ such that} \\ \ideald_2\divides \ideald \\ \Nrm (\ideald_1) > N^{1/3}
        }}
    \mu (\ideald_2 ) .
\end{align}
Note that the conditions $\Nrm (\ideald_2 )>N^{1/3}$, $\ideald_2\divides \ideald $ imply $\Nrm (\ideald )>N^{1/3}$.
Hence we are reduced to showing that the functions
\begin{align}
    L(\ideald_1 ):= \Lambda _K(\ideald_1 )
    \quad
    \te{and}
    \quad 
    T(\ideald ):= \sum _{\substack{
        \ideald_2 \te{ such that} \\ \ideald_2\divides \ideald \\ \Nrm (\ideald_2) > N^{1/3}
        }}
        \mu (\ideald_2 ) 
    \end{align}
obey the bounds $|L(\ideald_1 )|\le \log N$ and $|T(\ideald )|\le \tau (\ideald )$.
But both of them are obvious.
This completes the proof of the assertion for $\Lambda _K$.

We next show the assertion for $\LambdaCramer $, $\Lambda _{\mrm{Siegel},Q}$.
Since the coefficient $\frac{P(Q)^n}{\totient (P(Q))} $ has size 
$\ll _K \log N $ by Mertens' theorem, it suffices to show that the functions 
\begin{align}
    \idc &\mapsto 1[\idc +(P(Q)) =\OK ] \\ 
    \idc &\mapsto 1[\idc +(P(Q)) =\OK ] \cdot \psiSiegel ([\idc ])\cdot \Nrm(\idc )^{\betaSiegel -1} 
\end{align}
are of the asserted form.
Set $D:= N^{0.2}$ and let $(\lambda ^+ _\idd )_{\idd \in \Ideals_K}$ be the associated upper linear sieve coefficient from Lemma \ref{lem:sieve-theory} (2).
We know 
\begin{align}\label{eq:upper-linear-sieve-coefficient}
    1_{\idealc +P(Q)=(1)}
    \le 
    \sum _{\idd \divides P(Q) }\lambda ^+_\idd 1_{\idd \divides \idealc }
\end{align}
for all $\idealc $.
The right hand side is a Type I sum in the sense of Definition \ref{def:type-I/II}
by setting $A(\idd ):= \lambda ^+_\idd 1_{\idd \divides P(Q)}$
and noting that 
$\lambda ^+$ is supported on $\Ideals_{K,\le N^{0.2}}\subset \Ideals_{K,\le N^{2/3}}$ and has values in $\{ -1,0,1 \}$.

We want to apply Lemma \ref{lem:sieve-theory} to 
$a_\idn := 1_{\Nrm (\idn )\le N}$.
The theorem of density of ideals \ref{prop:density-of-ideals} 
says 
\begin{align}
    \sum _{\idd \divides \idn } 1_{\Nrm (\idn )\le N}
    &\overset{\idealm := \idn\idd\inv }{=} 
    \sum _{\idealm \in \Ideals _K} 1 _{\Nrm (\idealm )\le N/\Nrm (\idd )}
    \\ 
    &=
    \residue_{s=1} (\zeta _K(s)) 
    \frac{N}{\Nrm (\idd )}
    + O_K\left( \left(\frac{N}{\Nrm (\idd )}\right) ^{1-\frac{1}{n}} \right) ,
\end{align}
which means 
that if we set 
$g(\idd ):=1_{\idd \divides P(Q)} / \Nrm (\idd )$ and 
$X:= \residue _{s=1}(\zeta _K(s)) \cdot N$,
we have the estimate $r_\idd = O_K\left( \left(\frac{N}{\Nrm (\idd )}\right) ^{1-\frac{1}{n}} \right)$.

By Lemma \ref{lem:sieve-theory} (2), we have,
noting \eqref{eq:upper-linear-sieve-coefficient} and $\frac{\log D}{\log Q} \ge 
0.2 \sqrt{\log N}$ 
\begin{align}
    &\sum _{\idc , \Nrm (\idc )\le N} 
    \mgn{
        \left( \sum _{\idd \divides P(Q)} \lambda ^+_\idd 1_{\idd \divides \idc } \right)
        -    
        1_{\idc +P(Q)=(1)}
    }
    \\ 
    =
    &S^+(Q)-S(Q)
    \\ 
    =
    &O_K(Ne^{-0.2\sqrt{\log N}} + N^{1-\frac{0.8}{n}} )
    = O_K(Ne^{-\sqrt{\log N}} ) .
\end{align}
This expresses $1_{\idc +P(Q)=(1)}$ as the sum of a Type I sum and a negligible sum in the sense of Definition \ref{def:type-I/II}.

Multiplying by $\psiSiegel $, we see that the function 
\begin{align}
    \idc &\mapsto 1[\idc +(P(Q)) =\OK ] \cdot \psiSiegel ([\idc ])
\end{align}
is the sum of a twisted Type I sum $\idc \mapsto \psiSiegel ([\idc ])\sum _{\Nrm(\idd )\le D}\lambda ^+_{\idd } 1_{\idd \divides \idc } $
and a negligible sum.

We further multipy it by 
$\Nrm (\idc )^{\betaSiegel -1}$. 
The negligible part stays negligible because $\Nrm (\idc )^{\betaSiegel -1} \le 1$.
To handle the twisted Type I part,
just as before, we use the identity 
\begin{align}
    \Nrm (\idc )^{\betaSiegel -1}
    =
    N^{\betaSiegel -1 }
    +
    \sum _{t=1}^{N-1} (t^{\betaSiegel -1} - (t+1)^{\betaSiegel -1})\cdot 1_{\Nrm (\idc )\le t}
    .
\end{align}
Note that for any twisted Type I sum $S(- )$ and $t\le N$,
the product $S_{\le t}(\idc ):= 1[\Nrm (\idc )\le t]S(\idc )$ is also a twisted Type I sum.
Thus $\idc \mapsto \Nrm (\idc )^{\betaSiegel -1} S(\idc ) $ is written as a linear combination of 
\begin{align}
    S(\idc ):= \sum _{\Nrm(\idd )\le D}\lambda ^+_{\idd } 1_{\idd \divides \idc } \psiSiegel ([\idc ]) 
\end{align}
and the associated $S_{\le t}(\idc )$'s with coefficients 
$N^{\betaSiegel -1}$ ($\le 1$) and $(t^{\betaSiegel -1} - (t+1)^{\betaSiegel -1})$'s.
Since the sum of the absolute values of these coefficients is 
$N^{\betaSiegel -1}+( 1-N^{\betaSiegel -1})=1$, 
the proof for the case of $\LambdaSiegel $ is now complete.
\end{proof}

\section{The von Mangoldt function and its models are close: first reductions}\label{sec:Mangoldt-model-close}

We are finally ready to start the proof of one of our main goals, namely the proximity of the von Mangoldt function to its models in the Gowers norms:

\begin{theorem}  \label{thm:vonMang-Siegel-close!} 
    Let $\ideala \subset K $ be a non-zero fractional ideal and $N\gg _K 1$ be a large enough positive integer.
    There is a positive number $0<c_{s,n}<1$ such that for every convex set $\Omega \subset [\pm N]^n\subset \ideala $ with respect to a norm-length compatible basis,
    non-zero ideal $\idealq \subset \OK $ and a congruence class $a+\idealq \ideala  $, the following inequality holds:
    \begin{align}
            \lnorm{
            1_{\Omega \cap (a+\idealq \ideala )} \cdot 
            (\Lambda_{K}^{\ideala } - \LambdaSiegel ^{\ideala } )
            }_{U^{s+1}([\pm N]^n_{\ida })} 
            &\le 
            \exp (-(\log N)^{c_{s,n}}).
        \end{align}
    
\end{theorem}

The proof of Theorem \ref{thm:vonMang-Siegel-close!} will be finished in \S \ref{sec:Type-II}.

\begin{corollary}\label{cor:vonMang-Cramer-close!}
    Under the same setting as in Theorem \ref{thm:vonMang-Siegel-close!},
For any positive numbers $N\gg _{s,K} 1$, $A>1$ and convex set $\Omega \subset [\pm N]^n_{\ida }$, 
we have 
    \[
    \lnorm{
            1_{\Omega } \cdot 
            ( \Lambda_{K}^{\ideala } - \LambdaCramer ^{\ideala } )
            }_{U^{s+1}([\pm N]^n_{\ida })} 
            \ll _{A}
            (\log N)^{-A}.
    \]
\end{corollary}
\begin{proof}
    This is the combination of Theorem \ref{thm:vonMang-Siegel-close!} and 
    Proposition \ref{prop:Cramer-Siegel-close}.
\end{proof}

\subsection{Relation to asymptotic orthogonality with nilsequences}\label{sec:relation-to-asymptotic-orthogonality}
\label{sec:reduction-to-asymptotic-orthogonality}

Note the following consequence of Theorem \ref{thm:vonMang-Siegel-close!}.

\begin{theorem}\label{thm:asymptotic-orthogonality-with-nilsequences}
    Let $K$ be a number field of degree $n$ and $\ideala \subset K$ be a non-zero fractional ideal.
    Let $0\le s\le k$ be integers.
    There exists a positive number $0<\mydecayexponent _{s,k,n}<1/100$ such that 
    if $N\gg _{s,k,K} 1$ is sufficiently large, then for every choice of 
\begin{itemize}
    \item 
    convex set $\Omega \subset [\pm N]^n_\ideala $ with respect to a norm-length compatible basis $\ida \cong \Zn $, 
    \item congruence class $a+\idealq \ideala $ with respect to a non-zero ideal $\idealq \subset \OK $,
    \item nilmanifold $G/\Gamma $ of dimension $\le (\log N)^{\mydecayexponent _{s,k,n}}$, step $\le s$, degree $\le k$ and complexity $\le \exp ((\log N)^{\mydecayexponent _{s,k,n}})$,
    \item function $F\colon G/\Gamma \to \C $ with Lipschitz norm 
    $\le \exp ((\log N)^{\mydecayexponent _{s,k,n}})$,
\end{itemize}
    we have 
    \begin{align}\label{eq:asymptotic-orthogonality-with-nilsequences}
        \left\vert
               \Expec _{x\in [\pm N]^n_\ida }
               1_{\substack{\Omega \cap (a+\idealq \ideala )}} (x)
               \cdot 
               (\Lambda _K ^{\ideala }(x) - \LambdaSiegel ^{\ideala }(x))
               \cdot 
               F(g (x)\Gamma )
        \right\vert
        \\ 
        \le  
        \exp (-(\log N)^{\mydecayexponent _{s,k,n}}) .
    \end{align}
\end{theorem}
Since $s\le k$, the assertion is equivalent to the existence of the exponent $\mydecayexponent _{s,n}:= 
\mydecayexponent _{s,s,n}$, 
but for the proof it is convenient to distinguish the values of $s\le k$.
\begin{proof}[Proof of Theorem \ref{thm:asymptotic-orthogonality-with-nilsequences} assuming Theorem \ref{thm:vonMang-Siegel-close!}]
    If the inequality \eqref{eq:asymptotic-orthogonality-with-nilsequences} fail with some $0<\mydecayexponent <1/100$, 
    we apply Proposition \ref{prop:converse-to-inverse-theorem} to $\varepsilon\inv = M = \exp ((\log N)^\mydecayexponent )$ and $d=(\log N)^{\mydecayexponent }$ 
    to obtain
    \begin{align}
        \lnorm{1_{\Omega \cap (a+\idealq \ideala )} \cdot 
        (\Lambda _K ^{\ideala } - \LambdaSiegel ^{\ideala } ) }_{U^{s+1}([\pm N]^n_{\ida })}
        &\gg _{s,n} 
        \Bigl( \exp (-2(\log N)^{\mydecayexponent })  \Bigr) ^{(\log N)^{O_s(\mydecayexponent )}} 
        \\ &= \exp (- (\log N)^{O_s(\mydecayexponent )}).
    \end{align}
    Since we know that the left-hand side is smaller than $\exp (-(\log N)^{c_{s,n}})$ by Theorem \ref{thm:vonMang-Siegel-close!}, it follows $O_s(\mydecayexponent ) \ge c_{s,n}$.
    By contraposition, if $\mydecayexponent < c_{s,n}/O_s(1)$, then the inequality \eqref{eq:asymptotic-orthogonality-with-nilsequences} must hold.
\end{proof}

It turns out that Theorem \ref{thm:asymptotic-orthogonality-with-nilsequences} is actually equivalent to Theorem \ref{thm:vonMang-Siegel-close!} thanks to the Inverse Theorem of Gowers norms.

\begin{proof}[Proof that Theorem \ref{thm:asymptotic-orthogonality-with-nilsequences} implies Theorem \ref{thm:vonMang-Siegel-close!}]
    Let $0<c<1$ be a positive number to be specified later and suppose that the following inequality holds for some $N$:
    \begin{align}\label{eq:contradiction-assumption}
        \lnorm{
            1_{\Omega \cap (a+\idealq \ideala )} \cdot 
            (\Lambda_{K}^{\ideala } - \LambdaSiegel ^{\ideala } )
            }_{U^{s+1}([\pm N]^n_{\ida })} \ge \exp (-(\log N)^c) . 
    \end{align}
    By the Inverse Theorem \ref{thm:GowersInverse} of the Gowers norms, there exist a nilmanifold $G/\Gamma $ of dimension $\le (\log N)^{O_{s,n}(c )}$, 
    step and degree $\le s$ and 
    complexity $\le \exp ((\log N )^{O_{s,n}(c)})$, 
    a function $F\colon G/\Gamma \to \C $ with Lipschitz norm $\le \exp ((\log N )^{O_{s,n}(c)})$ and a polynomial map $g\colon \Z ^n \to G$ such that 
    \begin{align}
        \left\vert
               \Expec _{x\in [\pm N]^n_\ida }
               1_{\substack{\Omega \cap (a+\idealq \ideala )}}\cdot 
               (\Lambda _K ^{\ideala }(x) - \LambdaSiegel ^{\ideala }(x))
               F(g (x)\Gamma )
        \right\vert
        \ge
        \exp (-(\log N )^{O_{s,n}(c)}) .
    \end{align}
    Let $M>1$ be larger than all the $O_{s,n}(1)$ above.
    If $Mc \le \mydecayexponent _{s,n}$ and $N$ is large enough as in Theorem \ref{thm:asymptotic-orthogonality-with-nilsequences}, then the left-hand side must be smaller than $\exp (-(\log N)^{\mydecayexponent _{s,n}})$.
    It follows that $O_{s,n}(c) \ge \mydecayexponent _{s,n}\ge Mc $, contradicting $M>O_{s,n}(1)$.
    
    Therefore, for $c_{s,n}:=\mydecayexponent _{s,n}/M$, the assumption \eqref{eq:contradiction-assumption} cannot hold.
    That is, we have 
    \begin{align}
        \lnorm{
            1_{\Omega \cap (a+\idealq \ideala )} \cdot 
            (\Lambda_{K}^{\ideala } - \LambdaSiegel ^{\ideala } )
            }_{U^{s+1}([\pm N]^n_{\ida })} < \exp (-(\log N)^{c_{s,n}}) 
    \end{align}
    for all sufficiently large $N\gg _{s,K}1$.
\end{proof}


Having established the equivalence of Theorem \ref{thm:vonMang-Siegel-close!} and Theorem \ref{thm:asymptotic-orthogonality-with-nilsequences}, 
our goal is now to prove Theorem \ref{thm:asymptotic-orthogonality-with-nilsequences}.

\subsection{The step 0 case}\label{sec:step-0-case}

Suppose $\step (G)=0,$ i.e., $G$ is the trivial group. 
The map $g\colon \ideala \to G$ is the trivial map and thus the composite map $\ideala \xto g G \to G/\Gamma \xto F \C $ is a constant map whose value has size $\le 1$.
In this case, 
to prove Theorem \ref{thm:asymptotic-orthogonality-with-nilsequences}, it suffices to show that 
\begin{align}
    \frac{1}{N^n}\left\vert 
        \sum _{x\in \Omega \cap (a+\idealq \ideala )} 
        (\Lambda _K ^{\ideala }(x) - \LambdaSiegel ^{\ideala }(x))
    \right\vert 
\end{align}
exhibits a subpolynomial decay of the form $\exp (-(\log N)^\mydecayexponent )$ in $N$ uniformly in $\Omega ,a+\idq\ida $.

\ \ \  

If 
$\Nrm (\idq )< \exp ((\log N)^{1/200})$, 
then Proposition \ref{prop:essentially-Mitsui} (Mitsui's theorem) applies and gives a 
decay rate 
$\exp (-(\log N)^{1/201} )$. 

\ \ \ 

If 
$\Nrm (\idq )\ge \exp ((\log N)^{1/200})$, 
then the number of $x$'s involved in the sum is $\ll _K N^n/\Nrm (\idq )\le N^n \exp (-(\log N)^{1/200})$.
Since the size of the summand is $\ll _K \log N$ by Proposition \ref{prop:pointwise-bound}, we obtain the desired decay again with exponent 
$1/201$.

This proves the $s=0$ case of Theorem \ref{thm:asymptotic-orthogonality-with-nilsequences} with $\mydecayexponent _{0,k,n}:= 1/201$ say.

\leitfaden

\subsection{Induction on the step of nilpotence}\label{sec:induction-on-step}

Having established the base case $\step (G)=0$,
we turn to the induction part.

\begin{proposition}\label{prop:step-goes-down}
    Let $K$ be a number field of degree $n$. 
    Let $1\le s\le k$ be integers.
    Suppose we are given an exponent $0<\mydecayexponent <1/100$,
    positive integer $N\gg _{k,\rho , K}1$ 
    and data $\Omega ,a+\idq\ida $,
    $\ida\xto g G\surj G/\Gamma \xto F \C $
    satisfying the hypotheses of Theorem \ref{thm:asymptotic-orthogonality-with-nilsequences} with 
    step $\le s$, degree $\le k$ and 
    exponent $\mydecayexponent $
    but contradicting its conclusion:
    \begin{equation}\label{eq:negation-of-asymptotic-orthogonality}
        \left\vert
            \Expec _{x\in [\pm N]^n_\ida }
            1_{\substack{x\in \Omega \cap (a+\idealq \ideala )}}
            (\Lambda _K ^{\ideala }(x) - \LambdaSiegel ^{\ideala }(x))
            F(g (x)\Gamma )
            \right\vert
            >  
            \exp (-(\log N)^{\mydecayexponent }) .
        \end{equation}
    Then we can find new data $\Omega ',a'+\idq '\ida ,\ida \xto{g'}G'\surj G'/\Gamma '\xto{F'}\C $
    satisfying the hypotheses of Theorem \ref{thm:asymptotic-orthogonality-with-nilsequences}
    with 
    step $\le s-1$, degree $\le k$ and 
    exponent $O_{k,n}(\mydecayexponent )$
    that contradict its conclusion.
\end{proposition}

It is easy to see that Proposition \ref{prop:step-goes-down} implies Theorem \ref{thm:asymptotic-orthogonality-with-nilsequences}
because the $s=0$ case has already been settled in \S \ref{sec:step-0-case}.

The proof of Proposition \ref{prop:step-goes-down} will be finished 
in \S \ref{sec:Type-II}. 
See Figure \ref{fig:Leitfaden} for an overview of the proof.


\subsection{Easy reductions}\label{sec:easy-reductions}

First suppose $\Nrm (\idq )\ge \pseudopoly{3\mydecayexponent }$. By an argument similar to \S \ref{sec:step-0-case} the left-hand side of \eqref{eq:negation-of-asymptotic-orthogonality} would be bounded from above by 
$\pseudopoly{2 \mydecayexponent }$, a contradiction.
Therefore the conclusion of Theorem \ref{thm:asymptotic-orthogonality-with-nilsequences} would vacuously hold.
So we may assume $\Nrm (\idq )< \pseudopoly{3\mydecayexponent }$.

That said, by classical Fourier analysis (Proposition \ref{prop:Fourier} 
with 
$M:=\pseudopolydecay{10\mydecayexponent } $ and 
$Y:=M^2$), 
we can write the indicator function of $\Omega \cap (a+\idealq \ideala )$ in the form 
\begin{equation}\label{eq:indicator-function-Fourier-expansion}
    1_{\Omega \cap (a+\idealq \ideala )} 
    = 
    \left( \sum_{\theta \in \Theta } 
    c_\theta \cdot e^{2\pi i \theta (-)}  \right)
    +B+H,
\end{equation}
where 
$\Theta \subset \Hom (\ida ,\R )$ is a finite subset of size $< \pseudopoly{11\mydecayexponent } $,
$c_\theta \in \C $ has size $\le 1$ for every $\theta \in \Theta $, 
and the functions $B,H\colon [\pm N]^n_{\ida }\to \C $ satisfy 
\begin{itemize}
    \item 
    $\lnorm{B}_\infty \le \Nrm (\idq ) (< \pseudopoly{3\mydecayexponent })$, 
    \item 
    $\# \Supp (B)<\pseudopolydecay{9\mydecayexponent} N^n $ and 
    \item 
    $\lnorm{H}_\infty <\pseudopolydecay{9\mydecayexponent} $.
\end{itemize}

Substitute this into the left-hand side of \eqref{eq:negation-of-asymptotic-orthogonality}. 
The contribution of $B,H$ is negligible.
By the pigeonhole principle among the $\theta $'s, there is a frequency $\theta \in \Theta $ such that
\begin{equation}
    \left\vert 
        \sum _{x\in [\pm N]^n_{\ida }}
        (\Lambda_{K}^{\ideala } - \LambdaSiegel ^{\ideala }) (x) 
        \cdot e^{2\pi i \theta (x)} F(g(x)\Gamma )
    \right\vert 
    > \pseudopolydecay{12\mydecayexponent } N^n 
    .
\end{equation}
By the assumption $s\ge 1$, 
the function $x\mapsto e^{2\pi i \theta (x)} F(g(x)\Gamma )$ can be seen as a nilsequence involving the nilmanifold $(G/\Gamma )\times (\RZ )$ of step $s $, suject to the same type of complexity and \Lip norm bounds $<\exp ((\log N)^{12\mydecayexponent })   $.
Change the notation a little bit and write $F(g(x)\Gamma )$ for this new nilsequence.
This reduces us to the situation (without the indicator function $1_{\Omega \cap (a+\idq\ida )}$ and with exponent $12\mydecayexponent $)
\begin{equation}\label{eq:this-new-nilsequence}
    \left\vert 
        \sum _{x\in [\pm N]^n_{\ida }}
        (\Lambda_{K}^{\ideala } - \LambdaSiegel ^{\ideala }) (x) 
        \cdot F(g(x)\Gamma )
    \right\vert 
    > \pseudopolydecay{12\mydecayexponent } N^n 
    .
\end{equation}

Take the Fourier expansion of $F$ along $G_{(s)}$ in Proposition \ref{prop:Fourier-along-Gs} 
with $\delta := \exp (-(\log N)^{20\mydecayexponent })  $:
%
%
%
\begin{equation} 
    F=\left( \sum _{|\xi |\le \mypseudopolygrowth } F_\xi  \right) +E
    , 
\end{equation}
where $\xi \colon G_{(s)}\to \RZ $ are continuous homomorphisms, $F_\xi $ satisfies 
$\lnorm {F_\xi} _\infty + \mLip (F_\xi ) \le \mypseudopolygrowth $
and $E$ satisfies $\lnorm E _\infty < \exp (-(\log N)^{20\mydecayexponent }) $.

The contribution of $E$ to the left-hand side of \eqref{eq:this-new-nilsequence} is negligible. 
By the pigeonhole principle applied to the remaining part of \eqref{eq:this-new-nilsequence}, there is a vertical frequency $\xi \colon G_{(\nilcounter )}\to \RZ $ with $|\xi |\le \mypseudopolygrowth $ such that
\begin{equation}\label{eq:with-vertical-frequency-xi}
    \left\vert
        \sum  _{x\in [\pm N]^n_{\ida } } 
        (\Lambda_{K}^{\ideala } (x) - \LambdaSiegel ^{\ideala } (x) )
        F_\xi (g(x)\Gamma )
        \right\vert
        \ge \mypseudopolydecay N^n 
        .
    \end{equation} 
If this $\xi $ happens to be 
$\xi =0$, 
then $F_\xi $ can be seen as a function on the nilmanifold $\frac{G/G_{(\nilcounter )}}{\Gamma / (G_{(\nilcounter )} \cap \Gamma)}$ of step $\le s-1$ with the naturally induced filtration and \Malcev basis.
This achieves the assertion of Proposition \ref{prop:step-goes-down} in the case $\xi =0$.

From now on, assume $\xi \ne 0$ in \eqref{eq:with-vertical-frequency-xi}. 
In particular, $F_\xi $ has average $0$.

\subsection{Applying Vaughan decomposition}\label{sec:applying-Vaughan}

By the Vaughan decomposition (Proposition \ref{prop:Vaughan-decomposition}) of $\Lambda _K$ and $\LambdaSiegel $, 
we see that there is at least one (twisted) Type I or Type II sum $S(-)$ such that the next inequality holds:
\begin{equation}\label{eq:with-type-I-or-II}
    \left\vert
        \sum _{x\in [\pm N]^n_{\ida } } 
        S (x\ideala \inv )  
        F_\xi (g(x)\Gamma )
        \right\vert
        \ge \mypseudopolydecay N^n 
        .
\end{equation}
Write $\Phi (x):= F_\xi (g(x)\Gamma)$.
For the moment, what counts is that $\Phi $ satisfies $\lnorm{\Phi }_\infty \le \mypseudopolygrowth   $.

We are going to do case-by-case analysis depending on whether the above $S(-)$ is of (twisted) Type I or Type II.
Before that, however, we need a separate treatment when $\Nrm (\idealqSiegel )$ is too large.

\subsection{The case with large $\Nrm (\idqSiegel )$}\label{sec:case-of-large-q-Siegel}

Let us again be in the situation \eqref{eq:with-vertical-frequency-xi}.
Let $C>C'>1$ be two positive numbers at least three times larger than every $O_{k,n}(1)$ constant that appears from now to \eqref{eq:in-the-realm}
and assume $\Nrm (\idqSiegel )> \exp ((\log N)^{C\mydecayexponent })$.
By Proposition \ref{prop:Cramer-Siegel-close}, we know 
\begin{align}
    \lnorm{
        \LambdaSiegel ^{\ideala } - \LambdaCramer ^{\ideala } 
        }_{U^{s+1}([\pm N]^n_{\ida })}
    \ll_{s,K} 
    \Nrm (\idealqSiegel )^{-1/2^s}
     < \exp (-(\log N   )^{C'\mydecayexponent } ).
\end{align}
With this Gowers norm bound, Proposition \ref{prop:converse-to-inverse-theorem} (converse to the inverse theorem of Gowers norms) implies
\begin{multline}
    \left\vert 
        \sum _{x\in [\pm N]^n_{\ida } } 
        (\LambdaSiegel ^{\ideala } (x) - \LambdaCramer ^{\ideala } (x) F_\xi (g(x)) )
    \right\vert 
    <
    N^n
    \exp \left( -(\log N )^{C'\mydecayexponent /2 }  \right)
    .
\end{multline}
Note that the implied constant in the previous ``$\ll _{s,K}$'' has been absorbed in the decay factor by assuming $N\gg_{k,n,\mydecayexponent }1$.
This kind of procedure will not be mentioned henceforth.

Combine this inequality with \eqref{eq:with-vertical-frequency-xi} to obtain 
\begin{equation}
        \left\vert 
        \sum _{x\in [\pm N]^n_{\ida } } 
        (\Lambda_K ^{\ideala } (x) - \LambdaCramer ^{\ideala } (x) F_\xi (g(x)) )
    \right\vert 
    >
    \mypseudopolydecay N^n
    .
\end{equation}
Note that now we have the \Cramer model $\LambdaCramer ^{\ida }$ in the formula in place of the Siegel model $\LambdaSiegel ^{\ida }$.

Perform the Vaughan decomposition and get an inequality of the form \eqref{eq:with-type-I-or-II} again.  
Since the Vaughan decomopsitions of $\Lambda_K $ and $\LambdaCramer $ do not involve twisted Type I sums, the sum $S(-)$ in this case is either a Type I or Type II sum.
 
We conclude that the proof of Proposition \ref{prop:step-goes-down} in the realm
\begin{equation}\label{eq:in-the-realm}
    \Nrm (\idqSiegel ) > \exp ((\log N   )^{C\mydecayexponent })
\end{equation}
where $C=O_{k,n}(1)$ is large enough as above, is reduced to the case of (not twisted) Type I/II sums.

\section{(Twisted) Type I case}\label{sec:Type-I} 

By \S\S \ref{sec:applying-Vaughan}, \ref{sec:case-of-large-q-Siegel},
we are reduced to the situation \eqref{eq:with-type-I-or-II} where $S(-)$ is a (twisted) Type I/II sum, and when it is a twisted Type I sum we may assume 
\[\Nrm (\idealqSiegel )< \exp ((\log N)^{O_{k,n}(\mydecayexponent )}).\]
In this section we treat the (twisted) Type I sum.

Recall from \S \ref{sec:Vaughan} that a (twisted) Type I sum $S\colon \Ideals _{K,\le O_K(N^n)}\to \R $
is written as 
\[ S(\idealc )= 1 [\Nrm (\idealc )\le N']\cdot \psiSiegel (\idealc ) 
\sum _{ \substack{ \idealb \divides \idealc \\ \Nrm (\idealb )\le N^{2n/3} } } 
A(\idealb )  ,\] 
where $N'\le O_K(N^n)$, $|A(\idealb )|\le \log N$ for every ideal $\idealb $, and $\psiSiegel $ is the $Q$-Siegel character if it exists, and the constant function $1$ otherwise.

\subsection{\CauSch and manipulation of the range of sum}

After reordering the sum in \eqref{eq:with-type-I-or-II} we apply the triangle inequality. Also letting $\log N$ be absorbed into $\mypseudopolydecay $, we obtain (recall $\Phi (x):=F_\xi (g(x)\Gamma )$)
\begin{equation}
    \sum _{
        \substack{\idealb \\ \Nrm (\idealb )\le N^{2n/3}} }
        \frac 1{\Nrm (\idealb )^{1/2}}
        \cdot 
        \left\vert 
            \sum _{
                \substack{x\in [\pm N]^n_{\ida }  
                \cap \ideala \idealb 
                \\ 
                \Nrm (x\ideala \inv ) \le N' 
                }
            }
            \psiSiegel ([x\ideala \inv ])
            \Phi (x)
        \right\vert 
        >
        N^n\mypseudopolydecay   .
\end{equation}
Apply the Cauchy--Schwarz inequality for the sum over $\idealb $ and observe that the first factor becomes $\sum _{\idealb }\frac 1{\Nrm (\idealb )} =O_K(\log N) < \mypseudopolygrowth   $ by Lemma \ref{lem:divisor-moment}. This gives 
\begin{align}\label{eq:after-CS-in-type-I}
    \sum _{
        \substack{\idealb \\ \Nrm (\idealb )\le N^{2n/3}} }
        {\Nrm (\idealb )} 
        \cdot 
        \sum _{
            \substack{x\in [\pm N]^n_{\ida }  
            \ideala \idealb 
            \\ 
            \Nrm (x\ideala \inv ) \le N' 
            }
        }
        \left\vert 
            \psiSiegel ([x\ideala \inv ])
            \Phi (x)
        \right\vert ^2
        &>
        N^{2n}\mypseudopolydecay 
        .
\end{align}

Apply the pigeonhole principle to the dyadic decomposition 
\[ [1,N^{2n/3}]\subset [1,2]\cup [2,4]\cup [4,8]\cup \dots \] 
(consisting of $\ll _n \log N$ subintervals) of the range of $\Nrm (\idealb )$.
We find a $1\le D \le N^{2n/3}$ such that 
\[ 
\sum _{
        \substack{\idealb \\ D\le \Nrm (\idealb )\le 2D} }
        D
        \cdot 
        \sum _{
            \substack{x\in [\pm N]^n_{\ida }
            \cap \ideala \idealb 
            \\ 
            \Nrm (x\ideala \inv ) \le N' 
            }
        }
        \left\vert 
            \psiSiegel ([x\ideala \inv ])
            \Phi (x)
        \right\vert ^2
        >
        N^{2n}\mypseudopolydecay  . 
\]

There are $\ll _K D$ ideals $\idealb $ involved and the inner sum is $\ll _n (N^n/D)^2$ for each $\idealb $.
It follows that (enlarging the $O_{k,n}(1)$ constant as always)
there are $\ge \mypseudopolydecay  D$ many $\idealb $'s with $D\le \Nrm (\idealb )\le 2D $ such that

\[ 
\left\vert 
\sum _{
            \substack{x\in [\pm N]^n_{\ida }  
            \cap \ideala \idealb 
            \\ 
            \Nrm (x\ideala \inv ) \le N' 
            }
        }
            \psiSiegel ([x\ideala \inv ])
            \Phi (x)
        \right\vert
        \gg _K \mypseudopolydecay N^n /D . 
\]

For each of such $\idealb $, 
by the pigeonhole principle for the $\Nrm (\idqSiegel )$ different mod-$\idealqSiegel \ida\idb $ congruence classes, 
we can find a congruence class $a_\idealb + \idqSiegel\ideala\idb $ 
such that the following inequality holds:
\[ 
\left\vert 
\sum _{
            \substack{x\in [\pm N]^n_{\ida } 
            \cap (a_{\idb }+\idqSiegel \ideala \idealb )
            \\ 
            \Nrm (x\ideala \inv ) \le N' 
            }
        }
            \Phi (x)
        \right\vert
        \gg _K \mypseudopolydecay N^n /D . 
\]
Equip $\idqSiegel \ida\idb $
with a norm-length compatible basis. 
It defines a fundamental domain for the quotient $\ida / (\idqSiegel \ida\idb )$ of $\ida $.
Choose the lift $a_\idealb \in \ida $ 
of the congruence class $a_\idealb + \idqSiegel\ideala\idb $ from this fundamental domain.
Then we have inclusion 
\[ 
    [\pm N]^n_{\ida } \cap (a_\idealb + \idqSiegel \ideala \idealb )
    \subset a_\idealb + 
    \left[ 
        \pm O_K(1)
        \frac{N}
        { \Nrm (\idqSiegel \idb )^{1/n} 
        }  
    \right] _{\idqSiegel\ida\idb }  
\]
Write 
$N_\idealb 
:= 
\frac{N}
{\Nrm  (\idqSiegel\idb )^{1/n} 
}$. 
The above inequality becomes
\[ 
\left\vert 
\sum _{
            \substack{x\in  ([\pm N]^n_{\ida }-a_{\idb }) 
            \cap \left[ 
    \pm O_K(1){N_\idealb }  
\right]^n _{\idqSiegel\ida\idb } 
            \\ 
            \Nrm ((x+a_\idealb )\ideala \inv ) \le N' 
            }
        }
            \Phi (x+a_\idealb )
        \right\vert
        \gg _K \mypseudopolydecay N^n /D  . 
\] 
Cover the cube $\left[ 
    \pm O_K(1){N_\idealb }  
\right]^n _{\idealqSiegel\ida\idb } 
$ 
by translates of the smaller cube of side length $\mypseudopolydecay  N_\idealb $
$= \mypseudopolydecay N/D^{1/n}$.
By \S \ref{appendix:lattice-points}, the contribution from the smaller cubes intersecting the boundary of $([\pm N]^n_{\ida } -a_\idealb )\cap \{ x \mid \Nrm ((x+a_\idealb)\ideala\inv )\le N' \} $ is negligible.
By the pigeonhole principle applied to the $\mypseudopolygrowth $ many cubes, we can find a translate of the small cube represented by an $a_\idealb ' \in \left[\pm O_K(1)\frac{N}{D^{1/n}}\right] _{\ideala \idealb }$ such that 
\[  
\left\vert 
    \sum _{x\in \left[ 1, \mypseudopolydecay  N_\idealb\right]^n 
    _{\idealqSiegel\ida\idb } 
    }
    \Phi (x+a_\idealb ') 
\right\vert
\gg _K \mypseudopolydecay \frac{N^n}{D}  . \] 
Since $\# \left[ 1, \mypseudopolydecay  N_\idealb\right]^n 
= 
\mypseudopolydecay \frac{N^n}{
    \Nrm (\idqSiegel\idb ) 
}
< \frac{N^n}{D} 
$, 
in terms of average this says 
\begin{equation}\label{eq:in-terms-of-average}    
    \left\vert 
        \Expec _{x\in \left[ 1, \mypseudopolydecay  N_\idealb\right]^n
        _{\idqSiegel \ida\idb }
        }
        \Phi (x+a_\idealb ') 
        \right\vert
        \gg _K \mypseudopolydecay   .
    \end{equation} 
%
%
%
%

\subsection{Applying the equidistribution theorem}

Now we recall that $\Phi (x+a'_\idealb )=F_\xi (g(x+a'_\idealb )\Gamma )$ and apply Theorem \ref{thm:Leibman} to inequality \eqref{eq:in-terms-of-average} with the configuration
%
%
%
%
\[ \Z ^n \leadsto  
\idqSiegel \ideala \idealb,
\qquad  
\delta \leadsto \mypseudopolydecay ,
\qquad 
N \leadsto \mypseudopolydecay  N_\idealb .
\]
This gives some integer $r_\idealb \in \{ 1,\dots ,\dim (G/[G,G])\}$ 
and horizontal characters 
\[
\eta _{\etacounter ,\idealb}\colon G\to \R 
\] 
($\etacounter =1,\dots ,r_\idealb $) of sizes $\le \mypseudopolygrowth $ 
such that 
\[
\lnorm{\eta _{\etacounter ,\idealb}\circ g ((-)+a'_\idealb )}
_{\Cinfty [1,\mypseudopolydecay N_\idealb ]^n
_{\idqSiegel \ideala \idealb}
} 
\le \mypseudopolygrowth  
\]
for each $\etacounter =1,\dots ,r_\idealb $
and that if we define $G'_\idealb := \bigcap _{\etacounter =1}^{r_\idealb} \ker (\eta _{\etacounter ,\idealb})$ (which is $\mypseudopolygrowth   $-rational) then
$(G'_{\idealb})_{(\nilcounter )}$ is contained in $\ker \xi $. 

Now, one is tempted to consider the associated factorization $g=\varepsilon \cdot g_1\cdot \gamma $ given by Proposition \ref{prop:factorization}.
But it is too early.
It would only give us a factorization of $g$ as a polynomial map defined on the ideal 
$\idqSiegel \ideala \idealb$.
This is too sparse in $\ideala $.

Proposition \ref{prop:extrapolation} (extrapolation) applied to the affine-linear map 
$ \idqSiegel \ideala \idealb 
\to \ideala \idealb $ 
given by $x\mapsto x+a'_\idealb $ shows that there is an integer $1\le q_{\etacounter ,\idealb} \ll _K \Nrm (\idealqSiegel )^{O_{k,n}(1)} $ for each $\etacounter =1,\dots ,r_\idealb $ such that 
\begin{align} 
    \lnorm{(q_{\etacounter ,\idealb}\cdot \eta _{\etacounter ,\idealb}) \circ g}_{\Cinfty [1,\mypseudopolydecay N_\idealb ]^n_{\ideala\idealb }} 
    &\le 
    \mypseudopolygrowth    
    \Nrm (\idealqSiegel )^{O_{k,n}(1)}
    \\
    &= 
    \mypseudopolygrowth    
    .
\end{align}
By the definition of the smoothness norm (Definition \ref{def:smoothness-norm}), the same form of inequality holds if we replace the side length $\mypseudopolydecay N_{\idb }$ of the cube by $N/D^{1/n}$.

The horizontal characters $q_{\etacounter ,\idealb}\cdot \eta _{\etacounter ,\idb} $ have sizes $< \mypseudopolygrowth    $.
The number of horizontal characters obeying this bound is 
$< \mypseudopolygrowth   $. 
By the pigeonhole principle, 
we can pick $> \mypseudopolydecay D$ many ideals $\idealb $ with norm $\asymp D$ among the previously chosen
$\mypseudopolydecay D$ many so that 
we may assume $r:= r_\idealb $ is a common value for all such $\idealb $ and also
$\eta _\etacounter  := q_{\etacounter ,\idealb}\cdot \eta _{\etacounter ,\idealb}$ is common for each $\etacounter =1,\dots ,r$.
Also, since $\Cl (K)$ is finite, we may also assume that the ideals $\idealb $ chosen above belong to a common ideal class.
This brings us to the following situation.

\begin{situation}\label{situation:this-brings-us}
    In addition to the situation in \S \ref{sec:induction-on-step}, we have the following additional data: 
    \begin{itemize}
        \item 
        an integer $1\le r\le (\log N)^{\mydecayexponent } $,
        \item 
        horizontal characters $\eta _1,\dots ,\eta _r \colon G\to \R $ of sizes $\le \mypseudopolygrowth   $,
        \item 
        an integer $1\le D\le N^{2n/3}$,
        \item 
        at least $ \mypseudopolydecay  D$ many ideals 
        $\idealb $ with norms $\in [D,2D]$ in a common ideal class.
    \end{itemize}
    For each choice of $\etacounter \in \{ 1,\dots ,r \}$ and $\idealb $ appearing above, we have
    \[ \lnorm{\eta _{\etacounter }\circ (g|_{\ida\idb })}_{\Cinfty [1, N/D^{1/n} ]^n_{\ideala\idealb }} 
    < \mypseudopolygrowth    . \]
    Also, if we write 
    $G':= \bigcap _{\etacounter =1}^r 
    \ker (\eta _\etacounter ) $, 
    then $\xi $ is constant on $(G')_{(\nilcounter )}$.
\end{situation}

\subsection{Utilizing the parametrization of ideals}

Let $\idealb _0\subset \OK $ be a representative of the common ideal class of the chosen ideals $\idealb $.
Also choose a norm-length compatible fundamental domain $\calD \subset (K\otimes _{\Q}\R )\baci $ for the action of $\OK\baci $.
Then there is a bijection 
 \begin{equation}
    \begin{array}{ccc}
        \calD \cap \idealb_0\inv 
        &\cong 
        &\{ \idealb\subset \OK \mid [\idealb ]=[\idealb_0 ] \text{ in }\Cl (K) \}
        \\ 
        b & \mapsto & b\idealb_0
    \end{array}    
\end{equation} 

Choose a norm-length compatible basis of $\ideala \idealb_0 $
and transfer it to every $\ideala\idealb $ via the isomorphism
$\times b \colon \ideala\idealb_0 \isoto \ideala\idealb $ where $b$ corresponds to $\idealb $ in the above correspondence.
For each $\etacounter \in \{ 1,\dots ,r \}$, let $H^{\etacounter }$ be the composite map 

\[ H^{\etacounter }\colon \ideala\idealb_0 \times \idealb_0\inv \xto{(y,z)\mapsto yz} \ideala \xto g G\xto{ \eta _{\etacounter }}\R .  \] 

Write $\alpha =(\alpha _1,\dots ,\alpha_n)$ and $b=(b_1,\dots ,b_n)$ for the chosen bases.
For multi-indices $\undli =(i_1,\dots ,i_n)\in\N^n $, write $\alpha ^{\undli }:= \prod _{\nu =1}^n \alpha_\nu ^{i_\nu}\in \Z $.
The map $H^{\etacounter }$ can be written using coefficients 
\[ 
H^{\etacounter } (\alpha ,b)=
    \sum _{
        \undli ,\undlj \in \N^n
    }
h^{\etacounter }_{\undli ,\undlj } \alpha ^{\undli }b^{\undlj }, 
\quad
h^{\etacounter }_{\undli ,\undlj }\in \R . 
\]
Let $|\undli |:= i_1+\dots +i_n$. Note that $h^{\etacounter }_{\undli ,\undlj }\neq 0$ only when $|\undli |=|\undlj |$ because the map $\ideala\idealb _0 \times \idealb_0\inv \to \ideala $ is bilinear.
For each $b\in \idealb_0\inv $, the map 
\[ 
H^{\etacounter }(-,b)
\colon 
\ideala\idealb_0 \to \R, 
\quad 
\alpha \mapsto H^{\etacounter }(\alpha ,b) 
\]
is a polynomial map of degree $\le s$. 
Its $\undli \ordinalth $ coefficient is 
$\sum _{    
    \substack{\undlj \in \N^n \\ |\undlj |=|\undli | } }
    h_{\undli ,\undlj }^\etacounter b^{\undlj }
$.
Via the parametrization of ideals, Situation \ref{situation:this-brings-us} implies that for 
$> \mypseudopolydecay   D$ values of $b\in [\pm O_K(1)D^{1/n}]^n_{\idealb_0\inv } \cap \calD $, the polynomial map $H^{\etacounter }(-,b)$ satisfies 

\[
    \lnorm{H^{\etacounter }(-,b)}_{\Cinfty [1, N/D^{1/n} ]^n_{\ideala \idealb_0}}
    < \mypseudopolygrowth   .
\]

By the definition of $\lnorm{-}_{\Cinfty }$, this inequality is equivalent to saying that for each multi-degree $\undli \neq 0$, we have 
\begin{align} 
    \lnorm{\sum _{
            \substack{\undlj \in \N ^n \\ |\undlj |=|\undli |}
        } 
        h_{\undli ,\undlj }^\etacounter z^{\undlj }
        }_{\RZ }
    \le 
    \mypseudopolygrowth   
    \frac{D^{|\undli |/n}}{N^{|\undli |}} 
    .
\end{align} 

Let $H^{\etacounter }_{\undli }$ be the polynomial map 
$b\mapsto \sum _{
            \substack{\undlj \in \N ^n \\ |\undlj |=|\undli |}
        } 
        h_{\undli ,\undlj }^\etacounter b^{\undlj }
$.
Its $\undlj \ordinalth $ coefficient is $h_{\undli ,\undlj }^\etacounter $.
Since this inequality holds for $> \mypseudopolydecay  D $
values of $b$,
Proposition \ref{prop:Vinogradov} (Vinogradov lemma)\footnote{This is where we used the fact that there are enough number of $\idealb$!} 
gives an integer 
$0< q_{\undli }^\etacounter <\mypseudopolygrowth   $ such that 
\[ 
    \lnorm{
            q_{\undli }^\etacounter H^{\etacounter }_{\undli } 
        }
        _{\Cinfty [\pm D^{1/n}]^n_{\idealb_0}}
    <
    \mypseudopolygrowth   
    \frac{D^{|\undli |/n}}{N^{|\undli |}} .
\]  
Since the number of $\undli $ under consideration is $O_{k,n}(1)$, 
the same form of inequalities stay valid if we replace 
$q_{\undli }^\etacounter $'s 
by their least common multiple $q_{\etacounter }$.
By the definition of $\lnorm{ - }_{\Cinfty }$ we conclude 
\begin{equation}\label{eq:bound-coeff-hij}    
    \lnorm{q_{{\etacounter }} h_{\undli \undlj }^\etacounter}_{\RZ } 
    <
    \mypseudopolygrowth   
    \frac{1}{N^{|\undli |}} 
    \quad 
    \text{ for all }\undli ,\undlj \neq 0 .
\end{equation} 
(Recall that $h_{\undli ,\undlj }^{\etacounter}\in \R $ are coefficients of the map 
$H^{\etacounter }\colon \ideala\idealb_0 
\times 
\idealb_0\inv \to \ideala \to G\to \R $ 
and 
$h_{\undli ,\undlj }^{\etacounter}\neq 0$ only when $|\undli |=|\undlj |$.)

Let us specialize to $1\in \idealb_0\inv $.
This corresponds to considering the map 
\[ 
    q_{\etacounter } \eta _\etacounter \circ (g|_{\ideala\idealb_0 }) \colon \ideala\idealb_0 \subset \ideala \xto g G \xto{q_{\etacounter }\cdot \eta _\etacounter} \R .
\]
Because $1$ is represented by a vector with entries of sizes $O_K(1)$ with respect to the chosen basis of $\idealb_0$, 
the bound \eqref{eq:bound-coeff-hij} implies that the $\undli \ordinalth $ coefficient of 
$q_{\etacounter }\eta _\etacounter \circ (g|_{\ideala\idealb_0 }) $ has 
$\lnorm{-}_{\RZ }$-norm $< \mypseudopolygrowth   
    / {N^{|\undli |}} $ 
    for each $|\undli |$.
Equivalently, 
\[ 
    \lnorm{q_{\etacounter }\eta _\etacounter \circ (g|_{\ideala\idealb_0 })}_{\Cinfty [1,N]^n_{\ideala \idealb_0}} 
    <
    \mypseudopolygrowth   
    .
\] 
By Proposition \ref{prop:extrapolation} (extrapolation) applied to the inclusion $\ideala\idealb_0\inj \ideala $, up to modifying 
$q_{\etacounter}$ by a factor $O_{K,s}(1)$ this implies 
\[ 
    \lnorm{q_{\etacounter }\eta _\etacounter \circ g}_{\Cinfty [1,N]^n_{\ideala }} 
    <
    \mypseudopolygrowth   
    .
    \] 
\subsection{Taking the factorization associated to the horizontal characters}
\label{sec:taking-factorization-associated-to-Leng}

Now take the factorization $g=\varepsilon \cdot g_1 \cdot \gamma $ 
from Proposition \ref{prop:factorization}
associated to these horizontal characters $q _{1}\eta _1,\dots ,q _{r}\eta _r$:

\begin{itemize}
    \item $\varepsilon $ is $\Bigl( \mypseudopolygrowth   , [1,N]^n \Bigr) $-smooth;
    \item $g_1 $ takes values in $G_1 :=  \ker (\eta _1)\cap \dots \cap \ker ( \eta _r)$.
    Recall $G_1/(G_1\cap \ker (\xi ))$ with the induced filtration is a nilmanifold of step $<s$ by 
    Situation \ref{situation:this-brings-us};
    \item $\gamma $ is $\mypseudopolygrowth   $-rational,
\end{itemize}
and plug them into inequality \eqref{eq:with-vertical-frequency-xi}. This reads 
\begin{equation}\label{eq:after-factorization}
    \left\vert
        \sum _{x\in [\pm N]^n_{\ida } } 
        (\Lambda_{K}^{\ideala } (x) - \LambdaSiegel ^{\ideala } (x) )
        \cdot 
        F_\xi \Bigl( \varepsilon (x) g_1 (x)\gamma (x) \cdot \Gamma \Bigr)
        \right\vert
        \ge \mypseudopolydecay N^n 
        .
    \end{equation}
Let $1\le q_\gamma \le \mypseudopolygrowth   $ be the period of $\gamma (-)$.
By the pigeonhole principle for different mod $q_{\gamma }\ida $ classes, 
we obtain a version of \eqref{eq:after-factorization}
where the sum is restricted to a single congruence class $a'+q_{\gamma }\ida $. 
The value
$\gamma ':= \gamma (x)$ does not depend on $x\in a'+q_{\gamma }\ida $. 
Define the ideal $\idq '\subset \OK $ by
\[ \idq ':= q_\gamma \OK . \] 

If we let $\{ \gamma ' \}$ be the fractional part of $\gamma '$ (which lies in the fundamental domain of $G$ mod $\Gamma $ determined by the \Malcev basis), 
then $\gamma '\Gamma = \{ \gamma ' \}\Gamma $.
Therefore we may assume $\gamma '$ lies in this fundamental domain (this will be implicitly used for the complexity management).


Next, we chop the box $[\pm N]^n_{\ida }$ 
along coordinate hyperplanes into $\mypseudopolygrowth   $ different roughly equally sized small 
boxes $\Omega '$
(where we set the implied constant in $\mypseudopolydecay $ larger than before).
By the pigeonhole principle, our variant of the inequality \eqref{eq:after-factorization} remains valid with $[\pm N]^n_{\ida }$ 
replaced by one of these small convex sets $\Omega '$.
Choose a point $x_0'\in \Omega '\cap (a'+\idq '\ida )$ 
and write $\varepsilon _0':= \varepsilon (x_0')\in G $.
By the \Lip properties of $F_\xi $ and $\varepsilon (-)$, we have for all $x\in \Omega '$
\[
    F_\xi \Bigl( \varepsilon (x) g_1 (x) \gamma '\Gamma \Bigr)
    =
    F_\xi \Bigl( \varepsilon _0' g_1 (x) \gamma '\Gamma \Bigr)
    +
    \mypseudopolygrowth   ,
\]
where this $\mypseudopolydecay $ has larger implied constans than in \eqref{eq:after-factorization}.
It follows that our variant of \eqref{eq:after-factorization} remains valid with $\varepsilon (x)$ replaced by the constant $\varepsilon _0'\in G$. Namely:
\begin{equation}\label{eq:after-replacing-epsilon-by-constant}
    \left\vert
        \sum _{x\in \Omega ' \cap (a'+\idealq '\ideala )}
        (\Lambda_{K}^{\ideala } (x) - \LambdaSiegel ^{\ideala } (x) )
        \cdot 
        F_\xi \Bigl( \varepsilon _0' g_1 (x)\gamma ' \cdot \Gamma \Bigr)
        \right\vert
        \ge \mypseudopolydecay N^n 
        ,
    \end{equation}
with a different implied constant in $\mypseudopolydecay $ as always.

Consider the function $F':= F_\xi \Bigl( \varepsilon _0' \gamma ' \cdot (-)\cdot \Gamma \Bigr) \colon G/\Gamma \to \C $.
It still has vertical frequency $\xi $ (recall $G_{(\nilcounter )}$ is in the center of $G$),
and defines a function which we denote by the same symbol:
\[ F'\colon \frac{G/\ker \xi }{\Gamma /(\Gamma \cap \ker \xi ) } \to \C  .\] 
$F'$ has \Lip norm $\le \mypseudopolygrowth   $ with respect to the metric associated to the induced \Malcev basis.

Define the polynomial map $g'\colon \ida \to G/\ker \xi $ by 
\[ g'(x):= (\gamma ')\inv \cdot  g_1 (x)\cdot \gamma ' \mod \ker \xi . \]
This has values in the subgroup 
\begin{align}
    G' &:= (\gamma ')\inv G_1 \gamma ' \Big/ ( (\gamma ')\inv G_1 \gamma ' \cap \ker \xi ) ,
\end{align} 
which is conjugate to $G_1/(G_1\cap \ker \xi )$ by $\gamma '$ inside $G/\ker \xi $.
(Recall again that $\ker \xi \subset G_{(\nilcounter )}$ are in the center of $G$.)
Let $\Gamma ':= \Gamma \cap ((\gamma ')\inv G_1 \gamma ')/ (\Gamma \cap (((\gamma ')\inv G_1 \gamma '))\cap \ker \xi )$, which is a lattice in $G'$.
Then $G'/\Gamma ' $ with the induced filtration and \Malcev basis is a nilmanifold of step $<s$ (recall that this had been ensured by Leng's theorem \ref{thm:Leibman}) and complexity $\le \mypseudopolygrowth   $.

To summarize, we have arrived at the following situation. We have a polynomial map and a \Lip function involving the nilmanifold $G'/\Gamma '$ of complexity $\le \mypseudopolygrowth   $, 
$\dim (G)<(\log N   )^{\mydecayexponent }$:
\[ \ida \xto{ g'} G' \surj G'/\Gamma ' \xto{F'} \C  ,\quad \step (G') \le \nilcounter -1  ,\quad \lnorm{F'}_{\mLip } \le \mypseudopolygrowth   . \]
We also have a convex set $\Omega '\subset [\pm N]^n_{\ida } $ and a congruence class $a'+\idealq '\ideala $ such that the inequality 
\[ \left\vert 
    \sum _{x\in \Omega ' \cap (a'+\idealq '\ideala )}
    (\Lambda_{K}^{\ideala } (x) - \LambdaSiegel ^{\ideala } (x) )
    \cdot 
    F' \Bigl( g' (x) \cdot \Gamma ' \Bigr)
\right\vert 
\ge \mypseudopolydecay N^n
\]
holds. This establishes Proposition \ref{prop:step-goes-down}
in the (twisted) Type I case.

\ \ \ 

Before moving on to the Type II case, let us formulate part of our argument in \S \ref{sec:taking-factorization-associated-to-Leng} as a separate lemma for later use.
\begin{lemma}\label{lem:our-argument}
   In addition to the hypotheses of Proposition \ref{prop:step-goes-down},
    suppose that we are further given an integer $1\le r \le (\log N)^{\mydecayexponent }$, horizontal characters $\eta ^\etacounter \colon G\to \R $ for $\etacounter =1,\dots ,r $ of sizes $<\mypseudopolygrowth   $ such that
    \[ 
        \left( 
            \bigcap _{\etacounter =1}^r \ker (\eta^\etacounter )
        \right) 
        _{(\nilcounter )}\subset \ker (\xi )
    \] 
    and
    \begin{equation}
        \lnorm{\eta ^\etacounter \circ g }_{\Cinfty [1,N_{\ida }]^n_{\ideala}}
        < \mypseudopolygrowth   \qquad \te{for all } \etacounter =1,\dots ,r .
    \end{equation}
    %
    Then the conclusion of Proposition \ref{prop:step-goes-down} holds.
\end{lemma}

\section{Type II case}\label{sec:Type-II}

Next we consider the case where $S(-)$ in \eqref{eq:with-type-I-or-II} is a Type II sum.
This means 
\begin{equation}
    \left\vert 
    \sum _{x\in [\pm N]^n_{\ida } }
    \sum _{
        \substack{\idc \idd =x \ida \inv
        \\ 
        \Nrm (\idc ),\Nrm (\idd ) \ge N^{n/3} 
        }}
    L(\idc )T(\idd ) \Phi (x)
    \right\vert
    \ge 
    \mypseudopolydecay N^n , 
\end{equation}
where $L,T\colon \Ideals \to \C $ are functions satisfying $|L(\idc )|\le \log N$ and $|T(\idd )|\le \tau (\idd )$ for all $\idc ,\idd $ of norm $\ll _K N^n $, 
and we recall $\Phi (x):=F(g(x)\Gamma )$.

\subsection{\CauSch inequalities}

Eliminate $\idc $ by $\idc =x\ida \inv \idd \inv $ and reorganize the sum to obtain
(where we omit to write the conditions on $\idd ,x$ repeatedly):
\begin{align}
    &\mypseudopolydecay N^n 
    \\ 
    &\le 
    \left\vert 
    \sum _{\substack{
        \idd \subset \OK \\
        N^{n/3}\le \Nrm (\idd ) \ll _K N^{2n/3}
        }}
        \left(
        \sum _{
            \substack{
                x\in [\pm N]^n_{\ida }\cap \ida\idd \\ 
                \Nrm (x\ida\inv )\ge N^{n/3}\Nrm (\idd )
            }
        }
        T(\idd )
    L(x\ida \inv \idd \inv )
    \Phi (x)
    \right)
    \right\vert
    \\[10pt] 
    &=
    \left\vert
        \sum _{\substack{
        \idd 
        }}
        \frac{T(\idd )}{\Nrm (\idd )^{1/2}}
        \cdot 
        \left(
        \Nrm (\idd )^{1/2}
        \sum _{
            \substack{
                x 
            }
        }
        L(x\ida \inv \idd \inv )
        \Phi (x)
        \right)
    \right\vert
    .
\end{align}
By the \CauSch inequality for the sum $\sum _{\idd }$, this is bounded by 
\begin{equation}
    \le \left( \sum _{\idd } \frac{\tau (\idd )^2}{\Nrm (\idd )} \right)^{1/2}
    \cdot 
    \left( \sum _{\idd } \Nrm (\idd ) 
        \left\vert
            \sum _x   L(x\ida \inv \idd \inv ) \Phi (x)
        \right\vert ^2 
    \right)^{1/2} ,
\end{equation}
where the variables $\idd ,x$ run through the same ranges as above.

By Lemma \ref{lem:divisor-moment}, the first factor is bounded by a $\log N$ power, so we obtain
\begin{equation}
    \sum _\idd 
    \Nrm (\idd )
    \left\vert
        \sum _x L(x\ida \inv \idd \inv ) \Phi (x)  
    \right\vert ^2
    \ge \mypseudopolydecay N^{2n} .
\end{equation}
By the dyadic decomposition of the range of $\Nrm (\idd )$, we get some integer $D$ with
$N^{n/3}\le D \ll _K N^{2n/3} $
such that
\begin{equation}
    \sum _{
            \substack{\idd \\ \Nrm (\idd )\in [D,2D] }
        }
    D 
    \left\vert
        \sum _{
            \substack{
                x\in [\pm N]^n_{\ida }\cap \ida\idd \\ 
                \Nrm (x\ida\inv )\ge N^{n/3}\Nrm (\idd )
            }    
        }
        L(x\ida \inv \idd \inv ) \Phi (x)
    \right\vert 
    \ge 
    \mypseudopolydecay N^{2n} .
\end{equation}
We may also assume that the ideals $\idd $ appearing in the above sum belong to a common ideal class (say represented by $\idd _0$) by the pigeonhole principle.
Let $\calD \subset (K\otimes \R )\baci $ be an (always norm-length compatible) fundamental domain for the action of $\OK \baci $.
As $\idd _0\inv \cap \calD $ parametrizes the ideals $\idd $ in the ideal class of $\idd_0\inv $,
the above inequality can be rewritten as 
\begin{multline}\label{eq:type-II-before-Fourier-expansion}
    \sum _{
            \substack{
                d\in \idd_0\inv \cap \calD \\
                \Nrm (d\idd _0 )\in [D,2D]
            }
        }
    D
    \left\vert 
        \sum _{
            \substack{
                \alpha \in [\pm O_K(1)\frac N{D^{1/n}}]^n_{\ida\idd_0 }\cap 
                d\inv [\pm N]^n_{\ida } \\ 
                \Nrm (\alpha \ida\inv )\ge N^{n/3}\Nrm (\idd _0)
            }    
        }
        L(\alpha \ida \inv \idd _0 \inv ) \Phi (\alpha d)
    \right\vert ^2 
    \\ 
    \ge 
    \mypseudopolydecay N^{2n} 
    ,
\end{multline}
where we have let $d\in \idd_0\inv \cap \calD $ correspond to $\idd = d\idd _0 $
and $\alpha \in \ida\idd_0 $ to $x=\alpha d\in \ida\idd $.

For each $d$, the conditions for $\alpha $ defines a subset $\Omega _d $ in 
$ [\pm O_K(1)\frac N{D^{1/n}}]^n_{\ida\idd_0\otimes \R } $
whose boundary has a 
$\mLip (n,O_n(1),O_K\left( \frac N{D^{1/n}}\right) )$ 
boundary by \S \ref{appendix:lattice-points}.
We take a Fourier expansion of its indicator function as in \eqref{eq:indicator-function-Fourier-expansion} with $N\leadsto \frac{N}{D^{1/n}}$, $M:=\mypseudopolygrowth   $,
$Y:=M^2$ with a larger implied constant than above:
\begin{equation}\label{eq:Fourier-of-indicator-Type-II}
    1_{ \Omega _d
    }
    =
    \left( 
    \sum _{\theta \in \Theta }
    c_\theta ^{(d)}
    \cdot 
    e^{2\pi i \theta (-)} 
    \right) 
    +B^{(d)}
    +H ^{(d)}
    ,
\end{equation}
where: $\Theta $ 
is the following set of homomorphisms:
\[ \Theta := \Hom (\ida\idd_0 , \frac 1Y \Z / \Z )  \subset \Hom (\ida\idd_0 ,\RZ ); \] 
%
%
the coefficients satisfy $|c_\theta ^{(d)}|\le 1$; 
and the functions $B^{(d)}$, $H^{(d)}$ satisfy  
\begin{itemize}
    \item 
    $\lnorm{B^{(d)}}_\infty \le 1$, 
    \item 
    $\# \Supp (B^{(d)})< \mypseudopolydecay \frac{N^n}{D}$, 
    \item 
    $\lnorm{H^{(d)}}_\infty \le \mypseudopolydecay $.
\end{itemize}
Note that the index set $\Theta $ is independent of $d$ and its cardinality is $|\Theta | = Y^n<\mypseudopolygrowth   $.

Since we have given the new $\mypseudopolydecay $ a larger implied constant in its definition,
the contributions of $B^{(d)}, H^{(d)}$ to the left-hand side of 
\eqref{eq:type-II-before-Fourier-expansion} 
is negligible compared to the right-hand side.
By this and the triangle inequality of the $\ell ^2$-norm, we obtain 
\begin{align}\label{eq:l2-triangle-inequality}
    &\mypseudopolydecay N^{n} 
    \overset{\eqref{eq:indicator-function-Fourier-expansion}} 
    {\le}
    \\
    &\left(     
    \sum _{
            \substack{
                d\in \idd_0\inv \cap \calD \\
                \Nrm (d\idd _0 )\in [D,2D]
            }
        }
    D
    \left\vert 
        \sum _{\theta \in \Theta }
    \left(    
        \sum _{
            \substack{
                \alpha \in [\pm O_K(1)\frac N{D^{1/n}}]^n_{\ida\idd_0 } 
            }    
        }
        c^{(d)}_\theta \cdot 
        e^{2\pi i \theta (\alpha )}
        L(\alpha \ida \inv \idd _0 \inv ) \Phi (\alpha d)
    \right)
    \right\vert ^2
    \right)^{1/2} 
    \\  
    &\overset{\ell^2\text{-norm triangle}}{\le}
    \sum _{\theta \in \Theta }
    \left(     
    \sum _{
            \substack{
                d\in \idd_0\inv \cap \calD \\
                \Nrm (d\idd _0 )\in [D,2D]
            }
        }
    D
    \left\vert 
        \sum _{
            \substack{
                \alpha \in [\pm O_K(1)\frac N{D^{1/n}}]^n_{\ida\idd_0 } 
            }    
        }
        e^{2\pi i \theta (\alpha )}
        L(\alpha \ida \inv \idd _0 \inv ) \Phi (\alpha d)
    \right\vert ^2
    \right)^{1/2}
    .
\end{align}
(To see this, for each $\theta \in \Theta $ consider the vector $\bm v _\theta $ whose entries are indexed by $d$'s
and given by the inner sum $D^{1/2}c^{(d)}_\theta\sum _{\alpha }  e^{2\pi i \theta (\alpha )} L(\alpha \ida \inv \idd _0 \inv ) \Phi (\alpha d)$. 
The above inequality is nothing but $\lnorm{\sum _\theta \bm v_\theta }_2\le \sum _\theta \lnorm{ v_\theta}_2 $
together with $|c_\theta ^{(d)}|\le 1$.)

By the pigeonhole principle, there is some $\theta \in \Theta $ such that
the same type of inequality holds for this single $\theta $.
Now we expand the $|-|^2$ part, introducing the quantity 
$B(\alpha ,\alpha '):= e^{2\pi i \theta (\alpha -\alpha ')} L(\alpha \ida \inv \idd _0 \inv ) \overline{L(\alpha ' \ida \inv \idd _0 \inv )} $:
\begin{equation}\label{eq:expand-the-square}
    \sum _{\alpha ,\alpha '\in [\pm O_K(1)\frac N{D^{1/n}}]^n_{\ida\idd_0 } } 
    B(\alpha ,\alpha ')
     \Phi (\alpha d) \overline{\Phi (\alpha ' d)} .
\end{equation}
By the pointwise bound $|B(\alpha ,\alpha ')|\le (\log N)^2$,
Proposition \ref{prop:Cauchy-Schwarz} (a version of Cauchy-Schwarz) implies
\begin{equation}
    \left\vert
    \sum _{
            \substack{
                d\in 
                \calD 
                \cap 
                [\pm O_K(1)D^{1/n}]^n_{\idd_0\inv }
                \\
                \Nrm (d\idd _0 )\in [D,2D]
            }
        }
        D
        \cdot 
        \eqref{eq:expand-the-square}
        \right\vert
        \ll _K 
        (\log N)^2
        N^n
        \left\vert  
            \sum _{d,d'}\sum _{\alpha ,\alpha '}
            \Phi (\alpha d)\Phi (\alpha ' d') \overline{\Phi (\alpha 'd)\Phi (\alpha d')}
        \right\vert ^{1/2} ,
\end{equation}
where $d,d'$ and $\alpha ,\alpha '$ are subject to the same conditions as above.

Combined with the fact that \eqref{eq:l2-triangle-inequality} holds for our chosen $\theta $, we obtain
\begin{equation}
    \mypseudopolydecay N^{2n} 
    \le 
        \left\vert  
            \sum _{d,d'}\sum _{\alpha ,\alpha '}
            \Phi (\alpha d)\Phi (\alpha ' d') \overline{\Phi (\alpha 'd)\Phi (\alpha d')}
        \right\vert 
        =
        \sum _{d,d'}
        \left\vert 
            \sum _\alpha \Phi (\alpha d)\overline{\Phi (\alpha d')}
        \right\vert ^2
        .
\end{equation}
We know $\lnorm{\Phi }_\infty < \mypseudopolygrowth   $, the number of $\alpha $'s is $\ll _K N^n / D$ and there are $\ll _K D^2$ many $(d,d')$.
Thus with a larger implied constant in $\mypseudopolydecay $, it follows that there are $\gg _K \mypseudopolydecay D^2$ pairs $(d,d')$ such that 
\begin{equation}
    \left\vert 
            \sum _\alpha \Phi (\alpha d)\overline{\Phi (\alpha d')}
        \right\vert ^2
        \ge 
        \mypseudopolydecay \frac{N^{2n}}{D^2} .
\end{equation}
In terms of average, this says 
\begin{equation}\label{eq:Type-II-before-equidistribution}
    \left\vert 
            \Expec _{\alpha } F (g(\alpha d)\Gamma )\overline{F (g(\alpha d')\Gamma )}
        \right\vert 
        \ge 
        \mypseudopolydecay 
         \qquad \text{ for each } (d,d').
\end{equation}
\subsection{Applying equidistribution theory}

To use Leng's equidistribution theorem \ref{thm:Leibman},
we interpret inequality \eqref{eq:Type-II-before-equidistribution} 
in terms of the next polynomial map and \Lip function 
\begin{equation}
    \ida\idd_0 
    \xto{(g(-)d , g(-)d')} G\times G 
    \surj (G\times G)/(\Gamma \times \Gamma )
    \xto{F^\square :=\pr _1^*F(-)\cdot \pr _2^*\overline{F(-)}}
    \C .
\end{equation}
%
%
%
Here, $F^\square $ is the function $F^\square (z_1,z_2)=F(z_1)\cdot \overline{F(z_2)}$.
Note that $F^\square $ has vertical frequency 
$\pr _1^*\xi -\pr_2^*\xi $, i.e., for all 
$(g _1,g _2)\in G_{(\nilcounter )}\times G_{(\nilcounter )}$
and 
$(z_1,z_2)\in G\times G$ , we have 
\begin{equation}
    F^\square (g_1 z_1,g_2 z_2) = e^{2\pi i \xi (g _1) }\cdot e^{-2\pi i \xi (g _2) } 
    \cdot 
    F^\square ( z_1,z_2 ) .
\end{equation}
Leng's theorem gives us 
\begin{itemize}
    \item an integer $1\le r_{(d,d')}\le 2\dim (G/[G,G])$;
    \item a horizontal character $\eta _{(d,d')}^\etacounter \colon G\times G\to \R $ 
    of size $\le \mypseudopolygrowth   $
    for each $\etacounter =1,\dots ,r_{(d,d')}$. 
    Each of them is necessarily of the form 
    \[ 
        \eta _{(d,d')}^\etacounter =\eta _{(d,d')}^{\etacounter ,1}+\eta _{(d,d')}^{\etacounter ,2}
        \quad \colon \quad 
        (z_1,z_2)
        \mapsto 
        \eta _{(d,d')}^{\etacounter ,1}(z_1)
        +\eta _{(d,d')}^{\etacounter ,2}(z_2);  
    \]
    \item at least $ \mypseudopolydecay D ^2$ many pairs $(d,d')\in \calD $ with $\Nrm (d\idd_0),\Nrm (d'\idd_0)\in [D,2D]$; 
\end{itemize}
such that for each of these $(d,d')$:
\begin{enumerate}[(i)]
    \item 
    if we define $G_{(d,d')}:= \bigcap _{\etacounter =1}^{r_{(d,d')}} \ker (\eta ^{\etacounter }_{(d,d')}) \subset G\times G$
    we have the inclusion $(G_{(d,d')})_{(\nilcounter )}\subset \ker (\pr_1^*\xi - \pr_2^*\xi )$, and  
    \item we have the inequality:
    \begin{equation}
        \lnorm{ 
            \eta _{(d,d')}^\etacounter 
            \circ 
            \Bigl( g(-)d , g(-)d'\Bigr)
        }_{\Cinfty [\pm \frac N{D^{1/n}}]^n_{\ida\idd_0 }} 
        \le 
        \mypseudopolygrowth   
        .
    \end{equation}
\end{enumerate}
By the pigeonhole principle for the different values of $r_{(d,d')}$ ($<2(\log N  ) ^{ \mydecayexponent }$), by slightly changing the implied constants in $\mypseudopolydecay $,
we may assume that 
\[ r_{(d,d')}=r \]  
is independent of $(d,d')$ 
for those $> \mypseudopolydecay D^2$ many pairs $(d,d')$.
Also, as there are only $<\mypseudopolygrowth   $ many horizontal characters within the said size bound, we may assume that for each $\etacounter =1,\dots ,r $  
\[ \eta ^\etacounter _{(d,d')}=\eta ^\etacounter = \eta ^{\etacounter ,1}+\eta ^{\etacounter ,2} \] 
is independent of $(d,d')$.
In this notation, if we define $(G\times G)_1:= \bigcap _{\etacounter =1}^r \ker (\eta ^\etacounter ) \subset G\times G $, condition (i) above looks like: 
\begin{equation}\label{eq:GxG1-contained-in-kernel}
    \Bigl( (G\times G)_1 \Bigr) _{(\nilcounter )}
    \subset 
    \ker (\pr_1^*\xi - \pr_2^*\xi )
 . 
\end{equation}
%
%
Let us change the convention of notation and write $(d^1,d^2)$ for $(d,d')$.
We have 
\begin{equation}\label{eq:common-eta}
    \lnorm{\eta ^{\etacounter ,1}\circ g((-)d^1) + \eta^{\etacounter ,2}\circ (g((-)d^2)) }
    _{\Cinfty [\pm \frac N{D^{1/n}}]^n_{\ida\idd_0 }}
    \le \mypseudopolygrowth   
\end{equation}
for $> \mypseudopolydecay D^2$ many pairs $(d^1,d^2)$ and for all $\etacounter =1,\dots ,r $.

\subsection{Applying Vinogradov's lemma}

We fix $\etacounter \in \{ 1,\dots ,r \}$ for the moment. 
We will sometimes not reflect $\etacounter $ in the notation when we introduce quantities that depend on $\etacounter $. 

Let $H^1,H^2$ be the polynomial maps 
\begin{equation}
    H^{\oneortwo }\colon 
    \ida\idd_0 \times \idd_0\inv 
    \xto{\mathrm{mult}} 
    \ida 
    \xto{g} G
    \xto{\eta ^{\etacounter ,\oneortwo }}
    \R  
    \qquad (\oneortwo =1,2)  .
\end{equation}
If we write $\alpha =(\alpha _1,\dots ,\alpha _n)$ and $d=(d_1,\dots ,d_n)$
in some norm-length compatible bases of $\ida\idd_0 $ and $\idd_0\inv $,
these maps are written by their coefficients $h^\oneortwo _{\undli ,\undlj }\in \R $:
\begin{equation}\label{eq:H-polynomial-expansion}
    H^m(\alpha ,d)
    =
    \sum _{\undli ,\undlj } 
    h^m_{\undli ,\undlj }
    \alpha ^{\undli } 
    d^{\undlj },
    \qquad h^m_{\undli ,\undlj } =0 \text{ unless } |\undli |=|\undlj |
    . 
\end{equation}
For each $(d^1,d^2)$, 
the map $\alpha \mapsto \eta ^{\etacounter ,1}(g(\alpha d^1)) + \eta ^{\etacounter ,2}(g(\alpha d^2))$
is a polynomial map $\ida\idd_0 \to \R $ whose $\undli \ordinalth $ coefficient is
\begin{equation}
    \sum _{
        \substack{\undlj 
        \\ 
        |\undlj |=|\undli |}
        }
    \left(
        h^1_{\undli ,\undlj } (d^1)^{\undlj }
        +
        h^2_{\undli ,\undlj } (d^2)^{\undlj }
    \right) .
\end{equation}
By the definition of the $\Cinfty $-norm, the inequality \eqref{eq:common-eta} implies that for each $\undli \neq 0$,
\begin{equation}\label{eq:coefficients-bound-type-II}
    \lnorm{
        \sum _{
            \substack{\undlj 
            \\ 
            |\undlj |=|\undli |}
            }
        \left(
            h^1_{\undli ,\undlj } (d^1)^{\undlj }
            +
            h^2_{\undli ,\undlj } (d_0^2)^{\undlj }
        \right)
    }_{\RZ }
    \le 
    \mypseudopolygrowth   
    \frac{D^{|\undli |/n}}{N^{|\undli |}}
    .
\end{equation}
For a fixed $\undli $, the map $\idd_0\inv \oplus \idd_0 \inv \to \R $ given by 
\begin{equation}
    (d^1,d^2)\mapsto 
    \sum _{
        \substack{\undlj 
        \\ 
        |\undlj |=|\undli |}
        }
    h^1_{\undli ,\undlj } (d^1)^{\undlj }
    +
    h^2_{\undli ,\undlj } (d_0^2)^{\undlj }
\end{equation}
is homogeneous of degree $|\undli |$.
Since we know that the inequality \eqref{eq:coefficients-bound-type-II} holds 
for the
$> \mypseudopolydecay D^2 $ many pairs 
$(d^1,d^2)\in [\pm O_K(D^{1/n})]^n_{\idd_0\inv\oplus \idd_0\inv }$,
Proposition \ref{prop:Vinogradov} (polynomial Vinogradov) 
implies that 
there is a positive integer $1\le q_{\undli }\le \mypseudopolygrowth   $ such that
if we multiply the above polynomial map by $q_{\undli }$, the resulting map has $\Cinfty [\pm D^{1/n}]^n_{\idd_0\inv \oplus \idd_0\inv}$-norm $<\mypseudopolygrowth   \frac{D^{|\undli |/n}}{N^{|\undli |}}$.
In concrete terms, this means
\begin{equation}
    \lnorm{q_{\undli} h^\oneortwo _{\undli ,\undlj } }_{\RZ}
    <
    \mypseudopolygrowth   \frac{1}{N^{|\undli |}} \qquad \text{ for all } \undli ,\undlj \neq 0,\ \oneortwo =1,2.
\end{equation}
By replacing them by the least common multiple, we may assume that
\[ q_{\undli }=q \] 
are independent of $\undli $ and $\etacounter =1,\dots ,r$.

Recall from \eqref{eq:H-polynomial-expansion} that $h^\oneortwo _{\undli ,\undlj }$ are the coefficients of the polynomial map $H^\oneortwo \colon \ida\idd_0\times \idd_0\inv \to \R $.
Since $1\in \idd_0\inv $ is represented by a vector of size $O_K(1)$ in the chosen basis,
by substituting $d=1$ in the polynomial expansion of $H^\oneortwo $ and using the above inequality, 
we conclude that the $\undli \ordinalth $ coefficients of the map $qH^\oneortwo (-,1)=q\eta ^{\etacounter ,\oneortwo }\circ (g|_{\ida\idd_0})\colon \ida\idd_0\to \R $
have $\RZ $-norm $< \mypseudopolygrowth   \frac{1}{N^{|\undli |}}$ for all $\undli \neq 0$.
Equivalently, we have
\begin{equation}
    \lnorm{q\eta ^{\etacounter ,\oneortwo }\circ (g|_{\ida\idd_0})}_{\Cinfty [\pm N]^n_{\ida\idd_0 }}
    \le \mypseudopolygrowth   
    .   
\end{equation}
Since $\ida\idd_0\inj \ida $ has index $\le O_K(1)$, we conclude 
\begin{equation}\label{eq:almost-allows-us}
    \lnorm{q\eta ^{\etacounter ,\oneortwo }\circ g}_{\Cinfty [\pm N]^n_{\ida}}
    \le \mypseudopolygrowth   
    \qquad \te{ for all } \etacounter =1,\dots ,r,
    \quad \oneortwo =1,2.   
\end{equation}
\subsection{End of proof}

Now we let $\etacounter $ vary again.
Inequality \eqref{eq:almost-allows-us} almost allows us to apply Lemma \ref{lem:our-argument} to the $2r$ horizontal characters $\eta ^{\etacounter ,\oneortwo }$,
except that we still have to check the inclusion
\begin{equation}
    (G_1)_{(\nilcounter )} :=
    \left( 
        \bigcap _{\etacounter ,\oneortwo } \ker (\eta ^{\etacounter ,\oneortwo } )
    \right) 
    _{(\nilcounter )} 
    \subset 
    \ker (\xi ) .
\end{equation}
By the definitions we have $G_1\times G_1 \subset (G\times G)_1$ (see \eqref{eq:GxG1-contained-in-kernel} for the definition of $(G\times G)_1$).
So, a similar inclusion hold for their $\nilcounter \ordinalth $ commutator subgroups.
Consider the antidiagonal embedding (which is not a group homomorphism 
unless $G$ is commutative,
but it does not matter): 
\[ \Delta ' \colon\quad G\to G\times G ,\quad  x\mapsto (x,x\inv ).\]
By the above inclusions, this carries $(G_1)_{(\nilcounter )}$ into $((G\times G)_1)_{(\nilcounter )}$, which is contained in $ \ker (\pr_1^*\xi - \pr_2^*\xi )$ by \eqref{eq:GxG1-contained-in-kernel}.
Since the composite 
\[ 
    G_{(\nilcounter )}
    \xto{\Delta '} 
    G_{(\nilcounter )}\times G_{(\nilcounter )} 
    \xto{\pr_1^*\xi - \pr_2^*\xi }
    \R
\] 
is $2\cdot \xi $, it follows that $(G_1)_{(\nilcounter )}\subset \ker (\xi )$ as desired.

Now that we have checked the hypotheses of Lemma \ref{lem:our-argument} for the horizontal characters $\eta ^{\etacounter ,\oneortwo }$ ($\etacounter \in \{ 1,\dots ,r\} $, $\oneortwo =1,2$), 
the lemma creates data satisfying the conclusion of 
Proposition \ref{prop:step-goes-down}.
This establishes the Type II case of Proposition \ref{prop:step-goes-down}. 

\ \ \ 

    Combined with \S \ref{sec:Type-I}, this completes the proof of Proposition \ref{prop:step-goes-down}. 

\ \ \ 
    
As we saw in \S \ref{sec:induction-on-step}, Proposition \ref{prop:step-goes-down} implies Theorem \ref{thm:asymptotic-orthogonality-with-nilsequences}.
    Since we saw in \S \ref{sec:relation-to-asymptotic-orthogonality} that 
    this theorem was equivalent to Theorem \ref{thm:vonMang-Siegel-close!},
    this latter theorem has also been proven.

\section{Simultaneous prime values of affine-linear forms}\label{sec:prime-values}

Our next main result is Theorem \ref{thm:prime-values}, which generalizes the Green--Tao--Ziegler theorem \cite[Main Theorem (p.~1763)]{LinearEquations} \cite[Theorem 1.1]{MobiusNilsequences} \cite[Theorem 1.3]{GowersInverse} on simultaneous prime values of finite-complexity affine-linear forms $\Z^d \to \Z $ to the case where the target is the ring of integers $\OK $ of a general number field $K$.
In fact, for aesthetic and practical reasons it is better to allow the target to be also a non-zero fractional ideal.
This extended flexibility will be useful in \S \ref{sec:Hasse} via Appendix \ref{sec:justification-prime-values-localized}.
See also \cite[proof of Prop.~3.1]{ABHS}, 
\cite{Zywina2025a} for how this general form can be useful for the construction of abelian varieties with certain arithmetic properties.

Recall that for an affine linear map $\psi \colon \Zd \to \Zn $
written as $\psi (x_1,\dots ,x_d)=c_0+\sum _{i=1}^d c_{i}x_i$
(where $c_i\in \Zn $)
and a positive number $N>1$, we let $\lnorm{c_i}_{\infty }$ be the $\ell ^{\infty }$-norm and 
 define 
\begin{align}
    \lnorm{\psi }_N := 
    \max \left\{  \frac{\lnorm{c_0}_{\infty }}N  , \max _{1\le i\le d} \lnorm{c_i}_{\infty } \right\} .
\end{align}
This quantity will be considered when the target of the map is a fractional ideal $\ida \cong \Zn $ equipped with a basis.

Recall also when $\ida $ is a non-zero fractional ideal, we define $\Lambdaa (x):= x\ida\inv $ for $x\in \ida $.

\begin{theorem}[simultaneous prime values]
    \label{thm:prime-values}
    Let $\ida\cong \Zn $ be a non-zero fractional ideal of a number field $K $ equipped with a norm-length compatible basis.
    Let $t,d\ge 2$ be integer parameters
    and $A>1$ an arbitrary positive real number.
    Let $N>1$ be an integer large enough depending on $t,K,A$.
    Let 
    \[ \psi _1,\dots ,\psi _t\colon \Z ^d \to \ida \] 
    be affine-linear maps with $\lnorm{\psi _i}_N \le (\log N)^A$
    satisfying 
    \begin{quote}
        $\dot\psi _i |_{\ker (\dot\psi _j)}$ has finite cokernel for all $i\neq j$
        (which implies that $\psi _1,\dots ,\psi _t$ have finite {\em complexity} $\le s$ for some  $s\le t$ (Definition \ref{def:complexity})).
    \end{quote}
    Let $\Omega \subset [\pm N]^d $ be a convex set. Then  
        \begin{multline}\label{eq:prime-values-goal}
            \sum _{x\in \Omega \cap \Z^d}
            \left( \prod _{i=1}^t
            \Lambda ^{\ida}_K (\psi _i (x)) \right)
            =
            \frac{\Vol (\Omega )}{\residue\limits _{s=1}(\zeta_K (s))^t}
            \cdot
            \prod _{p } \beta _{p }(\psi _1,\dots ,\psi _t)
            \\ 
            +
            O_{t,d,K,A} (N^d (\log N)^{-A} )
            ,
        \end{multline}
        where 
        the quantities $\beta _p(\psi _1,\dots ,\psi _t)$ are defined by 
        \begin{align}
            \beta _{p }(\psi _1,\dots ,\psi _t):= 
            \left(\frac{p^n }{\totient (p) } \right)^t
            \Expec _{x\in (\Z /p\Z )^d}  1[\psi _i (x)\neq 0 \te{ in }\ida /\idp \ida \te{ for all $i$ and $\idp \divides p $}] 
            .
        \end{align}
    \end{theorem}
    
Before proving this, let us make a few remarks on the statement.
    
\begin{remark}
    The reader does not have to know what a {\em norm-length combatible basis} of $\ida $ (\S \ref{appendix:norm-length}) is 
    if they are willing to allow the implied constants to depend on the fractional ideal $\ida $ (such as when one is interested only in the case $\ida =\OK $). 
    They may take whatever basis they like.
    
    The independence of the implied constants is used in \S\S \ref{sec:Hasse}, \ref{sec:justification-prime-values-localized}, where we apply Theorem \ref{thm:prime-values} to unboundedly many fractional ideals $\ida $.
\end{remark}

\begin{remark}\label{rem:product-beta}
    The product $\prod _p \beta _p (\psi _1,\dots ,\psi _t)$ is either zero (precisely when one of the $\beta _p$'s is zero) or an absolutely convergent product (in particular a positive value) by \eqref{eq:1-g(p)_for_large_p}.
    We have $\beta _p(\psi _1,\dots ,\psi _t)\neq 0$ if and only if 
    there is an $x\in \Z^d$ such that 
    $\psi _i(x)\neq 0$ in $\ida /\idp\ida $ for all $\idp\divides p$ and all $i$,
    or equivalently $p\doesnotdivide \Nrm (\psi _i(x)\ida\inv )$ for all $i$.

    When $\Z^d$ carries an $\OK $-module structure with respect to which $\psi _i$ are all affine-$\OK $-linear, then as was mentioned in Proposition \ref{prop:sieve-theory} we have 
    a decomposition 
    \begin{align}
        \beta _p (\psi _1,\dots ,\psi _t)
        &=
        \prod _{\idp \divides p} \beta _{\idp } (\psi _1,\dots ,\psi _t)
        .
    \end{align}
    We have $\beta _\idp (\psi _1,\dots ,\psi _t)\neq 0$ if and only if 
    there is an $x\in \Z^d $ such that $\psi _i(x)\neq 0 $ in $\ida / \idp\ida $ for all $i$.
\end{remark}

\begin{remark}
    We would have liked to state and prove Theorem \ref{thm:prime-values} by only assuming that
    $\ker(\dot\psi _i)$'s do not contain each other (making no difference when $K=\Q$),
    so that the theorem would determine the correct asymptotic number of 
    constellations in prime numbers studied in our former work \cite{KMMSY}.
    This limitation comes from our current form of the von Neumann Theorem \ref{thm:vonNeumann}.
    It is unclear whether this limitation can be overcome just by some clever linear algebra trick or it is more deeply rooted.
\end{remark}

\begin{proof}[Proof of Theorem \ref{thm:prime-values}]
    
    Write $L:=(\log N)^A$.
    Note that $\psi _i (\Omega )$'s are contained in the cube $[\pm (d+1)L N]^n_{\ida }$.
    Let $B>1$ be a large number to be specified later.
    By Corollary \ref{cor:vonMang-Cramer-close!} we know 
    \[ \lnorm{\frac{\Lambdaa }{\log N} - \frac{\LambdaaCramer }{\log N} }
    _{U^{s+1}[\pm (d+1)L N]^n_{\ida }} 
    \ll _B (\log N)^{-B} . \]   
    Multiple applications of the generalized von Neumann theorem \ref{thm:vonNeumann} 
    where the functions $f_i$ are set to be 
    $\frac{\Lambdaa }{\log N}$, 
    $\left( \frac{\Lambdaa }{\log N}-\frac{\LambdaaCramer }{\log N}\right) $
    and 
    $\frac{\LambdaaCramer }{\log N}$
    imply that the left-hand side of \eqref{eq:prime-values-goal} is equal to
    \begin{equation}
        \sum _{x\in \Omega \cap \Z^d} \prod _{i=1}^t \LambdaaCramer (\psi _i (x))
        +
        O_{t,d,K,B} \left(N^d (\log N)^{t+Ad}(\log N)^{-B/(2d+2)} \right)
    \end{equation}
    The exponent of $\log N$ is arbitrarily negative when $B>1$ is arbitrarily large.
    Set $B:= 2(d+1)^2 A + 2(d+1)t $ to get $(\log N)^{-A}$, say. 
    By Proposition \ref{prop:sieve-theory}, the main term here is close to the main term of the right-hand side of \eqref{eq:prime-values-goal} with a much smaller error than aimed for.
    This completes the proof of Theorem \ref{thm:prime-values}.
\end{proof}

\section{Hasse principle of rational points for certain fibrations}\label{sec:Hasse}

Here we prove a Hasse principle result, Theorem \ref{thm:HSW-for-K}, which was formerly known over $\Q$ by Harpaz--Skorobogatov--Wittenberg \cite{Harpaz-Skorobogatov-Wittenberg}.
For this application, we need Theorem \ref{thm:prime-values-S}, a variant of Theorem \ref{thm:prime-values} for the localized ring $\OKS $.

\subsection{The localized variant of main theorem}

For a finite set $S$ of non-zero prime ideals of $\OK $,
define $\OKS $ to be the ring obtained by inverting all $\idp \in S$;
a pedantic formal definition of it is 
$\OKS = \Gamma (\SpecOK \setm S ,\calO )$, where $\calO $ is the structure sheaf.
A more concrete description is, 
since $\Cl (K)$ is torsion (even finite), some power of the ideal $\prod _{\idp \in S}\idp $
is a principal ideal, say $f\OK $.
Then one has $\OKS = \OK [f\inv ]$.

For a fractional ideal $\ida \subset 
K $ of $\OK $, write  
$\idas := \ida \otimes _{\OK }\OKS $. 

Let $\Ideals _{\OKS }$ be the multiplicative monoid of non-zero ideals of $\OKS $.
It is the free commutative monoid on the set of non-zero prime ideals of $\OKS $,
which is identified with the set of non-zero prime ideals of $\OK $ away from $S$.
Let $\Nrm =\Nrm _{S\inv }\colon \Ideals _{\OKS }\to \N $ be the norm map 
$\idb \mapsto |\OKS /\idb |$. We omit the subscript except when there is risk of confusion.

Define the \vMang function of the ring $\OKS $
\[ \Lambda _{\OKS }\colon \Ideals _{\OKS }\cup \{ (0) \}\to \R_{\ge 0}\]
by 
$\idb \mapsto \log (\Nrm (\idp )) $ 
if $\idb $ is a non-trivial power of a prime ideal $\idp $ of $\OKS $,
and $\idb \mapsto 0$ otherwise.
For a non-zero fractional ideal $\ida $ over $\OK $, define the function
\begin{align}
    \Lambda ^{\ida }_{\OKS }\colon \ida \to \R _{\ge 0}
\end{align}
by $x\mapsto \Lambda _{\OKS }(x\idas\inv )$.

A routine but tedious argument lets one deduce the following from our main results.

\begin{theorem}[simultaneous prime values, away from $S$]
    \label{thm:prime-values-S}
    Let $S$ be a finite set of non-zero prime ideals of $\OK $ and 
   $\ida \subset K $ be a non-zero fractional ideal of $\OK $.
    Let $t,d\ge 1$ be integers and 
   $\Psi = \{ \psi _1,\dots ,\psi _t\} $ be a set of affine-linear maps 
   \begin{align}
    \psi _i\colon \Z^d \to \idaS  
   \end{align} 
   having torsion cokernels such that for every pair $i\neq j$
   the restriction $\psi _i|_{\ker (\psi _j)}$ also has torsion cokernel.

   Let $A>1$ be an arbitrary positive number and $N\gg _{\Psi ,K,S,A} 1$. Let $\Omega \subset [\pm N]^d\subset \R^d $ be a set with boundary of \Lip class 
   $\mLip (d,Z,L_0N)$ for some $Z,L_0\ge 1$.
   ({\em See \S \ref{appendix:lattice-points} for this notion;
   convex sets have boundary of \Lip class $\mLip (d,1,O_d(N) )$.})
   
   Then we have 
   \begin{align}
        \sum _{\undl x \in \Omega \cap \Z ^d}
        \prod _{i=1}^t \Lambda _{\OKS }^{\ida}(\psi _i (\undl x ) )
        =
        C_{S,\Psi }\frac{\vol (\Omega )}{\residuezeta ^t}
        +O_{K,S,\Psi ,Z,L_0,A}(N^d (\log N)^{-A} )
        ,
   \end{align}
   where the coefficient $C_{S,\Psi }$ is as follows:
   \begin{align}
        C_{S,\Psi }
        := 
        \prod _{p } 
        \left( \frac{p^n}{\totient (p)} \right)^t 
        \left( 
            1- \frac{1}{p^d}
            \left|
                \bigcup _{i=1}^{t}\bigcup_
                {\substack{\idp \divides p \\ \te{\upshape not in }S }}  
                \psi _{i,\idp }\inv (0)
            \right|  
        \right)
        ,
   \end{align}
   where for each prime ideal $\idp $ not in $S$ with residue characteristic $p$, the map $\psi _{i,\idp }\colon \Z^d /p\Z ^d \to \ida /\idp \ida $ is the affine-linear map induced by $\psi _i$.
   The coefficient $C_{S,\Psi }$ is positive if the following condition is satisfied:
   \begin{align}\label{eq:Schinzel-condition}
    \parbox{10cm}{
    for all $p$, there is an $\undl x\in \Z^d $ such that $\Nrm (\psi _i(\undl x)\idas\inv )$ is prime to $p$ simultaneously for all $1\le i\le t$.
    }
\end{align}
When the source $\Z^d $ admits an $\OK $-module structure with respect to which $\psi _i$ are affine-$\OK $-linear, then the last condition is equivalent to:
\begin{align}
    \parbox{10cm}{
    for all $\idp \notin S$, there is an $\undl x \in \Z^d$
    such that $\psi _i(\undl x )$ are non-zero in $\ida /\idp\ida $ simultaneously for all $1\le i\le t$.}
\end{align}
\end{theorem}

The deduction of Theorem \ref{thm:prime-values-S} is not the main focus of this section, and thus is deferred until \S \ref{sec:justification-prime-values-localized}.

\subsection{Application to rational points}

Below, for a place $v$ of a number field $K$,
write $K_v$ for the completion.
Let $| - |_v\colon K_v\to \R_{\ge 0}$ be the associated absolute value, which we choose (say) as 
$|x|_{\idp }:= \Nrm (\idp )^{-v_\idp (x)}$ if $v=\idp $ is  non-archimedean and the usual absolute value of $\R, \C$
if $v$ is archimedean.

In any metric space, let us say that two elements are {\em $\varepsilon $-close} ($\varepsilon >0$)
if their distance is $\le \varepsilon $.

Here is an analogue of \cite[Proposition 1.2]{Harpaz-Skorobogatov-Wittenberg}.
The case $\tau _v=+\infty $ and $\tau _v=\infty $ for all real and complex places $v$ will suffice for the application in this section, but the following form is somewhat more natural to formulate.

\begin{proposition}[approximation by a prime element]
    \label{prop:HSW1.2}
    Let $S$ be a finite set of non-zero prime ideals of $\OK $
    and 
    suppose we are given the data of 
    \begin{itemize}
        \item 
        $(\lambda _{\idp} ,\mu _{\idp })\in K_{\idp }^{\oplus 2}$ for each $\idp \in S$, 
        \item 
        $\tau _v \in \R \cup \{ \pm \infty \}$ and $\sigma _v\in \{ \pm 1\}$ for each real place $v$
        and 
        \item 
        $\tau _v\in \C \cup \{\infty \} $ for each complex place $v$.
    \end{itemize}
    Let $e_1,\dots ,e_t \in \OKS $ be distinct elements.

    Then for any $0<\varepsilon <1$ there exist pairs $(\lambda ,\mu )\in (\OKS \nonzero )^{\oplus 2}$
    satisfying the following conditions:
    \begin{enumerate}
        \item for all $\idp \in S$, in the $\idp $-adic $\sup$ metric on $K_\idp ^{\oplus 2}$, the pair $(\lambda ,\mu )$ is $\varepsilon $-close to $(\lambda _\idp ,\mu _\idp )$;
        \item for each archimedean place $v$,
        the ratio $\lambda /\mu $ is $\varepsilon $-close to $\tau _v$.
        Here, when $v$ is real and $\tau _v=\pm \infty $,
        this is interpreted as $\lambda /\mu $ having the same sign as $\tau _v$ and $|\lambda /\mu |> 1/\varepsilon $.
        When $v$ is complex and $\tau _v=\infty $, then as $|\lambda /\mu |>1/\varepsilon $.
        \item $\mu $ has sign $\sigma _v$ in $K_v=\R $ for each real $v$.
        \item $\lambda -e_i\mu $ are prime elements of $\OKS $ for all $i$ and no two of them are associates under the $\OKS \baci $-action.
    \end{enumerate}
\end{proposition}
\begin{proof}
    By multiplying everything by some element of $\OKS \baci $
    having large valuations at primes $\idp \in S$,
    we may assume that $(\lambda _\idp ,\mu _\idp )\in \calO _{K_\idp }^{\oplus 2}$.
    Under this additional assumption, we claim that we can find 
    $(\lambda ,\mu )\in \OK ^{\oplus 2}$ (more restrictively than $\OKS ^{\oplus 2}$)
    satisfying the asserted condition. 

    Consider the following maps for $i=1,\dots ,t$:
    \begin{align}\label{eq:linear-map-lambda-mu}
        \OK ^{\oplus 2} &\to \OKS 
        \\ 
        (\lambda ,\mu )&\mapsto \lambda -e_i\mu .
    \end{align}
    Define the ideal $\ida :=\left( \prod _{\idp \in S}\idp  \right) ^e$ with $e\gg _{\varepsilon }1$ so that every element of $\ida $ has distance $<\varepsilon $ from $0$ in $K_\idp $ for all $\idp \in S$.
    Let $(\lambda _0,\mu _0)\in \OK ^{\oplus 2}$ be a representative of 
    the class
    \begin{align}
        (\lambda _\idp ,\mu _\idp )_{\idp \in S}
        \in \bigoplus _{\idp \in S} (\calO _{K_\idp }/ \idp ^e \calO _{K_\idp })^{\oplus 2}
        \cong (\OK /\ida )^{\oplus 2}.
    \end{align}
    Choose a basis $\ida \cong \Z^n$ and consider the following composite maps $\psi _i$ for $i=1,\dots ,t$:
    \begin{align}\label{eq:def-of-psi-application}
        \begin{array}[]{rcccl}
            \psi _i\colon \Z^{2n}\cong \ida^{\oplus 2} &\to &\OK^{\oplus 2}&\to &\OKS 
            \\ 
            (x,y) &\mapsto &(\lambda _0+x,\mu _0+y) &\mapsto &(\lambda _0+x)-e_i (\mu _0 +y).
        \end{array}
    \end{align}
    The condition \eqref{eq:Schinzel-condition} is satisfied by the collection of maps \eqref{eq:linear-map-lambda-mu}
    by taking $\lambda =1,\mu =0$,
    and therefore also by $\psi _1,\dots ,\psi _t$ because the affine-$\OK $-linear map 
    $\ida ^{\oplus 2}\to \OK ^{\oplus 2}$ in \eqref{eq:def-of-psi-application} becomes a bijection modulo $\idp $ for each $\idp \notin S $.

    One easily checks that the condition that the quotient 
    $(\lambda _0+x)/(\mu _0+y)$ be $\varepsilon $-close to $\tau _v$
    for all archimedean $v$ and $\mu $ has sign $\sigma _v$ for all real $v$ defines a set 
    $\Omega \subset [\pm N]^{2n}\subset \R^{2n}\cong \ida ^{\oplus 2}\otimes \R $
    with a $\mLip (2n,O(n),O_{S,\varepsilon }(N))$ boundary and 
    volume $\gg _{S,\varepsilon }N^{2n}$.

    We can now apply Theorem \ref{thm:prime-values-S} to this situation 
    and obtain ($\gg _{S,\varepsilon ,\{ e_i\}_i } N^{2n}/(\log N)^t$ many) pairs $(x,y)\in \ida ^{\oplus 2}$,
    which in turn determine $(\lambda ,\mu ):= (\lambda _0+x,\mu _0+y)$,
    satisfying the asserted conditions.
    Note that the event that two of the $\lambda -e_i\mu $'s 
    be associates under the $\OKS\baci $-action is asymptotically negligibly rare basically because $|\OKS\baci \cap \OKN |\ll _K (\log N)^{n+|S|-1}$.
\end{proof}

The passage from Proposition \ref{prop:HSW1.2} to the following statement is identical to 
\cite[Proposition 1.2 to Proposition 2.1]{Harpaz-Skorobogatov-Wittenberg}.

\begin{proposition}[{known for $K=\Q $ e.g.\ by \cite[Proposition 2.1]{Harpaz-Skorobogatov-Wittenberg}}]
    \label{prop:HSW2.1}
    Let $S$ be a finite set of non-zero prime ideals of $\OK $.
    Let $e_1,\dots ,e_t\in \OKS $ be distinct.
    Suppose we are given the data of 
    $\tau _\idp \in K_{\idp }$ for each $\idp \in S$.

    Then for any $0<\varepsilon <1$ there exist $\tau \in K$ and distinct non-zero prime ideals 
    $\idp _1,\dots ,\idp _t\notin S $ such that 
    \begin{enumerate}
        \item $\tau $ is $\varepsilon $-close to $\tau _\idp $ in $K_\idp $ for each $\idp \in S$;
        \item $\tau $ is totally positive and $|\tau |_v >1/\varepsilon $ for all archimedean place $v$;
        \item for all $i=1,\dots ,t$, we have $v_{\idp _i} (\tau -e_i ) =1$;
        \item for all $i=1,\dots ,t$ and $\idp \notin S\cup \{ \idp _i\}$, we have $v_{\idp }(\tau -e_i)\le 0$;
        \item for any $i=1,\dots ,t$ and any cyclic extension $L/K$ unramified outside $S$,
        we have
        \begin{align}
            \oname{inv}_{\idp_i} (L/K,\tau -e_i)+\sum _{\idp \in S} \oname{inv}_{\idp } (L/K,\tau _{\idp }-e_i)=0\in \QZ .
        \end{align}
        In particular, if $\sum _{\idp \in S} \oname{inv}_{\idp } (L/K,\tau _{\idp }-e_i)=0$, then $\idp _i$ splits completely in $L/K$.
    \end{enumerate}
\end{proposition}

Thanks to Propositions \ref{prop:HSW1.2} and \ref{prop:HSW2.1}
which play the roles of \cite[Propositions 1.2 and 2.1]{Harpaz-Skorobogatov-Wittenberg},
all the results in \cite{Harpaz-Skorobogatov-Wittenberg}
now hold over an arbitrary number field $K$ rather than just $\Q $.

For example, we have the following theorem.
There, let $\mathbf{A}_K$ be the ring of \adeles of $K$.
Write $\oname{Br}_{\mathrm{vert}}(X) := \pi^*\oname{Br}(K(\mathbb{P}^1))\cap \oname{Br}(X) \subset \oname{Br}(K(X))$ for the so-called {\em vertical Brauer group} of an integral variety $X$
and set 
\begin{align}
    &X(\mathbf{A}_K)^{\oname{Br}_{\mathrm{vert}}}:= 
    \\ 
    &\{ x=(x_v)_v\in X(\mathbf{A}_K) \mid 
    \te{for all }\omega \in \oname{Br_{\mathrm{vert}}(X)},\ 
    \sum _{v\te{: all places}} \oname{inv} _v(\omega |_{x_v} ) =0 \te{ in }\QZ  \} .
\end{align}
\begin{theorem}[{known for $K=\Q $ by \cite[Theorem 3.1]{Harpaz-Skorobogatov-Wittenberg}}]
    \label{thm:HSW-for-K}
    Let $X$ be a geometrically integral variety over $K$
    equipped with a smooth surjective morphism 
    $\pi \colon X\to \mathbb{P}^1$.
    Assume:
    \begin{itemize}
        \item except for finitely many fibers $X_1,\dots ,X_t$ over $\mathbb{P}^1(K)$,
        each fiber of $\pi $ contains a geometrically integral irreducible component;
        \item the fiber $X_i$ ($1\le i\le t$) has an irreducible component $X_i^0$ such that the algebraic closure 
        of $K$ in its function field
        is abelian over $K$.
    \end{itemize}
    Then $\mathbb{P}^1(K)\cap \pi (X(\mathbf{A}_K) )$ is dense in 
    $\pi (X(\mathbf{A}_K)^{\oname{Br}_{\mathrm{vert}}}) \subset \mathbb{P}^1(\mathbf{A}_K)$.
\end{theorem}
\begin{proof}
    The proof of \cite[Theorem 3.1]{Harpaz-Skorobogatov-Wittenberg} carries over
    because it uses nothing particular about $\Q $
    except \cite[Proposition 2.1]{Harpaz-Skorobogatov-Wittenberg}, 
    now available for general number fields as Proposition \ref{prop:HSW2.1}.
    Perhaps it might seem that Step 2 of \loccit.\ is using the fact that $\Q $ has only one archimedean place, 
    but one can as usual use small deformation in $X(K_v)$'s (for $v$ archimedean) by the inverse function theorem and the density of the inclusion $K\inj K\otimes \R $ to achieve the situation where the archimedean part $(M_v)_{v}$ of the ad{\'e}lic point to approximate lies over the point $\infty \in \mathbb{P}^1(K)$.
\end{proof}

\begin{remark}
    Over $\Q$, stronger statements are shown in Harpaz--Wittenberg \cite{Harpaz-Wittenberg} using additive combinatorics results of Browning--Matthiesen \cite{Browning-Matthiesen}. 
    Therefore for this particular application, it might be more efficient to proceed by way of a number field version of Browning--Matthiesen.
\end{remark}

\appendix
\renewcommand{\theequation}{\Alph{section}.\arabic{equation}} 

\section{Results from the nilsequence theory}\label{appendix:nilsequences}

Here we list some advanced results from the theory of nilsequences that we use in this paper.
Some of them had not been formulated in the literature in the form we need (mainly, statements for functions on $\Z ^n$ as opposed to $\Z $), so we have to state and verify them in more general forms.
The reader new to this subject will benefit from sources like \cite{GreenVideos, HigherOrderFourier}, where background ideas and classic results behind these advanced results are comprehensibly explained.
See Definition \ref{def:Gowers-norm} and \S \ref{sec:nilsequences} for the recollection of basic notions.

Recall that the \Lip norm $\lnorm{-}_{\mLip }$ of a $\C$-valued function on a metric space is the sum of the $L^\infty $ norm and the \Lip constant (if finite).

\subsection{The inverse theory of Gowers norms}\label{sec:Gowers-inverse}

For positive integers $N\in \N$ recall that we write $[\pm N]^n\subset \Z ^n$ for the product of intervals $[-N,N]^n$.

Let us state a version of the Inverse Theorem for Gowers norms.
This type of statement was first conjectured (in public) by Green--Tao \cite{LinearEquations}
and subsequently proven by Green--Tao--Ziegler \cite{GowersInverse}.
There had been a major quantitative improvement by Manners \cite{Manners}
before the version below was established.

\begin{theorem}[Gowers Inverse Theorem on $\Zn $ of Leng--Sah--Sawhney] 
    \label{thm:GowersInverse}
    Let 
    $n,s,N\in \Z _{\ge 1} $.
    Let $0<\delta \ll _{n,s} 1/2$
    and $N\gg _{n,s}1$. 

    Suppose $f\colon [\pm N]^n\to \C $ is a function such that 
    \begin{align}
        \lnorm {f}_{U^{s+1}([\pm N]^n)} \ge \delta .
    \end{align}
    Then there exist:
    \begin{itemize}
        \item a filtered nilmanifold $G/\Gamma $ of degree $\le s$, dimension 
        $d \le (\log (1/\delta ))^{O_{n,s}(1)}$, 
        complexity 
        $M \le \exp ((\log (1/\delta ))^{O_{n,s}(1)})$;
        \item a $1$-bounded function $F\colon G/ \Gamma \to \C $ of Lipschitz constant 
        $\le \exp ((\log (1/\delta ))^{O_{n,s}(1)})$;
        \item a polynomial map $g\colon \Z ^n\to G$
    \end{itemize}
    such that 
    \begin{align}
        \left|  
            \Expec _{x\in [\pm N]^n} f(x)F(g(x)) 
        \right|
        \ge 
        \left(
            \exp (\log (1/\delta )^{O_{n,s}(1)})
        \right)\inv .
    \end{align}
\end{theorem}
\begin{proof}
    For $n=1$, this is exactly \cite[Theorem 1.2]{QuasipolynomialInverseGowers}.

    The assertion for $n>1$ can be deduced from the case of $n=1$ by the following application of ``Kronecker's substitution''.\footnote{This argument also appears for example in Terence Tao's blog post:
    \url{https://terrytao.wordpress.com/2015/07/24/}}

    Let $\phi \colon \Z ^n\to \Z $ be the map 
    $(x_1,\dots ,x_n)\mapsto \sum _{i=1}^n (5N)^{i-1} x_i $.
    This maps the cube $[\pm N]^n $ injectively into the interval $[\pm (5N)^n]\subset \Z $.
    Define 
    \begin{align}
        f' \colon [\pm (5N)^n]\to \C  
    \end{align}
    by setting the value equal to that of $f$
    on the image of $[\pm N]^n$, and zero otherwise.
    By the invariance of Gowers norms under Freiman isomorphisms \cite[p.\,1817]{LinearEquations} we know that the ratios of Gowers norms 
    \begin{align}\label{eq:ratios-of-Gowers-norms}
        \frac{\lnorm {f '} _{U^{s+1}([\pm (5N)^n])}}
            {\lnorm {f } _{U^{s+1}([\pm N]^n)}}
    \end{align}
    are bounded from above and below by positive numbers depending only on $n,s$.
    Then by the theorem for $n=1$ with slightly compromised exponents (to make up for the difference caused by the ratios \eqref{eq:ratios-of-Gowers-norms}) one can find data
    \[ \Z \xto{g'}G\surj G/\Gamma \xto{F}\C \]
    satisfying 
    \begin{align}
        \left|  
            \Expec _{x\in [\pm (5N)^n]} f'(x)F(g'(x)) 
            \right|
            \ge \exp (
                \log (1/\delta )^{O_{s}(1)}
            ) \inv .
        \end{align}
    Define a polynomial map $g\colon \Z ^n\to G$ by the composition (recall that every additive map $\Z ^n\to \Z $ is a polynomial map)
    \begin{align}
        g\colon \Z ^n \xto{\phi } \Z \xto {g'} G.
    \end{align}
    As we know that the expectations computed on $[\pm N]^n$ and $[\pm (5N)^n]$ are not so different:
    \begin{align}
        \left(\Expec _{x\in [\pm N]^n} f(x)F(g(x)) \middle) \quad \middle/ \quad \middle( \Expec _{x\in [\pm (5N)^n]} f'(x)F(g'(x)) \right) 
        = \frac{10N^n+1}{(2N+1)^n} \asymp _n 1,
    \end{align}
    the assertion for $n$ follows.
\end{proof}

%

\begin{proposition}[the converse to Inverse Theorem]\label{prop:converse-to-inverse-theorem}
    Let $0< \varepsilon <1/2 $, $M>10$ be positive numbers
    and $d,s\ge 2$ be positive integers.
    Let $G/\Gamma $ be a nilmanifold of dimension $\le d$, complexity $\le M$ and both degree and step $\le s$.
    Suppose furthermore we are given: a function $F\colon G/\Gamma \to \C $ satisfying 
    $\lnorm{F}_{\mLip} \le M$,
    a polynomial map $g\colon \Z^n \to G$ and a $1$-bounded function $f\colon [\pm N]^n \to \C $ such that
    \begin{align}
        \left|  
            \Expec _{x\in [\pm N]^n} f(x)F(g(x)) 
        \right|
        \ge \varepsilon .
    \end{align}  
    Then we have 
    \begin{align}
        \lnorm {f}_{U^{s+1}([\pm N]^n)} 
        \gg _{s,n}
        (\varepsilon /M )^{d^{O_{s}(1)}}. 
    \end{align}  
\end{proposition}
\begin{proof}
    This is \cite[Lemma B.5 (p.~83)]{QuasipolynomialInverseGowers}
    at least if $n=1$.
    The only point where the proof might depend on the fact that the source is $\Z $ rather than $\Z ^n$ is 
    the formula 
    ``$\Expec _{h\in [\pm N]} \lnorm{\Delta _h f}^{2^2}_{U^s[N]}\ll _s \lnorm{f}^{2^{s+1}}_{U^{s+1}[N]}$''
    towards the end of the proof. 
    But the same fact holds for $f\colon [\pm N]^n \to \C $ easily from the definitions, with an implied constant depending on $n$ as well. 
\end{proof}

\subsection{The equidistribution theory of polynomial sequences}\label{sec:appendix_nilsequences}

\begin{theorem}[{Quantitative Leibman Theorem on $\Zn$ of Leng}]
    \label{thm:Leibman}
    
    Let $0<\delta <1/10$, $M\ge 1$ be positive numbers and $N,k,s,d,n\ge 1$ be positive integers.

    Let $G/\Gamma $ be an $s$-step nilmanifold of dimension $d$, step $ s$, degree $ k$ and complexity $\le M$.
    
    Let $F\colon G/\Gamma \to \bbC $ be an $M$-\Lip function such that there is a continuous homomorphism $\xi \colon G_{(s)}\to \bbR $ of size $\mgn{\xi }\le M/\delta $ such that for all $a\in G_{(s)}$ and $b\in G$, we have 
    \[ F(a b)= e^{2\pi i \xi (a)}\cdot F(b) . \] 

    Let $g\colon \Z ^n\to G$ be a polynomial sequence and 
    suppose we have 
    \[ \mgn{ \Expec _{x\in [\pm N]^n } F(g(x))  } \ge \delta . \] 
    Then either:
    \begin{itemize}
        \item we have $N\le (M/\delta )^{O_{k,n}(d^{O_{k,n}(1)})}$;
        \item we can find an integer $0\le r\le \dim (G/[G,G])$ and horizontal characters 
        $\eta _1,\dots ,\eta _r\colon G\to \R $ of size $\le (M/\delta )^{O_{k,n}(d^{O_{k,n}(1)})}$ such that for all $i=1,\dots ,r$, 
        \[ \lnorm{\eta _i\circ g}_{C^\infty ([\pm N]^n)} \le (M/\delta )^{O_{k,n}(d^{O_{k,n}(1)})} \] 
        and such that if we define $G':= \bigcap _{i=1}^r \ker (\eta _i )$, its $s \ordinalth $ iterated commutator subgroup satisfies:
        \[  G'_{(s)}\subset \ker (\xi ). \] 
        (Recall that the iterated commutator subgroups of $G$ for us are $G_{(1)}:=G$, $G_{(2)}:=[G,G]$, $G_{(3)}:=[G,G_{(2)}]$ and so on.)
    \end{itemize}
\end{theorem}
\begin{proof}
    This is a particular case of \cite[Theorem 4]{EfficientEquidistribution}.
\end{proof}

The following associated Factorization Theorem is also useful. 
The reason why we do not cite the combined statement of Theorem \ref{thm:Leibman} and Proposition \ref{prop:factorization} is that in our argument,
the horizontal characters obtained by Theorem \ref{thm:Leibman} have to go through some preparatory treatments before 
they are fed into Proposition \ref{prop:factorization}.

Before stating Factorization, let us recall a few more notions from 
{\cite[Definitions 1.7, 9.1]{PolynomialOrbits}}.

\begin{definition}\label{def:(M,N)-smooth}
    Let $G/\Gamma $ be a nilmanifold of degree $\le k$ and let $M>10 $ be a positive number.
   
    \begin{enumerate}
        \item 
        A polynomial map $\varepsilon \colon \Zn \to G$ is said to be 
        {\em $(M,[\pm N]^n)$-smooth}
        ($M,N\ge 3$) 
        if for all $x\in [\pm N]^n$ we have 
        $d(\varepsilon (n),1_G)\le M$
        and $d(\varepsilon (x),\varepsilon (x+e_i))\le M/N$ for all $1\le i\le n$,
        where $e_i\in \Zn $ is the $i\ordinalth$ standard vector.
        \item 
        A polynomial map $\gamma \colon \Zn \to G$ is said to be {\em $M$-rational}
        if for every $x\in \Zn$, the coset $\gamma (x)\Gamma \subset G$ contains an element $g_x $ such that $g_x^{M'_x} =1$ for some $0<M'_x\le M$.
        When the \Mal basis has complexity $\le M$,
        modulo raising $M$ to the 
        $\dim (G)^{O_k(1)}$
        st power,
        this is the same as saying that
        its values via the \Mal coordinates $G\cong \R ^{\dim G}$
        are all in $\frac{1}{M'}\Z ^{\dim G}$ for some $0<M'\le M$
        \cite[Lemma A.11]{PolynomialOrbits}.
    \end{enumerate}
\end{definition}

\begin{proposition}[{factorization associated to horizontal characters}]
    \label{prop:factorization}
    Let $G/\Gamma $ be a nilmanifold of dimension $\le d$, degree $\le k$, step $\le s$ and complexity $\le M$.
    Let $g\colon \Zn \to G$ be a polynomial map and $N\ge 1$ be a positive integer.
    
    Given linearly independent horizontal characters $\eta _1,\dots ,\eta _r\colon G\to \R $ of size $\le L$ (for some $L>1$) satisfying 
    $\lnorm{\eta _i\circ g}_{C^\infty ([\pm N]^n )} \le \delta\inv $ for some $0<\delta <1$ and for all $i=1,\dots ,r$, then we have a factorization 
    \begin{align}
        g = \varepsilon \cdot g_1 \cdot \gamma ,
    \end{align}
    where 
    \begin{itemize}
        \item $\varepsilon $ is $((ML/\delta )^{O_{k,n}(d^{O_{k,n}(1)})}, [\pm N]^n )$-smooth;
        \item $g_1$ has values in $G_1:= \bigcap _{i=1}^r \ker (\eta _i )$, which has complexity $\le (ML/\delta )^{O_{k,n}(d^{O_{k,n}(1)})}$;
        \item $\gamma $ is $(ML)^{O_{k,n}(d^{O_{k,n}(1)})}$-rational.
    \end{itemize}
    Furthermore, if $g$ satisfies $g(0)=1_G$, then we can choose $\varepsilon , g_1, \gamma $ such that $\varepsilon (0)=g_1(0)=\gamma (0)=1_G$.
\end{proposition}
\begin{proof}
    This is \cite[Lemma A.1]{EfficientEquidistribution}.
\end{proof}

We need the following general results about polynomial maps $\Z ^n\to \R $ as well.

\begin{proposition}[{polynomial Vinogradov theorem, \cite[Proposition 2.3]{PolynomialOrbitsErratum}}]
    \label{prop:Vinogradov}
    Let $g\colon \Z ^n\to \R $ be a polynomial map of degree $\le k$, and $0<\delta <1$ and $0<\varepsilon < \frac{\delta}{10} $ be positive numbers.

    If the bound 
    \begin{align}
        \lnorm {g(x)}_{\RZ }\le \varepsilon 
    \end{align}
    holds for at least $\delta N^n$ values of $x\in [\pm N]^n$,
    then there is some $0< Q< \delta ^{-O_{k,n}(1)}$ such that 
    \begin{align}
        \lnorm{Qg}_{C^\infty ([\pm N]^n)} < \delta ^{-O_{k,n}(1)} \varepsilon .
    \end{align}
\end{proposition}

\begin{proposition}[{extrapolation, \cite[Lemma 8.4]{PolynomialOrbitsErratum}}]
    \label{prop:extrapolation}
    Let $k,n\ge 1$ be integers.
    Let $q_0,N\ge 1$ be integers, $q_0$ being sufficiently large depending on $k,n$, and $g\colon \Z ^n\to \R $ be a polynomial map of degree $\le k$.

    Let $\undl a = (a_1,\dots ,a_n), \undl b= (b_1,\dots ,b_n)\in \Z ^n $ be vectors satisfying
    \begin{align}
        |a_i|\le q_0N, \quad |b_i|\le q_0 .
    \end{align}
    Then there is $0<q<q_0^{O_{k,n}(1)}$ such that we have 
    \begin{align}
        \lnorm {qg}_{C^\infty ([\pm N]^n)} < q_0^{O_{k,n}(1)} \cdot \lnorm{g(\undl a + \undl b \cdot (-))}_{C^\infty ([\pm N]^n)} 
        .
    \end{align}
\end{proposition}
\begin{proof}
    This statement differs from \cite[Lemma 8.4]{PolynomialOrbitsErratum} in that we allow $|a_i|\le q_0N$
    rather than $|a_i|\le q_0$ as in \loccit 
    But the proof carries over and is in any case easy.
\end{proof}

\begin{proposition}[Fourier expansion along $G_{(s)}$]\label{prop:Fourier-along-Gs}
    Let $G/\Gamma $ be a nilmanifold of dimension $\le d$, degree $\le k$, step $\le s$ and complexity $\le M$.
    Let $F\colon G/\Gamma \to \C $ be a function with \Lip norm $\le L$.
    Then for any $0<\delta <1/100$, there is a decomposition as a function on $G/\Gamma $:
    \begin{align}
        F= \left( \sum _{|\xi |\le (ML /\delta )^{O_{k}(d^{O_{}(1)})}} F_\xi \right)  + E ,
    \end{align}
    where 
    \begin{itemize}
        \item each $\xi \colon G_{(s)}\to \RZ $ is a continuous homomorphism of indicated size;
        \item $F_\xi \colon G/\Gamma \to \C $ has vertical frequency $\xi $, i.e., it satisfies
        \[ F_\xi (a b) = e^{2\pi i \xi (a)}\cdot F_\xi (b) \] 
        for all $a\in G_{(s)}$ and $b\in G$;
        \item $\lnorm{F_\xi }_{\mLip } < (ML/\delta )^{O_k(d^{O(1)})}$;
        \item the remainder term $E\colon G/\Gamma \to \C $ satisfies 
        $\lnorm{E}_\infty \le \delta $.
        \end{itemize}
\end{proposition}
\begin{proof}
    This is \cite[Lemma A.6 and Remark]{EfficientEquidistribution} specialized to $H=G_{(s)}$ 
    combined with the Remark in \loccit 
\end{proof}

Let us also recall the following variant of the Cauchy--Schwarz inequality.

\begin{proposition}[A version of Cauchy--Schwarz]
    \label{prop:Cauchy-Schwarz}
    Let $A,B$ be finite sets and let $f\colon A\to \C $, $F\colon A\times B\to \C $ be functions.
    Then 
    \begin{align}
        \mgn{
            \sum _{a\in A}\sum _{b\in B} f(a)F(a,b) 
            }
        \le 
        \max_{a\in A} |f(a)|
        (\# A)^{1/2}
        \mgn{
            \sum _{a\in A} \sum _{b,b'\in B} F(a,b)\ol{ F(a,b')} 
        }^{1/2}
    \end{align}
\end{proposition}
\begin{proof}
    This is taken from \cite[Lemma A.10]{QuadraticUniformity} 
    but the coefficient is not explicated in \loccit 
    We have to clarify what the coefficient depend on because we apply it unboundedly many times.
    
    By the triangle inequality 
    \begin{align}
        \mgn{
            \sum _{a\in A}\sum _{b\in B} f(a)F(a,b) 
            }
        \le \sum _{a\in A} 
        |f(a)|
        \mgn{
            \sum _{b\in B} F(a,b) 
        }
        \le 
        \max _{a\in A}|f(a)|
        \sum _{a\in A}
        \mgn{\sum _{b\in B} F(a,b) }
        .
    \end{align}
    By \CauSch for the functions 
    $1$ and $a\mapsto \mgn{\sum _{b\in B} F(a,b) }$ on $A$, we have 
    \begin{align}
        \sum _{a\in A}
        \mgn{\sum _{b\in B} F(a,b) }
        &\le 
        (\# A )^{1/2}
        \left(
            \sum _{a\in A}
            \mgn{\sum _{b\in B} F(a,b)}^2
        \right)^{1/2}
        \\
        &=
        (\# A )^{1/2}
        \left(
            \sum _{a\in A}
            \sum _{b,b'\in B} F(a,b)\ol{ F(a,b')}
        \right)^{1/2}
        .
    \end{align}
    These two inequalities imply the assertion of Proposition \ref{prop:Cauchy-Schwarz}.
\end{proof}

%
%

\section{Recollections from algebraic number theory}\label{appendix:algebraic-number-theory}

\subsection{Mertens' theorem}

\begin{proposition}[{Mertens' theorem for number fields; \cite[Theorem 2]{Rosen}}]
    \label{thm:Mertens-Rosen}
    For positive numbers $x>1$ we have:
    \begin{align}
        \prod _{
            \substack{\idealp \\ \Nrm (\idealp )<x}
        }
        \left( 1-\frac{1}{\Nrm (\idealp )} \right) \inv 
        =
        e^{\gamma }\cdot \residue _{s=1}(\zeta _K(s))
        \cdot \log x + O_K(1) .
    \end{align}
    Here, $\gamma =\lim _{n\to +\infty } \left( \sum _{k=1}^n \frac 1k - \log n \right)$ is Euler's constant.
\end{proposition}

\subsection{Sieve theory for ideals}\label{appendix:sieve}

In the proof of Proposition \ref{prop:sieve-theory}, we used part (1) of the following lemma (for the usual context of integers $\Z $).
Part (2) was used when we performed a Vaughan-type decomposition of 
the \Cramer/Siegel models of $\Lambda _K$ in Proposition \ref{prop:Vaughan-decomposition}. 

Let $\Ideals _K ^{\mrm{sq\mathchar`-free}}$ be the set of square-free elements 
(i.e., the products $\idealn = \idealp _1\cdots \idealp _r$ of distinct prime ideals, $r\ge 0$).
Choose a total order ``$<$'' on the set of non-zero prime ideals, compatible with the norm in the sense that 
\begin{align}\label{eq:def-of-compatibility-order-and-norm}
    \idealp_1 < \idealp _2 \te{ implies }\Nrm (\idealp_1)\le \Nrm (\idealp_2).
\end{align}
For an ideal $\idealn $ and a real numbers $z$, let us write $\idealn <z$ (resp.\ $\idn \le z$, \dots)
to mean $\Nrm (\idealn )<z$ (resp.\ $\Nrm (\idn )\le z$, \dots).
For $z>2$, let $\mathfrak P(z)$ be the ideal 
\begin{align}
    \mathfrak P(z):= \prod _{\idealp <z}\idealp .
\end{align}

\begin{lemma}[Fundamental Lemma]
    \label{lem:sieve-theory}
    Let $z>0$, $\kappa >0$ be real numbers and let $D\ge z^{9\kappa +1}$.

    Let $g\colon \Ideals _K^{\mrm{sq\mathchar`-free}}
    \to [0,1)\subset \R $ be a multiplicative function.
    Assume that we have a bound 
    \begin{align}\label{eq:sieve-hyp-bound}
        \prod _{w\le \idealp<z } (1-g(\idealp)) \inv \le M\left( \frac{\log z}{\log w} \right) ^{\kappa } 
    \end{align} 
    for some constant $M>0$ and all $2\le w < z$.

    Let $(a_\idn )_{\idn \in\Ideals_K}$ be a family of non-negative real numbers, zero except for finitely many $\idn $,
    and $X>0$ be a real number.
    For each ideal $\idd \le D$ dividing $\mathfrak P(z)$, let $r_\idd $ be the ``remainder'' defined by the relation
    \begin{align}\label{eq:coefficient-X}
        \sum _{\substack{
            \idn \in\Ideals_K \\ \idd \divides \idn 
            } }
            a_{\idn } 
            =
            X g(\idd ) + r_\idd .
    \end{align} 
        Then:
        \begin{enumerate}
            \item 
            Write $S(z):= \sum _{\substack{
                \idn \in\Ideals _K  
                }}
                a_\idn 
                1[{\idn +\mathfrak P(z)=(1)}]\in \R_{\ge 0}$. 
                We have 
            \[
                S(z)
                    =
                    X\left( \prod _{0<\idealp <z} (1-g(\idealp )) \right)
                    (1+O_K(M^{10} 
                    e^{9\kappa - \frac{\log D}{\log z}} ))
                    + O_K\left( \sum _{
                        \substack{ 1\le \idd \le D \\ \idd |\mathfrak P(z) }
                        } 
                        |r_\idd |
                        \right) .
                        \]
            \item There are functions (depending on $D$)
                $\lambda ^\pm \colon \Ideals _K \to \{ -1,0,1 \}$;
                $\idd \mapsto \lambda ^{\pm }_\idd $
                satisfying:
                \begin{itemize}
                    \item $\lambda _{(1)}^\pm =1$;
                    \item $\lambda ^+$ is supported on $\Ideals_{K,\le D}$;
                    \item for all $\idn \in \Ideals _K$, we have 
                    \[ \sum _{\idd \divides \idn } \lambda _\idd ^- \le 1_{\idn =(1)} \le \sum _{\idd \divides \idn }\lambda _\idd ^+ ;\]
                    \item Write $S^\pm (z ):= \sum _{\idn\in\Ideals_K} a_\idn \sum _{\idd \divides (\idn + \mathfrak P(z))}\lambda ^\pm _\idd $. We have 
                    \begin{align}
                        \mgn{
                            S(z) - S^\pm (z)
                        }
                        <
                        O_K(e^{9\kappa -\frac{\log D}{\log z}  } M^{10} )\cdot 
                        X\cdot \prod _{\idealp <z} (1-g(\idealp ))
                        +
                        \sum _{\substack{
                            \idd \divides \mathfrak P(z) \\ \idd <D }} |r_\idd | .
                    \end{align}
                \end{itemize}
        \end{enumerate}
\end{lemma}
\begin{proof}
    When $K=\Q$, this is \cite[pp.~153--159]{Iwaniec-Kowalski},
    notably (6.19), (6.20), Theorem 6.1, Corollary 6.2 and Lemma 6.3 in \loccit 
    (We are writing $D$ for their $y$.)
    Consult also \cite[\S\S 6.4--6.5]{opera-cribro} for clearer definitions and more careful proofs.
    
    One can check, slightly amuzingly, that the definitions and arguments in the above references continue to make sense and work 
    for any free commutative monoid $\mathscr I$ generated by a countable set, equipped with:
    \begin{itemize}
        \item a multiplicative function (the ``norm'')
            $\Nrm \colon \mathscr I \to \N $
            such that for every $n \in \N $,
            the inverse image $\Nrm \inv (n )\subset \mathscr I$ is finite;
        \item a total order on the set of irreducible elements (``prime ideals'')
            compatible with the norm in the sense of \eqref{eq:def-of-compatibility-order-and-norm}.
    \end{itemize}
    In particular the proofs in the above references apply to $\mathscr I=\Ideals _K$ with the usual norm and the chosen ordering.
\end{proof}

\subsection{Hecke characters}\label{sec:Hecke_chars}

We set up notation relatd to Hecke characters.
See \cite[\S 2.3]{KaiMit} for references.

Let $\Ideles $ be the group of \ideles of $K$.
For non-zero ideals $\idq =\prod _{\idealp } \idealp ^{n_\idealp }\subset \OK $,
denote by $C_K(\idq )$ the following quotient of the \idele class group $\Ideles /K\baci $:
\begin{align}
    C_K(\idq ):= \left( (\KR )\baci \times \bigoplus _{\idealp \divides \idq} K\baci / (1+\idealp^{n_\idealp}\mathcal O_{K,\idealp })
    \times 
    \bigoplus _{\idealp\notdivide\idq }\Z \right) 
    / K\baci .
\end{align}
An inclusion $\idq \subset \idq '$, or $\idq '\divides \idq $, induces a surjection 
$C_K(\idq )\surj C_K(\idq ')$.
When a non-zero fractional ideal $\idc $ is coprime to $\idq $,
we have a well-defined class 
\begin{align}
    [\idc ]\in C_K(\idq ).
\end{align}

A {\it Hecke character} of the number field $K$ is a continuous homomorphism 
from the \idele class group $\Ideles /K\baci $ to the group of complex numbers of absolute value $1$:
\begin{align}
    \psi\colon \Ideles /K\baci \to S^1:=\{ z\in\C\baci \mid |z|=1 \} .
\end{align}
The continuity implies that 
$\psi $ factors through $C_K(\idq )$
for some non-zero ideal $\idq $.
In this case $\psi $ is called a {\it mod $\idq $ Hecke character}.
There is a maximum ideal $\idq $ (or {\it minimum} in terms of divisibility) such that this factorization happens;
we call it the {\it conductor} of $\psi $.

For a non-zero fractional ideal $\ideala $, 
let us write
\begin{align}
    (\ideala /\idq \ideala )\baci 
    &:=
    \{ x\in \ideala /\idq\ida \mid x\te{ generates }\ideala/\idealq\ideala \te{ as an $\OK$-module} \} ,
    \\ 
    \ida ^{\times ,\idq } &:= \ida \times _{\ida /\idq\ida }(\ideala /\idq \ideala )\baci 
    \\
    &= 
    \{ x\in \ideala \mid x\te{ generates }\ideala/\idealq\ideala \te{ as an $\OK$-module} \} .
\end{align}
There is a natural map of sets,
which we also denote by $\incl $: 
\begin{align}
    \incl \colon (\ideala / \idq \ideala )\baci
    \to 
    C_K(\idq ),
\end{align}
induced by the inclusion maps $(\ida / \idq\ida )\baci \to K\baci / (1+\idp ^{n_\idp }\OKp )$
for $\idp \divides \idq $, 
and the constant maps
$(\ida /\idq\ida )\baci \to \{ v_{\idp }(\ida ) \} \subset \Z $ for $\idp \notdivide \idq $.
(We wrote $v_{\idp }(\ida )$ to mean the exponent of $\idp $ in the prime decomposition of $\ida $.)
Let us write 
\begin{align}
    \psi ^{(\ida /\idq\ida )\baci } := \psi \circ \incl \colon (\ida /\idq\ida )\baci 
    \to C_K(\idq )\to S^1 .
\end{align}
Note that the archimedean part of $\psi $ has no contribution to $\psi ^{(\ida /\idq \ida )\baci }$.

By the {\it restriction $\psi ^{\ida }\colon \ida \to \C $ of $\psi $
to $\ida $} we mean the zero-extension of the following composite map:
\begin{align}
    \psi ^{\ida }:
    \ideala ^{\times ,\idq }
    \xto{(\incl , \can )} 
    (\KR )\baci \times (\ideala  / \idq \ideala )\baci 
    \xto{\incl + \incl } 
    C_K(\idq )\xto{\psi } S^1\subset \C .
\end{align}
For $x\in \ideala ^{\times ,\idq }$, the class $[x\inv \ida ]\in C_K(\idq )$ is well defined, and 
by subtracting the diagonal image of $x$ in $C_K(\idq )$, we have 
$\psi ^{\ida }(x)
= \psi ([x\inv \ida ])
= \psi ([x\ida\inv ])\inv $; for $x\in \ida \setminus \ida ^{\times ,\idq }$ we have $\psi^{\ida }(x)=0$.

Suppose that $\psi $ has discrete image (as in the case of a Siegel character) and let $\Omega \subset (\KR )\baci $ be a connected component. 
Then by continuity, the image of the composite 
$\Omega \subset (\KR )\baci \subset \Ideles \xto \psi S^1 $
consists of one point, which we write $\psi (\Omega )\in S^1 $. 
Then we have 
\begin{align}
    \psi ^{\ida }(x) = \psi (\Omega )\cdot \psi ^{(\ida /\idq\ida )\baci }(x)
    \quad \te{ for all $x\in \ida ^{\times ,\idq }\cap \Omega $},
\end{align}
and in particular the restriction of $\psi ^{\ida } $ to $\ida \cap \Omega $ factors through $\ida /\idq\ida $.
(The resulting map $\ida /\idq\ida \to \C $ depends on $\Omega $.)

\section{Subsets with \Lip boundaries and lattice points}\label{appendix:lattice-points}

The intuition that the number of lattice points in a subset of $\R ^n$
should be well approximated by its volume is deeply embedded in our genes.
In this section we recall a quantitative result to this effect.

By convention, subsets of the Euclidean spaces $\R^n$ are given the usual $l^2$-metric 
(also called the Euclidean metric)
and its restriction in this section.

A {\em lattice} $\Lambda $ in $\R^n$ means a discrete subgroup isomorphic to $\Z^n$.
Its {\em determinant} $\det (\Lambda )$ is defined to be the volume of the compact quotient $\R ^n / \Lambda $.
We write $\lambda _1$ for the smallest length of non-zero elements:
\begin{align}
    \lambda _1 := \min \{ |v| : v\in \Lambda\nonzero \} 
    .
\end{align} 
The choice of notation is based on the fact that it is the first of the so-called {\em successive minima} of $\Lambda $.

\begin{definition}[subsets of certain \Lip classes]
    \label{def:Lip_class}
    Let $Z\ge 1$ be an integer and $L> 0$ be a positive real number.
    A subset $S\subset \R ^n$ is said to be of {\em \Lip class 
    $\mLip (n,Z,L)$}
    if there are 
    maps 
    \begin{align}
        \psi _i \colon 
        [0,1]^{n-1} 
        \to \R ^n, 
        \quad i=1,\dots ,Z
    \end{align}
    which have \Lip constant $\le L$
    and cover $S$ in the sense that $S\subset \bigcup _{i=1}^Z \phi _i 
    ([0,1]^{n-1})$. 
\end{definition}
This is 
taken from 
\cite[Definition 2.2]{Widmer2011}.
The same notion is written as  
$\mLip (n,1,Z,L)$ in \cite[Definition 2.1]{Widmer2010}. 

Note that when we have finitely many sets $\Omega_i\subset \R^n$
whose boundary $\partial \Omega_i$ is of \Lip class $\mLip (n,Z_i,L_i)$,
then the intersection $\bigcap _i \Omega_i$ has boundary 
of \Lip class $\mLip (n,\sum _i Z_i,\max _i L_i)$ at worst.

\begin{proposition}\label{prop:Widmer}
     Let $\Lambda \subset \R ^n$ be a lattice.
    \begin{enumerate} 
        \item Let $v_1,\dots ,v_n\in \Lambda $ be a basis.
        Write 
        \begin{align}
            F= \{ a_1v_1+\dots +a_nv_n \mid 0\le a_i<1 \te{ for all $1\le i\le n$} \}
            \subset \R ^n
        \end{align}
        for the associated fundamental domain.
        Let $C\ge 1$ be a positive number such that 
        \begin{align}
            |v_i| \le C \lambda _1  \quad \te{ for all $i=1,\dots ,n$}
            .
        \end{align}
        Then for every subset $S\subset \R^n$ of Lipschitz class 
        $\mLip (n,Z,L)$ with $L\ge \lambda _1$, we have 
        \begin{align}
            \#\{ v\in \Lambda \mid (v+F )\cap S \neq \varnothing \}
            \ll _n 
            Z C^{n^2}\left( \frac{L}{\lambda _1} \right) ^{n-1}
            .
        \end{align}
        \item 
        Let $\Omega \subset \R^n$ be a closed set whose boundary $\partial \Omega $ is of \Lip class $\mLip (n,Z,L)$
        with $L\ge \lambda _1$. Then we have 
        \begin{align}
            \mgn{
                \# (\Omega \cap \Lambda ) - \frac{\vol (\Omega )}{|\det (\Lambda )|}
                }
                \ll_{n}
                Z \left( \frac{L}{\lambda _1} \right) ^{n-1}
            \end{align}
    \end{enumerate}
\end{proposition}
\begin{proof}
    These are specialized forms of \cite[Proposition 5.2 and Corollary 5.3]{Widmer2010}.
\end{proof}

Results about subsets with \Lip boundaries apply in particular to convex sets thanks to the following result (though in this particular case the proofs are easier).

\begin{theorem}\label{thm:convex-Lip-boundary}
    Every convex set in the cube $[- N,N]^n$ have boundary of \Lip class $\mLip (n,1,O_n(N))$.
\end{theorem}
\begin{proof}
    This is 
\cite[Theorem 2.6]{Widmer2011}.
\end{proof}

We are often in a situation where we partition a box $[-N,N]^n$ into small boxes of side length $\varepsilon N$ for some $0<\varepsilon <1$ (with respect to a basis of a lattice).
The next statement is taylored for this situation, which can be obtained by rescaling Proposition \ref{prop:Widmer}.
\begin{corollary}\label{cor:Widmer-rescaled}  
    Let $\Lambda \subset \R ^n$ be a lattice. Choose its basis $v_1,\dots ,v_n$ and define $F\subset \R ^n$, $\lambda _1>0$, $C\ge 1$ as in Proposition \ref{prop:Widmer}.    
    Let $N>1$ and $1/N <\varepsilon <1$.
    Let $S\subset \R^n$ a set of \Lip class $\mLip (n,Z,L_0N )$ with $L_0\ge \varepsilon \lambda _1$.
    Consider the set 
    \begin{align}
        X(\{ v_i \} _i ,N,S,\varepsilon ):= \{ v\in \Lambda \mid (v+\varepsilon N \cdot F )\cap S \neq \varnothing \} .
    \end{align}
    \begin{enumerate} 
        \item We have 
        \begin{align}
            \# X(\{ v_i \} _{1\le i\le n} ,N,S,\varepsilon ) 
        \ll _n 
        Z C^{n^2}\left( \frac{L_0}{\varepsilon\lambda _1} \right) ^{n-1}
        ,
    \end{align}  
    so that  
    \begin{align}
        \# \bigcup _{v \in X(\{ v_i \} _i ,N,S,\varepsilon )}
        (v + \varepsilon N \cdot F )\cap \Lambda 
        &\ll _n 
        Z C^{n^2}\left( \frac{L_0}{\varepsilon\lambda _1} \right) ^{n-1}\cdot (\varepsilon N)^n
        \\ 
        &\ll _{n,Z,C,L_0,\lambda _1} \varepsilon N^{n} .
    \end{align}  
    \item In particular, if $\Omega \subset [-N,N]^n$ is a convex set,
        we have 
    \begin{align}
        \# \bigcup _{v \in X(\{ v_i \} _i ,N,\partial \Omega ,\varepsilon )}
        (v + \varepsilon N \cdot F )\cap \Lambda 
        \ll _n 
        \frac{C^{n^2}}{\lambda _1^{n-1}} \varepsilon N^n
        \ll _{n,C,\lambda _1} \varepsilon N^{n} .
    \end{align}
\end{enumerate}
\end{corollary}

Besides convex sets,
results in this section are used mainly to subsets with the following type of boundary:
\begin{lemma}[The sets with constant norm are Lipschitz]
    \label{lem:bounded-by-norm-Lipschitz}
    Let $N\ge 2$ be a positive real number and set 
    \begin{align}
        \C _{\le N} := \{ z\in \C \mid |z|\le N \}
        .
    \end{align}
    Let $r_1,r_2\ge 0$ be integers and consider the norm function
    \begin{align}
        \Nrm \colon \R ^{r_1}\times \C^{r_2} &\to \R 
        \\ 
        (x,z)=(x_1,\dots ,x_{r_1};z_1,\dots ,z_{r_2})&\mapsto \prod _{i=1}^{r_1}x_i \prod _{j=1}^{r_2} |z_j|^2.
    \end{align}
    Then for any $a\in \R$ the subset 
    \begin{align}
        S_{\le N}(a):= \{ (x,z)\in [- N,N]^{r_1} 
        \times \C ^{r_2}_{\le N} \mid \Nrm (x,z)=a \} 
    \end{align}
    is in \Lip class $\mLip (r_1+2r_2,O_{r_1,r_2}(1),O_{r_1,r_2}(N))$:
    in other words, there are 
    continuous maps $\phi _i$ ($1\le i\le O_{r_1,r_2}(1)$)
    \begin{align}
        [0,1]^{r_1+2r_2-1} 
        \xrightarrow{\phi _i}
        \R ^{r_1}\times \C ^{r_2}
    \end{align}
    with \Lip constants $\le O_{r_1,r_2}(N)$ such that 
    \begin{align}
        S_{\le N}(a) \subset \bigcup _{i} \phi _i ([0,1]^{r_1+2r_2-1}). 
    \end{align}
\end{lemma}
\begin{proof}
    This is elementary but we include a proof for the sake of completeness.
    The norm map $\Nrm $ factors as 
    \begin{equation}
        \begin{array}{rcccl}
            \R^{r_1}\times \C^{r_2}
            &\to 
            &\R^{r_1}\times \R^{r_2}_{\ge 0}
            &\xto{\overline \Nrm }
            &\R 
            \\ 
            (x,z)&\mapsto & (x;|z_1|,\dots ,|z_{r_2}|),
            \\ 
            &&(x,y)&\mapsto &\prod\limits _{i=1}^{r_1}x_i \prod\limits _{j=1}^{r_2} y_j^2
        \end{array}
    \end{equation}
    and for any $a\in R$ if we set 
    \begin{align}
        \overline S_{\le N}(a) := \{ (x,y)\in [-N,N]^{r_1}_{\le N}\times [0,N]^{r_2} \mid \overline \Nrm  (x,y)=a \}
        \subset \R^{r_1}\times \R^{r_2}_{\ge 0}
        ,
    \end{align}
    we can recover $S_{\le N}(a)$ by 
    \begin{multline}
        S_{\le N}(a) = 
        \\ 
        \{ (x;y_1e^{2\pi i \theta _1},\dots ,y_{r_2}e^{2\pi i \theta _{r_2}} ) \in \R^{r_1}\times \C^{r_2} \mid (x,y)\in \overline S (a),(\theta _j)_{1\le j\le r_2}\in [0,1]^{r_2} \}
        .
    \end{multline}
    This reduces the problem to showing that 
    $\overline S_{\le N}(a) $ is of \Lip class $\mLip (r_1+r_2,O_{r_1,r_2}(1),O_{r_1,r_2}(N))$.
    While it is easy to construct an explicit covering with this \Lip property,
    here we choose to describe a quick proof using the case of the boundary of a convex set.
    
    By the symmetry $\pm x_i$ one can restrict attention to $\overline S_{\le N}(a)\cap [0,N]^{r_1+r_2}$
    and assume $a\ge 0$.
    Since $\overline N $ is a convex function on $[0,N]^{r_1+r_2}$, the set 
    $\Omega _{\le N}(a):=  \{ x\in [0,N]^{r_1+r_2} \mid \overline \Nrm  (x)\ge a \} $
    is a convex set.
    As $\partial \Omega _{\le N}(a)$ 
    is known to be of \Lip class $\mLip (r_1+r_2,1,O_{r_1+r_2}(N))$ by \cite[Theorem 2.6]{Widmer2011}
    and we have $\overline S_{\le N}(a)\subset \partial\Omega _{\le N}(a)$,
    we conclude that $\overline S_{\le N}(a)$ is of the same \Lip class.
    This completes the proof.
\end{proof}


The next proposition gives a Fourier approximation of the indicator function of a set with \Lip boundary.
In this paper we use it only when $Z,L_0 \ll _n 1$.

\begin{proposition}\label{prop:Fourier}
    Let $\Omega \subset [\pm N]^n $ be a set with \Lip boundary of class $\mLip (n,Z,L_0N)$ 
    and $1<M,Y<N$ be positive integers.
    Then 
    \begin{enumerate} 
        \item 
        There is a decomposition of the indicator function $1_\Omega \colon [\pm N] ^n\cap \Z ^n\to \C $ of $\Omega $ into a sum of the form
        \begin{align}\label{eq:decomposition-of-indicator-function}
            1_\Omega  
            = \left( \sum _{\theta \in \Theta } c_\theta \cdot e^{2\pi i \theta (-) } \right)
            + G +H,
        \end{align}
        where 
        \begin{itemize}
            \item $\Theta := \Hom (\Z ^n ,\frac{1}{Y}\Z /\Z ) $ is a finite set of size $Y^n$;
            \item 
            $|c_\theta |\le 1$ for all $\theta\in\Theta $;
            \item 
            $\lnorm G _\infty \le 1$ and 
            $\#\oname{supp} (G) \ll _n  \frac{ZL_0^{n-1}}{M}N^n$;
            \item 
            $\lnorm H _{\infty } = O_n(M\frac{\log Y}{Y})$.
        \end{itemize}
        \item Furthermore, suppose a lattice $\Lambda \subset \Zn $ 
        and a residue class $a\in \Zn / \Lambda $ is given. 
        Then the indicator function $1_{\Omega \cap (a+\Lambda )}$ is written as a sum of the same form \eqref{eq:decomposition-of-indicator-function},
        except now $\Theta := \{ \theta \in \Hom (\Z ^n,\QZ ) \ : \ \theta (\Lambda )\subset \frac{1}{Y}\Z /\Z \}$ has $Y^n \cdot \# (\Zn /\Lambda )$ elements.
    \end{enumerate}
\end{proposition}
\begin{proof}
    Approximate the indicator function $1_\Omega $ by the function 
    \[ x\mapsto \max \left\{ 1-\frac{M}{N} \cdot d(x,\Omega ) , 0 \right\} \]
    of \Lip constant $\le \frac{M}{N}$.
    The difference of this function and $1_\Omega $ is supported on the set  
    $\{ x\in [\pm N]^n\cap \Z ^n \ : \ 0<d(x,\Omega )< \frac{N}{M}  \}$,
    which has $\ll _n \frac{ZL_0^{n-1}}{M} {N^n}$ elements by Proposition \ref{prop:Widmer} (1)
    (set $\Lambda :=\Z ^n$, $\lambda _1=1$, $C=1$, $\varepsilon :=1/M $).
    To obtain part (1), apply the Fourier approximation of this \Lip function in 
    \cite[Lemma A.9]{QuadraticUniformity} (reproduced as Proposition \ref{proposition:Fourier} below) after rescaling by $2N+2$.

    For part (2),
    by Fourier analysis on the finite abelian group $\Zn /\Lambda $,
    the indicator function $1_{a+\Lambda }\colon \Zn \to \C $ can be written as 
    a $\C $-linear combination
    \begin{align}
        1_{a+\Lambda } = \sum _{\eta \in \widehat{\Zn / \Lambda }} c_\eta \cdot e^{2\pi i \eta (-)} ,
    \end{align}
    where $\widehat{\Zn / \Lambda }=\Hom (\Zn / \Lambda ,\QZ )$ is the dual abelian group.
    To conclude, multiply \eqref{eq:decomposition-of-indicator-function} by $1_{a+\Lambda }$.
\end{proof}

\section{The von Neumann theorem for the Gowers norms}\label{sec:vonNeumann}

The (generalized) von Neumann theorem for the Gowers norms 
gives a criterion for when a sum of the form 
$\sum _{\bm x\in [-N,N]^n\cap \Omega }\prod _{i=1}^t f_i (\phi _i (\bm x )) $ is small, 
where $\Omega \subset \R^n$ is a convex set, $\phi _i\colon \Z^d \to \Z^n$ are affine-linear maps and 
$f_i\colon \Z^n\to \C $ are functions.

We state it (Theorem \ref{thm:vonNeumann}) after explaining 
the {\em complexity}%
\footnote{
    Though \cite[Exercise 1.3.23 (p.~59)]{HigherOrderFourier}
    states the von Neumann theorem for an arbitrary finite abelian group $G$, the definition of the complexity of a set of affine-linear maps $\psi _i \colon G^d \to G$ only makes sense 
    when each $\psi _i $ is the sum of multiplication-by-integer maps $G \to G$, which is always the case when $G$ is cyclic.
    (Notice that the definition of the complexity in \loccit requires a lift $\psi _{i,\Z } \colon \Z ^d \to \Z  $ of $\psi _i$ in the sense that $\psi _{i,\Z }\otimes _{\Z }G = \psi _i $.)
} 
$s\ge 0$ of the set of linear maps $\dot\phi _1,\dots ,\dot\phi _t$.



In this section we will use bold face letters such as $\bm x$ to denote tuples instead of underlined ones $\undl x$ for better readability.

\subsection{Preliminary notions}

\begin{definition}[complexity of a set of linear maps]
    \label{def:complexity}
    Let $\psi _1,\dots ,\psi _t\colon \bZ ^d\to \bZ ^n $
    linear maps.
    \begin{enumerate}
        \item
        They are said to have
        {\it (\CauSch) complexity $\le s$
        at the index $i\in \br{1,\dots ,t}$}
        if there is a partition $\br{1,\dots ,\widehat i , \dots , t} = \bigcup _{k=1}^{s+1} C_k $
        such that for all $k\in \br{1,\dots ,s+1}$,
        the restriction of $\psi _i$:
        \[\psi _i\colon  \bigcap _{j\in C_k} \ker (\psi _j) \to \bZ ^n \]
        has rank $n$.
        (For this to hold, of course each $\ker (\psi _i)$ must have rank $\ge n$.)
    \item They are said to have {\em complexity} $\le s$ if they do at every index.
    \item
    They are said to be {\it in $s$-normal form at the index $i$}
    if there are disjoint sets $B_1,\dots ,B_{s+1}$ (or less members) of standard vectors of $\bZ ^d$
    such that
    \begin{enumerate}
        \item for all $j\neq i$, there is a $1\le k\le s+1$ such that $\psi _j$ is the zero map on $\bZ ^{B_k}$;
        \item for all $1\le k\le s+1$, the restriction of $\psi _i$:
        \[
            \psi _i \colon \bZ ^{B_k} \to \bZ ^n
        \]
        has rank $n$.
    \end{enumerate}
    \item There are obvious $\Fp $-linear analogues of the above notions for linear maps $\psi _1,\dots ,\psi _t\colon \Fp ^d \to \Fp ^n$. 
    \end{enumerate}
\end{definition}

Note that if the maps $\br{ \psi _j}_{1\le j\le t}$ are in $s$-normal form at $i$, then they have \CauSch complexity $\le s$ at $i$
because we can take $C_k:= \br{j \middle| {\psi _j }_{| \bZ ^{B_k}} =0 }$.

A set of $s+1$ linear maps $\Z ^d \to \Z ^n$ either has \CauSch complexity $\le s$ or does not have finite complexity.
The linear maps have finite complexity if and only if $\psi _i (\ker \psi _j)\subset \Z ^n$ has finite cokernel for all distinct $i,j\in \br{1,\dots ,s+1}$.

An $\OK $-linear map $\OK ^d \to \OK $ has finite cokernel if and only if it is non-zero.
Therefore a set of $s+1$ linear maps $\OK ^d \to \OK $ has \CauSch complexity $\le s$ (via an isomorphism $\OK \cong \Z ^n$) if and only if their kernels do not contain each other.

The next fact can be used to reduce problems involving linear maps of \CauSch complexity $\le s$ 
to the case of $s$-normal forms.

    \begin{proposition}\label{proposition:can-be-made-into-normal-form}
        Let $\psi _1,\dots ,\psi _t\colon \bZ ^d\to \bZ ^n $ be
        linear maps
        having
        \CauSch complexity $\le s$
        at some index $i\in \br{1,\dots ,t}$.
    
        Then there is a surjection
        $\bZ ^{d'}\to \bZ ^d$
        with $d' \le (s+2)d$
        such that the composites $\psi _i '\colon \bZ ^{d'}\to \bZ ^d \to \bZ ^n $
        ($i=1,\dots ,t$) are in $s$-normal form at the index $i$.
    
        The same holds for $\Fp $-linear maps $\psi _1,\dots ,\psi _t\colon \Fp ^d \to \Fp ^n$.
    \end{proposition}
    \begin{proof}
        For each index set $C_k$ from the definition of complexity, let $d_k$ ($\le d$) be the rank of the $\Z$-module $\bigcap _{j\in C_k}\ker (\psi _j)$
        and consider the map $\bZ ^{d_k}\xrightarrow[\cong]{} \bigcap _{j\in C_k}\ker (\psi _j)\subset \bZ ^d$. 
        Then the surjection
        \[
            \bZ ^{d}\oplus \bigoplus _{k=1}^{s+1} \bZ ^{d_k}
            \longrightarrow \bZ ^{d}
        \]
        satisfies the desired properties.
        The proof for the $\Fp $-linear case is the same.
    \end{proof}

    \subsection{The statement of the generalized von Neumann theorem}

    The next is a direct generalization of
    \cite[Proposition 7.1]{LinearEquations}
    to the case where the target of the linear maps is $\mathbb Z^n$
    rather than $\mathbb Z$.
    Let the notation $\lnorm {\psi } _N \le L$ mean that the affine-linear map 
    $\psi \colon \Z ^d \to \Z ^n$ has linear coefficients with $l^\infty $-norm $\le L$ 
    and constant term with $l^\infty $-norm $\le NL$.
    When $\psi $ is a linear map, we simply write $\lnorm \psi \le L$. 
    
    \newcommand{\statementVonNeumann}{
        Let $t,d,L,n>0$ be positive integer parameters.
    
        Let $\psi _1,\dots ,\psi _t \colon \bZ ^d \to \bZ ^n$ be affine-linear maps,
        satisfying $\lnorm {\psi _i} _N \le L$ and 
        having \CauSch complexity $\le s$.
        Let $f_1,\dots ,f_t\colon \bZ ^n\to \bC $ be functions such that 
        \[ \lnorm{ f_i} _{U^{s+1}([-(d+1)LN,(d+1)LN]^n)} \le \delta \] 
        for some $i\in \{ 1,\dots ,t \}$ and positive $\delta >0$.
        Let $\Omega \subset [-N,N]^d \subset \bR ^d$ be a convex body.
    
        Then for any positive values $M,Y>0$ we have
        \begin{equation}
            \mgn{ \Expec _{[-N,N]^d} 1_\Omega  \prod _{i=1}^t (f_i\circ \psi _i) }
            \ll _{d,n,s}
            L^d
            \paren{\delta Y^d 
            +\frac 1M +M\frac{\log Y}Y
            }.  
        \end{equation}
        }
        
    \begin{theorem}
        [generalized von Neumann theorem]
        \label{thm:vonNeumann}
        \statementVonNeumann
    In particular, setting $M:=\left( Y/ \log Y\right)^{1/2}$ and $Y:= (\delta \inv )^{1/(d+\frac 12)}$, the left hand side is bounded as 
        $
            \ll _{d,n,s} 
            L^d
            \cdot \delta ^{1/(2d+2)}
            .$
    \end{theorem}

    \begin{strategy}\label{strategy}
    The proof will go through several reduction steps:
    \begin{enumerate}
        \item We go modulo $p$, where $p$ is a prime number comparable in size to $(d+1)LN$;
        this step account for the factor $L^d$ in the right hand side.
        \item In order to dispose of the indicator function $1_{(\Omega \cap \Z^d \oname{mod}\, p )}$, we first replace it with functions with \Lip\ constants $\le M$.
        This account for the term $\frac 1M$.
        \item We write the \Lip\ function(s) as a Fourier series containing at most $O_d(X^d)$ terms and an error term.
        The error term accounts for the term $M\frac{\log X}{X}$
        and the number of terms in the Fourier series accounts for the factor $X^d$.
        \item The trigonometric functions in the Fourier series can be absorbed into the functions $f_i$
        so that we can ignore them.
        \item We are reduced to the estimation of averages of the form $\Expec _{\Fp ^d} \prod _{i=1}^t f_i\circ \psi _i $.
        This will be bounded in size by $\delta $.
     \end{enumerate}
    \end{strategy}

The following fact about the Gowers norm will be used in an intermediate step, Proposition \ref{proposition:mod-p-without-Lip}.
    
\begin{proposition}[Gowers--\CauSch inequality]\label{prop:Gowers-Cauchy-Schwarz}
    Let $(X_\alpha )_{\alpha \in A}$ be a family of non-empty finite sets.
    For each subset $B\subset A$, write $X_B:=\prod _{\alpha \in B} X_\alpha$
    and suppose that we are given a $1$-bounded function $f_B\colon X_B\to \C $.
    Denote also by $f_B$ the composite with the projection $X_A \to X_B\xto{f_B}\C $.
    \begin{enumerate}[(1)]
        \item 
        We have
    \begin{multline}\label{eq:Gowers-Cauchy-Schwarz}
            \left\vert \Expec _{x\in X_A} \prod _{B\subset A} f_B(x) 
        \right\vert \le \lnorm{f_A}_{\square (X_A)}
        \\ 
        := \left( \Expec _{x^{(0)},x^{(1)} \in X_A} 
        \left[ \prod _{\omega \in \{ 0,1 \} ^{A}}
        \mcal{C}^{|\omega |}f_A((x^{(\omega (\alpha ))}_\alpha )_{\alpha \in A}) \right] \right) ^{1/2^{|A|}} . 
    \end{multline}  
        \item Suppose that $X_\alpha $ are all (finite) abelian groups and that there are a function $f$ on another finite abelian group $Y$
        and surjective homomorphisms $\pi _\alpha \colon X_\alpha \surj Y $ such that 
        $f_A$ can be factored as 
        $X_A \xto{\sum _\alpha \pi_\alpha } Y \xto f \C $.
        Then the right-hand side of \eqref{eq:Gowers-Cauchy-Schwarz} equals $\lnorm{f}_{U^{|A|}(Y)}$.
    \end{enumerate}
\end{proposition}
\begin{proof}
    Part (1) is \cite[(B.7) on p.~1814]{LinearEquations}.
    Part (2) follows by elementary change of variables.
    We only mention the change of variables and leave the verification of the equality to the reader.
    Consider 
    $y:= \sum _{\alpha \in A} \pi _\alpha (x_\alpha ^{(0)})\in Y$ and 
    \[ \bm h = (h_\alpha )_{\alpha \in A} := \Bigl( \pi _\alpha (x_\alpha ^{(1)}-x_\alpha ^{(0)}) \Bigr) _{\alpha \in A} \in Y . \]
    Then one can verify 
    \[ \Expec _{x^{(0)},x^{(1)} \in X_A} 
    \left[ \prod _{\omega \in \{ 0,1 \} ^{A}}
        \mcal{C}^{|\omega |}f_A((x^{(\omega (\alpha ))}_\alpha )_{\alpha \in A}) \right] 
        = 
        \Expec _{\substack{y\in Y \\ \bm h\in Y^{A} }} 
        \left[ \prod _{\omega \in \{ 0,1\} ^A } \mcal{C}^{|\omega |} f(y + \omega \cdot \bm h ) \right] . \]    
        By definition, this is $\lnorm{f}_{U^{|A|}(Y)}^{2^{|A|}}$.
\end{proof}

\subsection{Proof of the von Neumann theorem}

We first treat 
an $\Fp $-linear version of Proposition \ref{thm:vonNeumann}
where 
the convex set $\Omega $ is not present

\begin{proposition}\label{proposition:mod-p-without-Lip}
    Let $t,d,n>0$ positive integer parameters.
    
    Let $\psi _1,\dots ,\psi _t\colon \Fp ^d \to \Fp ^n$ be affine-linear maps
    having \CauSch\ complexity $\le s$ for some integer $s\ge 0$
    and 
    $f_1,\dots ,f_t\colon \Fp ^n\to \bC $ be $1$-bounded functions.

    Assume $\lnorm{ f_i} _{U^{s+1}(\Fp ^n)} \le \delta $ for some $i\in \{ 1,\dots , t \} $ and positive $\delta >0$.

Then 
we have
    \begin{equation}
        \mgn{ \Expec _{\Fp ^d} \prod _{i=1}^t (f_i\circ \psi _i) }
        \le 
            \delta  
        .
    \end{equation}
\end{proposition}
\begin{proof}
    By symmetry in the indices $i=1,\dots ,t$, we may assume $\lnorm{ f_1} _{U^{s+1}(\Fp ^n)} \le \delta $.
    By Proposition \ref{proposition:can-be-made-into-normal-form}, 
    we may assume $\psi _i$ are in $s$-normal form at the index $1$
    because pre-composing a surjection $\Fp ^{d'}\surj \Fp ^d$ to the maps $\psi _i$ does not change the average in question, and $\delta $ is totally unrelated to $\psi _i$'s.

    Using the disjoint sets
    $B_1,\dots ,B_{s+1}$ (or less members) of standard vectors from Definition \ref{def:complexity},
    write
    \[
        X_k\coloneqq  \Fp ^{B_k} \text{ and }X\coloneqq \bigoplus _{k=1}^{s+1} X_k .
    \]
    Also, let $V$ be the subspace spanned by those standard vectors not in $\bigcup _{k=1}^{s+1} B_k $.
    We have a decomposition
    \[
         \Fp ^d = X_1\oplus \dots \oplus X_{s+1} \oplus V
    \]
    and we know
    \begin{enumerate}
        \item for all $j\in  \br{2,\dots ,t}$, the map $\dot\psi _j$ is zero on some $X_k$;
        \item for all $1\le k\le s+1$, the restriction $\dot\psi _1 \colon X_k\to \Fp ^n$ is surjective. 
    \end{enumerate}
    Let us use the symbol $x_k$ to denote elements of $X_k$ and $v$ to denote elements of $V$.
    We denote elements of $X=X_1\oplus \dots \oplus X_{s+1}$ by the symbol $\bm x =(x_1,\dots ,x_{s+1})$.

We can rewrite the average $\Expec _{\Fp ^d}$ as a nested one $\Expec _{v\in V} \Expec _{\bm x\in X}$.
By the triangle inequality we have
\begin{equation}\label{eq:rewrite-the-average}
    \mgn{\Expec _{\Fp ^d} \prod _{j=1}^t f_j\circ \psi _j  }
    \le \Expec _{v\in V}
    \mgn{
        \Expec _{\bm x \in X}
        \prod _{j=1}^t f_j( \psi _j (\bm x +v))
    }
    .
\end{equation}
Note that since $\psi _j$ is affine-linear, we have $f_j(\psi _j (\bm x +v))=f_j(\dot\psi _j (\bm x)+\psi _j(v))$.
For each $v\in V$, 
we apply the \CauSchGow\ inequality (Proposition \ref{prop:Gowers-Cauchy-Schwarz}) 
to the situation
\begin{itemize}
    \item $A=\br{1,\dots ,s+1}$;
    for each subset $B\subset A$ we consider $X_B= \prod_{k \in B} X_k$;
    \item consider the subset $\Omega (j)\subset A$ given by 
    \[ \Omega (j):= \Bigl\{ k \in \{ 1,\dots ,s+1 \}  \midsep \psi _{j}| X_k \neq 0 \Bigr\} ; \]
    then by condition (1) above, we know $\Omega (j)\subsetneq A$ for $j\neq 1$; 
    the map $\bm x\mapsto f_j(\dot\psi _j (\bm x) +\psi _j(v)) )$ can be seen as a function on $X_{\Omega (j)}$;
    \item set $f_A\colon X_A\to \C $ to be $f_1( \dot\psi _1 (-)+\psi _1(v) )$;
    for each $B\subsetneq A$ set 
    \[ f_B := \prod _{\substack{ j\te{ such that}\\ \Omega (j)=B}} f_j(\dot\psi _j (-)+\psi _j(v)) \colon X_B \to \C . \]
    \item as has been noted in (2) above, the restriction $\dot\psi _1 \colon X_k \surj \Fp ^n$ is surjective for each $1\le k\le s+1$.
\end{itemize}
This implies that the previous value \eqref{eq:rewrite-the-average} is bounded as
\begin{equation}\label{eq:after-CauSchGow-inequality}
    \le 
    \Expec _{v\in V}
        \lnorm{f_1 ((-)+\psi _1(v))
        } _{
            U^{s+1}(\Fp ^n)
            }
    .
\end{equation}
As the Gowers norm is translation-invariant, this equals $\lnorm{f_1}_{U^{s+1}(\Fp ^n)}$,
which is $\le \delta $ by assumption.
\end{proof}

Now we want to upgrade Proposition \ref{proposition:mod-p-without-Lip}
by adding a \Lip\ function into the picture.
Embed the abelian group $\F _p$ into $\RZ $ by
\[
    \begin{array}{rcl}
        \Fp = \bZ / p\bZ &\cong & \paren{\frac 1p \bZ } / \bZ
        \hookrightarrow \RZ
        \\
        x \mod p\bZ  &\leftrightarrow & \frac xp \mod \bZ
    \end{array},
\]
through which each function on $(\RZ )^d$ 
gives a function on $\Fp ^d$. 

Equip $(\RZ )^d$ with the usual Riemannian metric. 
This allows us to talk about \Lip\ constants of functions on $(\RZ )^d$.

We use the following Fourier approximation of \Lip\ functions.

\begin{proposition}[{see \cite[Lemma A.9]{QuadraticUniformity}}]
    \label{proposition:Fourier}
    Let $M,Y >0$ be arbitrary positive numbers.
    If $F\colon \paren{\RZ }^d \to \bC $ is a 1-bounded \Lip\ function 
    with \Lip\ constant $\le M$,
    then there is a decomposition
    \[
        F  =
        \left( \sum _{\theta \in \Theta } c_{\theta }\cdot \etheta \right) 
        + H,
    \]
    where $\Theta := [-Y,Y]^d \subset \Z ^d \cong \Hom ((\RZ )^d, \RZ ) $ is a subset of size $O_d(Y^d)$,
    the coefficients satisfy $\mgn{c_{\theta }} \le 1 $,
    and $H\colon (\RZ )^d \to \bC $ is a function that satisfies
    $\lnorm H _{\infty } \le O_d\paren{M\cdot \frac {\log Y}{Y}}$.
\end{proposition}

\begin{proposition}\label{proposition:mod-p-with-Lip}
    Let $t,d,n>0 $ be positive integers.

    Let $M,Y>0$ be further positive parameters.
    Let $F\colon (\RZ )^d\to \bC $ be a 1-bounded (i.e., values have absolute value $\le 1$) \Lip\ function with \Lip\ constant $\le M$.

    Let $\psi _1,\dots ,\psi _t \colon \Fp ^d \to \Fp ^n$ be surjective affine-linear maps
    with \CauSch\ complexity $\le s$.
    Let $f_1,\dots ,f_t\colon \Fp ^n \to \bC $ be 1-bounded functions.
    Assume $\Gnorm{f_i}{\Fp ^n} \le \delta $ for some index $i$ and $ \delta >0$.

    Then we have
    \[
        \mgn{
            \Expec _{\Fp ^d} F \cdot \prod _{i=1}^t f_i \circ (\psi _i )
            }
            \ll _d
            \delta X^d 
            +M\cdot\frac{\log X}{X}
            .
    \]
\end{proposition}

\begin{proof}[Proof of Proposition \ref{proposition:mod-p-with-Lip}]
Taking the approximation of $F$ in Proposition \ref{proposition:Fourier}, we have
\begin{align}\label{eq:after-Fourier}
    \mgn{\Expec _{\Fp ^d} F\cdot \prod _{i=1}^t (f_i\circ \psi _i )}
    \le &
    \sum _{\theta \in \Theta }
    \mgn{
        \Expec _{\Fp ^d} \etheta \prod _{i=1}^t (f_i\circ \psi _i )
    }
    \\
    &+
    \mgn{
        \Expec _{\Fp ^d}  H\cdot \prod _{i=1}^t (f_i\circ \psi _i )
    }
    .
\end{align}
The second term in the right-hand side can be bounded by the triangle inequality and pointwise bound:
\begin{align}
    \mgn{
        \Expec _{\Fp ^d}  H\cdot \prod _{i=1}^t (f_i\circ \psi _i )
    }
    \le 
    \Expec _{x\in \Fp ^d} \mgn{H (x) \cdot \prod _{i=1}^t (f_i\circ \psi _i (x) ) }
    \ll _d &
    M\frac{\log X}X .
\end{align}
So the second term in the right-hand side of \eqref{eq:after-Fourier}
is an $O_d\paren{M\frac{\log X}{X}}$. 

Next, we show that the first term in the right-hand side of \eqref{eq:after-Fourier} is an $O_d(\delta X^d)$.
There are two cases to consider depending on $\theta $.
To formulate them, note that 
the homomorphism $\theta \colon \paren{\frac 1p \bZ /\bZ }^d \to \RZ $ 
has image in $\frac 1p \bZ /\bZ \cong \Fp $ and 
can be regarded as 
an $\Fp $-linear map $\theta \colon \Fp ^d \to \Fp $.

\paragraph{Case 1} $\ker (\theta )\not\supseteq \bigcap _{i=1}^t \ker (\dot\psi _i)$.

Take a 1-dimensional $\Fp $-subspace $\ell \subset \bigcap _{i=1}^t \ker (\dot\psi _i)$ not containedd in $\ker (\theta )$.
Take a direct sum decomposition $\Fp ^d =\ell \oplus V$.
By writing the average $\Expec _{\Fp ^d}$ as the successive average over $\bmv\in V$ and $\bmw\in \ell $,
we have
\begin{align}
    \Expec _{\Fp ^d} \paren{
        \etheta \cdot \prod _{i=1}^t (f_i\circ \psi _i) 
        }
    &=
        \Expec _{\bmv \in V}\paren{
            \Expec _{\bmw \in \ell } e^{2\pi i \theta (\bmv+\bmw )} \prod _{i=1}^t f_i (\psi _{i} (\bmv +\bmw ) )
            }
    .
\end{align}
By the choice of $\ell $, the quantity $\psi _i(\bmv +\bmw )$ does not depend on $\bmw $ and $e^{2\pi i \theta (\bmv+\bmw )}$ does depend on $\bmw $.
It follows that each of the inner averages $\Expec _{\bmw \in \ell }$ is zero, completing the treatment of Case 1.

\paragraph{Case 2}
$\ker (\theta )\supseteq \bigcap _{i=1}^t \ker (\dot\psi _i)$.

Note the following elementary fact.
\begin{claim}\label{claim:elementary-linear-sh**}
    Under the assumption of Case 2 (that $\ker (\theta )\supseteq \bigcap _{i=1}^t \ker (\dot\psi _i)$), there are linear maps $\eta _i\colon \Fp ^n \to \Fp $ such that
    $\theta = \sum _{i=1}^t \eta _i\circ \dot\psi _i $.
\end{claim}
\begin{proof}
    Let $\pr _j \colon \Fp ^n \to \Fp $ be the $j\ordinalth$ projection.
    The assumption is equivalent to the membership relation 
    $\theta \in \langle \pr _j \circ \dot\psi _i \rangle _{i,j}$
    where the right-hand side is the $\Fp$-span in the dual space $(\Fp^n)^* $.
    So there are $a_{i,j}\in \Fp $ such that
    $\theta = \sum _{i,j} a_{i,j} \pr _j \circ \dot\psi _i$.
    Then $\eta _i := \sum _{j} a_{i,j} \pr _j$ give desired maps, completing the proof of Claim \ref{claim:elementary-linear-sh**}.
\end{proof}

Using Claim \ref{claim:elementary-linear-sh**}, consider the functions
$f_i'\colon \Fp ^n\to \bC $ defined by
\[
    f_i' (\bm y) =  f_i(\bm y)\cdot \exp \paren{2\pi i\,  \eta _i (\bm y)/p}.
\]
For every $\bm x\in \Fp ^d$ we have
\begin{align}
    \prod _{i=1}^t f_i'(\psi _i(\bm x))
    =&
    \paren{\prod _{i=1}^t f_i(\psi _i(\bm x))}
    \cdot \exp \paren{\frac{2\pi i}{p} \sum _{i=1}^t \eta _i (\psi _i (\bm x)) }
.\end{align}
By the choice of $\eta _i$, we have
\[ \sum _{i=1}^t \eta _i (\psi _i (\bm x))
= \sum _{i=1}^t  \eta _i (\dot\psi _i (\bm x) +\psi _i(0))
= \theta (\bm x ) + \sum _{i=1}^t \eta _i(\psi _i(0)) \te{ in }\Fp (\overset{\times 1/p}{\inj} \RZ ).\]
The second term is independent of $\bm x$.
Therefore we have (with a complex number $u$ of absolute value $1$):
\[
    \prod _{i=1}^t f_i'(\psi _i(\bm x))
    =
    \left( \prod _{i=1}^t f_i(\psi _i(\bm x))\right)  \cdot e^{2\pi i \theta (\bm x)} \cdot u
    .
\]
We apply Proposition \ref{proposition:mod-p-without-Lip} to the functions $f_i'$, noting that  
$\Gnorm{f_i'}{\Fp ^n} = \Gnorm{f_i}{\Fp ^n}$ by \cite[(B.11) in Appendix B]{LinearEquations}.
This gives us 
\begin{align}
    \mgn{
        \Expec _{\Fp ^d} \etheta \prod _{i=1}^t (f_i\circ \psi _i )
        }
        =&
        \mgn{
            \Expec _{\Fp ^d} \prod _{i=1}^t (f_i'\circ \psi _i) 
            }
        \le
        \delta .
        \end{align}

        By Cases 1 and 2, for every $\theta $ we have obtained
        $\mgn{
            \Expec _{\Fp ^d} \etheta \prod _{i=1}^t (f_i\circ \psi _i )
            }
            \le \delta $.
    Noting $\mgn \Theta =O_d(Y^d)$, it follows that
    \[
        \sum _{\theta \in \Theta } \mgn{
            \Expec _{\Fp ^d} \etheta \prod _{i=1}^t (f_i\circ \psi _i )
            }
            = O_d( Y^d ) \cdot \delta ,
    \]
completing the proof of Proposition \ref{proposition:mod-p-with-Lip}.
\end{proof}

Now we consider the situation where the sum is taken over a given convex body
$\Omega \subset (\RZ )^d$.
As always, we embed $\Fp ^d\inj (\RZ )^d$.

\begin{proposition}\label{proposition:mod-p-with-convex}
    Let $t,d,n>0 $ be positive integers.

    Let $\Omega \subset [-1/3,1/3]^d $ be a convex body
    which we regard as a subset of $(\RZ )^d$.
    Let $1_\Omega \colon (\RZ )^d\to \br{0,1}$ be its indicator function.

    Let $\psi _1,\dots ,\psi _t \colon \bZ ^d \to \bZ ^n$ be surjective affine-linear maps
        with \CauSch\ complexity $\le s$.
    Let $f_1,\dots ,f_t\colon \Fp ^n \to \bC $ be 1-bounded functions. 
    
    Assume $\Gnorm{f_i}{\Fp ^n} \le \delta $ for some index $i$ and $ \delta >0$.

    Then, for any positive numbers $M,Y>0$,
    we have
    \[
        \mgn{
            \Expec _{\Fp ^d} 1_\Omega  \cdot \prod _{i=1}^t f_i \circ (\psi _i \otimes \Fp )
            }
            \ll _d
            \delta Y^d 
             + M \frac{\log Y}Y
            +  \frac 1M
            .
    \]
\end{proposition}

\begin{proof}
    There is a decomposition
    \[
        1_\Omega  = F_M + G_M
    \]
    where $F_M$, $G_M$ are $1$-bounded functions
    such that $F_M$ has \Lip\ constant $\le M$
    and $G_M$ has support with volume $O_d(1/M)$ in $(\RZ )^d$.
    (For example, writing $d$ for the metric on $(\RZ )^d$, we can take $F_M (x):= \min \br{0,1-M\cdot d(\Omega ,x)}$,
    $G_M:=1_\Omega -F_M$ as in \cite[p.~1820]{LinearEquations}.)

    By the triangle inequality, we have
    \begin{equation}\label{eq:F_M-and-G_M}
        \mgn{
            \Expec _{\Fp ^d} 1_\Omega  \cdot \prod _{i=1}^t (f_i \circ \psi _i )
            }
            \le
            \mgn{
                \Expec _{\Fp ^d} F_M \cdot \prod _{i=1}^t (f_i \circ \psi _i )
            }
            +
            \mgn{
            \Expec _{\Fp ^d} G_M \cdot \prod _{i=1}^t (f_i \circ \psi _i )
            }
            .
        \end{equation}
    By Proposition \ref{proposition:mod-p-with-Lip},
    the first term in the right-hand side is
    $\ll _d \delta Y^d  +M \frac{\log Y}Y$.

    The second term in the right-hand side is bounded by the triangle inequality as:
    \begin{equation}\label{eq:bounded-by-G_M'}
        \le \mgn{
            \Expec _{\Fp^d} G_M
            }
            \le \mathrm{vol}(\Supp G_M)=O_d(1/M)
            .
        \end{equation}

    Combining the bounds that we have obtained for the right-hand side of \eqref{eq:F_M-and-G_M}, we conclude
    \[
        \mgn{
            \Expec _{\Fp ^d} 1_\Omega  \cdot \prod _{i=1}^t f_i \circ (\psi _i \otimes \Fp )
            }
            \ll _d
            \delta Y^d 
            +
            M \frac{\log Y}Y
            +  \frac 1M.
    \]
    This completes the proof of Proposition \ref{proposition:mod-p-with-convex}.
\end{proof}

\begin{remark}\label{rem:more-weakly-than-convexity1}
    Assume, more weakly than convexity, that $\Omega \subset [-1/3,1/3]^{d}$ regarded as a subset of $\R ^d$ has boundary of \Lip class $\mLip (d,Z,L)$ for some $Z,L\ge 1$.
    (See Definition \ref{def:Lip_class} for this condition.)
    Then the conclusion of Proposition \ref{proposition:mod-p-with-convex}
    stays valid if we replace the term $1/M$ by $ZL^{d-1}/M$.
    This is because the bound $\vol (\Supp (G'_M))\ll _d 1/M$
    gets replaced by $\ll _d ZL^{d-1}/M$
    as any set of \Lip class $\mLip (d,Z,L)$ has $1/M$-neighborhood of volume $\ll _d ZL^{d-1}/M$.
    However, we do not need this generality.
\end{remark}

We are now ready to prove our main result of this section.
For the convenience of reference, let us reproduce verbatim the 
statement:

\theoremstyle{plain}
\newtheorem*{vonNeumannTheorem}{Theorem \ref{thm:vonNeumann}}
\begin{vonNeumannTheorem}[reproduced]
\statementVonNeumann
\end{vonNeumannTheorem}
\begin{proof}[Proof of Theorem \ref{thm:vonNeumann}]
By the assumption $\lnorm{\psi _i}_N \le L$, we have
\[ \psi _i ([\pm N]^d)  \subset [\pm (d+1)LN ]^n =: \mcal B . \]
Let $p$ be a prime number satisfying $2(d+1)LN +2\le p \le 4(d+1)LN +4$,
which exists by Bertrand's postulate (Chebyshev's theorem).

The composite map $\mcal B \subset \Z ^n \surj \Fp ^n $ is injective. Thus we can define the functions 
$1_{\mcal B}f_i $ on $\Fp ^n$ whose value at $x\in \Fp ^n $ is $0$ if $x$ is not in the image of $\mcal B$, and $f_i(\wtil x)$ if $x$ is the image of an $\wtil x \in \mcal B$ (necessarily unique).
By \cite[Lemma B.5]{LinearEquations}, we have
\begin{equation}\label{eq:Gowers-norms-comparable}
    \Gnorm{1_{\mcal B}f_i}{\Fp ^n} = O_{s,n}(1) \Gnorm{f_i}{[\pm (d+1)LN]^n} .
\end{equation}
As $[\pm N]^d \inj \Fp ^d$ is also an injection, we can regard $\Omega \cap \bZ ^d$ as a subset of $\Fp ^d$.
Let us write $\Omega \cap \Fp ^d$ when it is convenient to distinguish.
By the commutative diagram
\[
    \xymatrix{
\Omega \cap \bZ ^d \ar@{}|{\subset }[r]
            \ar@{}|{\Vert}[d]
    & [\pm N]^d\ar[r]^(0.4){\psi _i}
            \ar@{^{(}->}[d]
    & \mcal B  
    \ar@{^{(}->}[d]^\iota
            \ar[dr]^{f_i}
            \ar@{}|{\subset }[r]
    & \Z ^n
    \\
\Omega \cap \Fp ^d \ar@{}|{\subset}[r]
    &\Fp ^d \ar[r]^{\psi _i\otimes \Fp }
    & \Fp ^n \ar[r]^{1_{\mcal B}f_i}
    & \bC ,
    }
\]
we can rewrite the average in question in the following way.
\begin{align}\label{eq:from-Z-to-Fp}
    \mgn{ \Expec _{[-N,N]^d} 1_\Omega  \prod _{i=1}^t (f_i\circ \psi _i) }
    =&
    \frac{1}{(2N+1)^d}
    \mgn{
    \sum _{\bmx \in \Omega  } \prod _{i=1}^t f_i (\psi _i (\bmx ))
    } 
    \\
    =&
    \frac{1}{(2N+1)^d}
    \mgn{
    \sum _{\bmx \in \Omega  } \prod _{i=1}^t 1_{\mcal B}f_i (\psi _i\otimes \Fp  (\bmx ))
    } 
    \\
    =&
    \frac{p^d}{(2N+1)^d}
    \mgn{
        \Expec _{\Fp ^d} 1_{\Omega \cap \Fp ^d} \prod _{i=1}^t (1_{\mcal B}f_i )\circ (\psi _i\otimes \Fp )
    } 
    \\
    \ll _{d}&
    L^d
    \mgn{
        \Expec _{\Fp ^d} 1_{\Omega \cap \Fp ^d} \prod _{i=1}^t (1_{\mcal B}f_i )\circ (\psi _i\otimes \Fp )
    } 
    ,
    \end{align}
where the last inequality follows from the choice of $p\ll _d LN $.
Let 
\[ 
    \delta _{\Fp } := \max _{1\le i\le t }\Gnorm{1_{\mcal B}f_i}{\Fp ^n}.
\]
The right-most term in \eqref{eq:from-Z-to-Fp} can be bounded thanks to Proposition \ref{proposition:mod-p-with-convex}:
\[
    \ll _{d}
    \paren{
        \delta _{\Fp } Y^d 
        +\frac 1M + M\frac{\log Y}Y
    } 
    .
\]
By \eqref{eq:Gowers-norms-comparable}, we have $\delta _{\Fp } 
\ll_{n,s} 
\delta $. 
This completes the proof of Theorem \ref{thm:vonNeumann}.
\end{proof}

\begin{remark}\label{rem:more-weakly-than-convexity2}
    If we assume, more weakly than convexity, that $\Omega \subset [-N,N]^d$ has boundary of \Lip class $\mLip (d,Z,LN)$
    for some $Z,L\ge 1$ (though this generality is not needed for us), then the conclusion of Theorem \ref{thm:vonNeumann}
    stays valid if we replace the term $1/M$ by $ZL^{d-1}/M$.
    See Remark \ref{rem:more-weakly-than-convexity1}.
\end{remark}

\section{Verification of the localized variant of main results}\label{sec:justification-prime-values-localized}

In this appendix we supply a proof of Theorem \ref{thm:prime-values-S}, which is about simultaneous prime values in $\OKS $.

The idea the following:
suppose for example that we want to estimate the sum $\sum _{n\in [1,N]} \Lambda _{\Z [1/6]} (n)$,
where $\Lambda _{\Z [1/6]}$ is a variant of the \vMang function which assigns to $n$ the value $\log p$ if $n = 2^a3^b p^k $ for some $a,b\ge 0$, $k\ge 1$ and prime $p\neq 2,3$, and $0$ otherwise. 
We can morally estimate the sum as follows, where in practice the sum in $a,b\ge 0$ should be truncated at some large enough threshold:
\begin{align}\label{eq:idea-Lambda-1/6}
    \sum _{n\in [1,N]} \Lambda _{\Z [1/6]} (n)
    &=
    \sum _{a,b\ge 0} 
    \left( 
        \sum _{\substack{
            n\in [1,2^{-a}3^{-b}N]\cap \Z  
            \\ 
            \te{coprime to }6
            } }
    \Lambda _{\Z [1/6]} (2^a3^b n) \right)
    \\ 
    &\fallingdotseq 
    \sum _{a,b\ge 0} 
    \left( \sum _{n\in [1,2^{-a}3^{-b}N]}
    \Lambda ( n) \right)
    ,
\end{align}
then we apply Prime Number Theorem:
\begin{align}
    &\fallingdotseq 
    \sum _{a,b\ge 0} 
    2^{-a}3^{-b}N 
    \\ 
    &= 
    N\sum _{a\ge 0} 2^{-a} \sum _{b\ge 0} 3^{-b}
    =N\cdot  2\cdot \frac{3}{2} = 3N.
\end{align}
In general one finds 
$\sum _{n\in [1,N]} \Lambda _{\Z [1/f]} (n) \fallingdotseq (f / {\varphi (f)})\cdot N $.


\subsection{Preliminary discussions}

Under the setting of Theorem \ref{thm:prime-values-S},
let $\ida _0 \subset \idas $ be the (finitely generated) $\OK $-submodule generated by the subset $\bigcup _{i=1}^t \psi _i(\Z ^d) $.
When condition \eqref{eq:Schinzel-condition} is satisfied, we have $\ida_0[S\inv ]=\idas $.
As the statement of Theorem \ref{thm:prime-values-S} only depends on $\idas $, we may assume that $\ida _0=\ida $, in particular the maps $\psi _i$ land in $\ida \subset \idas $.
By the torsion cokernel assumptions, the set of maps 
\begin{align}
    \psi _i\colon \Z^d \to \ida \quad (i=1,\dots ,t) 
\end{align}
has finite complexity (say $s\ge 1$)
in the sense of Appendix \ref{sec:vonNeumann}.

We want to repeat the arguments of \S \ref{sec:prime-values} for the function $\Lambda ^{\ida}_{\OKS }$ in place of $\Lambda ^{\ida }_K$.

\subsection{The \Cramer model in the localized setting}\label{sec:Cramer-OKS}

Let us continue to use the quantities 
$\exp ((\log N)^{1/100})< Q\exp ((\log N)^{1/3}) $ 
and $P(Q)=\prod _{p<Q} p \in \N $, which vary with $N>1$.
We may assume that $N\gg _S 1$ is large enough that $Q$ is larger than all the residue characteristics of the primes $\idp \in S$.

\begin{definition}
    Let $\Lambda _{\OKS ,Q}\colon \Ideals _{\OKS }\cup \{ (0)\} \to \R_{\ge 0}$ be the {\em $Q$-\Cramer model} of the \vMang function defined by 
    \begin{align}
        \idc \mapsto 
        \begin{cases}
            \disp \frac{P(Q)^n}{\totient (P(Q))}\cdot \frac{1}{\residuezeta }
            & \te{if $\idc +(P(Q))=(1)$ in $\OKS $},
            \\ 
            0 &\te{if not}. 
        \end{cases}
    \end{align}
    For a non-zero fractional ideal $\ida $ of $\OK $, define $\Lambda ^{\ida }_{\OKS ,Q }\colon \ida \to \R _{\ge 0}$ 
    by $x\mapsto \Lambda _{\OKS ,Q}(x\ida \inv )$.
\end{definition}
We are not using the decoration ``\Cramer'' anymore because we will not introduce other models.

Our trick, as we sugggested in \eqref{eq:idea-Lambda-1/6}, is to deduce everything from the case of $\OK $ using the following sum description.
Note that for our purposes, we may replace the \vMang function $\Lambda _{\OKS }^{\ida }$
by the function $\Lambda _{\OKS }^{\prime \ida }$
which is defined using the modification $\Lambda '_{\OKS} $ of $\Lambda _{\OKS}$ where we set the values at higher prime powers $\idp ^{k}$ ($k\ge 2$) to be zero because the contribution of such prime ideals is negligible.

For $\idq \in \IdealsK $, let us say it is {\em supported in $S$} and write 
\begin{align}
    |\idq | \subset S
\end{align}
when all the prime divisors of $\idq $ are in $S$.
We say $\idq $ is {\em coprime to} $S$ if no prime divisor of $\idq $ is in $S$.

\begin{lemma}\label{lem:Lambda-OKS-as-sum}
    For all non-zero fractional ideal $\ida $ over $\OK $,
    the following 
    equality of functions on $\ida $ holds:
    \begin{align}\label{eq:Cramer-as-sum}
        \Lambda ^{\ida }_{\OKS ,Q} 
        &= 
            \sum _{
                \substack{
                    \idq \in \IdealsK \\ |\idq | \subset S 
                }
                } 
                1_{\idq\ida }\Lambda ^{\idq\ida }_{\mCramer ,Q}
                .
    \end{align}
    
    Also, $\Lambda ^{\prime \ida }_{\OKS }$ admits a decomposition:
    \begin{align}\label{eq:vMang-as-sum}
        \Lambda ^{\prime \ida }_{\OKS } 
        &= 
        \left( 
        \sum _{
            \substack{
                \idq \in \IdealsK \\ |\idq | \subset S 
            }
            } 
            1_{\idq\ida }\Lambda ^{\prime \idq\ida }_{K}
            \right) 
            +H
            ,
    \end{align}
    where $H$ satisfies $\lnorm H _\infty \ll _{K,S} 1$ and $ \Supp (H)\subset \{ x\in\ida \midsep |x\ida\inv |\subset S \}$ so that $\# ([\pm N]^n_\ida \cap \Supp (H))\ll _K (\log N)^{|S|+n-1}$.

    Furthermore, when we restrict the domains of the functions 
    $\Lambda ^{\ida }_{\OKS ,Q} ,\Lambda ^{\prime \ida }_{\OKS } $ 
    to the cube $[\pm N]^n_\ida $,
    then
    the sums in \eqref{eq:Cramer-as-sum}\eqref{eq:vMang-as-sum} can be restricted to those $\idq $'s having norms 
    $\ll_K N^n$ without changing the values. 
    
\end{lemma}
\begin{proof}
    For $x\in \ida $, 
    For the \Cramer model, the value of both sides of \eqref{eq:Cramer-as-sum} are non-zero precisely when the ideal $ x\ida\inv + P(Q) $ of $\OK $ is supported on $S$ and the non-zero values are both $\frac{P(Q)^n}{\totient (P(Q))}$.

    For the equality \eqref{eq:vMang-as-sum}, 
    the value $\Lambda ^{\prime \ida }_{\OKS }(x)$ at $x\in \ida $ is non-zero precisely when 
    $x\ida\inv $ has the form $\idq \idp $ where $\idq $ is an ideal of $\OK $ satisfying $|\idq |\subset S$ and $\idp $ is a prime ideal away from $S$.
    When this is the case, the value of $ 1_{\idq\ida }\Lambda ^{\idq\ida }_K $ at $x$ is non-zero precisely for this $\idq $ and the non-zero values of both sides coincide and equal $\log \Nrm (\idp )$.
    The sum in the right-hand side is non-zero precisely when $x\ida\inv $ has the form $\idq\idp $ where $\idq $ is supported on $S$ and $\idp $ is {\em any} prime ideal.
    The event that the left-hand side is zero but the sum in the right-hand side is non-zero happens precisely when $x\ida\inv $ is a non-unit ideal supported on $S$; 
    the value of the sum in the right-hand side is then $\sum _{\idp \in S \te{ dividing }x\ida\inv } \log \Nrm (\idp )$, which is $\le |S|\max _{\idp\in S} \log \Nrm (\idp )$.
    Define $H$ to take this value at such $x$ and $0$ elsewhere.

For the bound of the size of $\# \Supp (H)$,
    there are $\ll _K (\log N)^{|S|} $ ideals supported on $S$ and having norm $\ll_K N^n $.
    For each such ideal, there are $\ll _K (\log N)^{n-1} $ points corresponding to it in the cube 
    $[\pm N]^n_\ida $ by the consideration of the $\OK\baci $-action.

    The claim about the range of the sums $\sum _\idq $ holds because 
    the existence of a non-zero element $x\in [\pm N]^n_{\ida }\cap \idq\ida$ implies 
    $\Nrm (\idq )\le \Nrm (x\ida\inv )\ll _K N^n$.
\end{proof}

\subsection{Localized \vMang and \Cramer are close}\label{sec:Fubini}

Now we are in a position to state and prove an intermediate step toward Theorem \ref{thm:prime-values-S}, namely the comparison of the localized \vMang function to its \Cramer model:

\begin{proposition}\label{prop:vMang-close-to-Cramer-OKS}
    For any $N\gg _{s,K,S}1$, any non-zero fractional ideal $\ida $ and any positive number $A>1$,
    we have 
    \begin{align}
        \lnorm{
                \Lambda ^{\ida }_{\OKS }
                -
                \Lambda ^{\ida }_{\OKS ,Q}
        }_{U^{s+1} [\pm N]^n_{\ida }  }
        \\ 
        \ll _{s,K,S,A}
        (\log N)^{-A}
        .
    \end{align}
\end{proposition}

For the proof, the following unnormalized Gowers norm is more convenient.
Let $f\colon Z\to \C $ be a function on an abelian group $Z$ and $A\subset Z$ be a finite subset.
Its {\em unnormalized Gowers uniformity norm} $\lnorm{ f} _{\wtil U^{s+1}(A)}\ge 0$ is defined by:
\begin{align}
    (\lnorm{ f} _{\wtil U^{s+1}(A)} )^{2^{s+1}}
    := 
    \sum _{x\in Z,\undl h\in Z^{s+1}} 
    \prod _{\omega \in \{ 0,1 \}^{s+1}} 
    \calC ^{|\omega |} (1_A f) (x+\omega\cdot \undl h).
\end{align}
The relationship of the normalized (Definition \ref{def:Gowers-norm}) and unnormalized norms is:
\begin{align}
    \lnorm{f}_{U^{s+1}(A)} = \lnorm{f}_{\wtil U^{s+1}(A) } / \lnorm{1}_{\wtil U^{s+1}(A)} .
\end{align}

Let us record a crude upper bound for $\lnorm{f}_{\wtil U^{s+1}(A)}$.
For the value $(1_A f)(x+\omega \cdot \underline h )$ to be non-zero for all $\omega \in \{ 0,1\}^{s+1}$
(hence {\it a fortiori} for all $\omega $ having only one non-zero entry),
we have to have $x\in A\cap \Supp (f)$ and $\underline h \in ( (A\cap \Supp (f)) -x)^{s+1}$.
This implies that there are at most $|A\cap \Supp (f)|^{s+2}$ choices for the tuple $(x,underline h)$.
It follows that 
\begin{equation}\label{eq:crude-bound-unnormalized-Gowers}
    \lnorm{f}_{\wtil U^{s+1}(A)} \le \lnorm{f}_\infty \cdot |A\cap \Supp (f)|^{(s+2)/2^{s+1}}
    \quad (\le \lnorm f _\infty \cdot |A\cap \Supp (f)|).
\end{equation}
When $A=\Omega \subset \Z^n$ comes from a convex set in $\R^n$ with non-empty interior, 
one can easily see that the set 
\[ \{ (x,\undl h)\in (\Z^{n})^{s+2} \midsep \forall \omega \in \{ 0,1\}^{s+1},\ x+\omega \cdot \undl h \in \Omega \} \] 
also comes from a convex set in $(\R ^n)^{s+2}$ with non-empty interior.
It follows that 
\begin{equation}\label{eq:unnormalized-norm-of-1}
    \lnorm{1}_{\wtil U^{s+1}(N\Omega \cap \Zn )} \asymp _{s,\Omega } N^{n(s+2)/2^{s+1}}.
\end{equation}
In particular we have 
 $   \lnorm{1}_{\wtil U^{s+1}([\pm N]^n)} \asymp _{s,n} N^{n(s+2)/2^{s+1}}.$

\begin{proof}[Proof of Proposition \ref{prop:vMang-close-to-Cramer-OKS}]
    Our goal is to prove 
    \[
        \lnorm{
            \Lambda ^{\ida }_{\OKS }
            -
            \Lambda ^{\ida }_{\OKS ,Q}
        }_{\wtil U^{s+1} [\pm N]^n_{\ida }  }
        \ll _{s,K,S,A}
        N^{n(s+2)/2^{s+1}} (\log N)^{-A}
        .
    \]
    By the triangle inequality applied to the sum decomposition in Proposition \ref{lem:Lambda-OKS-as-sum},
    we have
    \begin{align}\label{eq:triangle-inequality-localized}
        &\lnorm{
            \Lambda ^{\ida }_{\OKS }
            -
            \Lambda ^{\ida }_{\OKS ,Q}
        }_{\wtil U^{s+1} [\pm N]^n_{\ida }  }
        \\ 
        &\le 
        \sum _{
            \substack{
                \idq \in \IdealsK \\ |\idp | \subset S \\ \Nrm (\idq )\ll _K N^n
            }
            } 
        \lnorm{
            1_{\idq\ida }
            \left(
                \Lambda ^{\prime \idq\ida }_{K}
                -
                \Lambda ^{\idq\ida }_{\mCramer ,Q}
            \right)
        }_{\wtil U^{s+1} [\pm N]^n_{\ida }  }
        +
        \lnorm{H}_{\wtil U^{s+1} [\pm N]^n_{\ida }  }
        .
    \end{align}
By the crude bound \eqref{eq:crude-bound-unnormalized-Gowers} coupled with the knowledge about $\lnorm H _\infty $ and $\Supp (H)$ from Lemma \ref{lem:Lambda-OKS-as-sum},    
one sees that $\lnorm H _{\wtil U^{s+1} [\pm N]^n_{\ida }  } $ is negligible.    

With respect to a norm-length compatible basis of $\idq\ida $, we have 
$[\pm N]^n_{\ida }\cap \idq\ida \subset [\pm O_K(1) \frac{N}{\Nrm(\idq)^{1/n}}]^n_{\idq\ida }$.
By the definitions of (un)normalized Gowers norms,
Corollary \ref{cor:vonMang-Cramer-close!}
and 
\eqref{eq:crude-bound-unnormalized-Gowers}, 
we have the bounds (mind the difference of $\wtil U^{s+1}$ and $U^{s+1}$)
\begin{multline}    
    \lnorm{
        1_{\idq\ida }
        \left(
            \Lambda ^{\prime \idq\ida }_{K}
            -
            \Lambda ^{\idq\ida }_{\mCramer ,Q}
            \right)
            }_{\wtil U^{s+1} [\pm N]^n_{\ida }  } 
            \\ 
            = \lnorm{1}_{\wtil U^{s+1}[\pm N]^n_{\ida }\cap \idq\ida }\cdot \lnorm{
                1_{[\pm N]^n_{\ida }\cap \idq\ida}
                \left(
                    \Lambda ^{\prime \idq\ida }_{K}
                    -
                    \Lambda ^{\idq\ida }_{\mCramer ,Q}
                    \right)
                    }_{U^{s+1} [\pm O_K(1) \frac{N}{\Nrm(\idq)^{1/n}}]^n_{\idq\ida }  }
    \\
                    \ll_{s,K,S,A}
    \begin{cases}
        \left( \frac{N^n}{\Nrm (\idq )} \right) ^{(s+2)/2^{s+1}} (\log N)^{-A}
        & \te{if } \Nrm (\idq ) \le N^{n/2}, 
        \\ 
        (N^{n/2})^{(s+2)/2^{s+1}} 
        \log N 
        & \te{if } N^{n/2}\le \Nrm (\idq ) \ll_K N^{n}.
    \end{cases}
                \end{multline} 
Thus the sum in \eqref{eq:triangle-inequality-localized} is bounded as 
\begin{multline}
    \ll _{s,K,S,A}
    (N^n)^{(s+2)/2^{s+1}}(\log N)^{-A}
    \left( \sum _{
        \substack{
            \idq \in \IdealsK \\ |\idq | \subset S \\ \Nrm (\idq )< N^{n/2}
        }
        }
    \left( \frac{1}{\Nrm (\idq )} \right)^{(s+2)/2^{s+1}}
    \right)
    \\ 
    +
    (N^{n/2})^{(s+2)/2^{s+1}}\log N 
    \left(
        \sum _{\substack{\idq \\ |\idq |\subset S \\ N^{n/2}\le \Nrm (\idq ) \ll_K N^n }}
        1
    \right)
    .
\end{multline}
The first sum $\sum _\idq $ is $\ll _{K,S} 1$ and the second sum $\sum _\idq $ is $\ll _{K} (\log N )^{|S|}$.
Therefore the entire expression is $\ll_{s,K,S,A} N^{n(s+2)/2^{s+1}} (\log N)^{-A }$,
completing the proof.    
\end{proof}
    
\begin{remark}
    Note how useful it was in the above proof that the thresholds and constants in Corollary \ref{cor:vonMang-Cramer-close!} are uniform in $\ida $ (in the current context, $\idq\ida $) as long as we stick to norm-length compatible bases.
\end{remark}

\begin{proof}[Proof of Theorem \ref{thm:prime-values-S}]
    The proof is similar to that of Theorem \ref{thm:prime-values}.
    By the generalized von Neumann theorem (Theorem \ref{thm:vonNeumann}) applicable by Proposition \ref{prop:vMang-close-to-Cramer-OKS}, 
    we are reduced to the evaluation of the sum 
    \begin{equation}
        \sum _{x\in \Omega \cap [\pm N]^n_\ida } 
        \left( 
        \prod _{i=1}^t \Lambda ^{\ida }_{\OKS ,Q} (\psi _i (x))
        \right) .
    \end{equation}
    Namely we are going to establish an analogue of Proposition \ref{prop:sieve-theory}.
    For this purpose, 
    we apply the general machinery of sieve theory (Lemma \ref{lem:sieve-theory}) to the monoid $\Ideals _{\OKS }$ equipped with the norm function $\idd \mapsto \Nrm (\idd):= |\OKS / \idd |$.
    
    For each square-free ideal $\idd $ of $\OKS $, set 
    \[ 
        g(\idd ):= 
        \left\{ x\in (\Z /\Nrm (\idd ))^d 
        \midsep 
        \te{the ideal} 
        \prod _{i=1}^t (\psi _i (x)\ida \inv ) 
        \te{ is divisible by } \idd
        \right\} 
        / \Nrm (\idd )^d 
        .
    \]
Note that this defining condition depends only on the residue class of $x$ modulo $\Nrm (\idd )$
so $g(\idd )$ is actually well defined.
By Chinese Remainder Theorem, $g$ is multiplicative.
Also, by the definitions $1-g(\idp )= \left( \frac{\totient (\idp )}{\Nrm (\idp )} \right) ^t \beta _{\idp }$.

Next, for ideals (not necessarily square-free) $\idn $ of $\OKS $, define 
\begin{align}
    \Sigma _\idn 
    &:= \{ x\in \Omega \cap \Z ^d \midsep 
    \te{the equality of ideals }
    \prod _{i=1}^t (\psi _i (x)\ida \inv ) = \idn 
    \te{ holds}
    \} , 
    \\ 
    a_{\idn } 
    &:= \# \Sigma _\idn  \in \N _{\ge 0}.
\end{align}
For each square-free ideal $\idd $ dividing 
\[ \mathfrak P (Q) := 
\prod _{\substack{\idp \te{ prime in }\OKS \\ \Nrm (\idp ) <Q }} \idp ,
\] 
we need to approximate the value 
\begin{align}
    &\sum _{\substack{\idn \in \Ideals _{\OKS
    } \\ \te{multiple of }\idd }}
    a_{\idn }
    \\ 
    &= 
    \#
    \left\{ x\in \Omega \cap \Z^d \midsep \te{the ideal }\prod _{i=1}^t (\psi _i (x)\ida \inv ) = 
    \te{ is divisible by }\idd \te{ in }\OKS \right\} .
\end{align}
By packing $\Omega \cap \Z^d $ by translates of the cube $[ 1, \Nrm (\idd )]^d$ and treating the contribution of the boundary as an error,
one sees that this is approximated as 
\begin{align}
    = \vol (\Omega ) \cdot g(\idd ) 
    + O_{\Psi ,Z,L} (N^{d-1} \Nrm (\idd ) ) .
\end{align}
Then Lemma \ref{lem:sieve-theory} gives an asymptotic formula for
\begin{multline}
    S = \sum _{\substack{\idn \in \Ideals _{\OKS } \\ \te{coprime to } \mathfrak P (Q) }}
    a_{\idn }
    \\ =
    \{ x\in \Omega \cap \Z^d \midsep 
    \te{the ideal }
    \prod _{i=1}^t (\psi _i (x)\ida \inv ) \te{ is coprime to } 
    \mathfrak P (Q) 
    \te{ in }\OKS 
    \} ,
\end{multline}
which by definitions equals (where $\idp $ runs over the same range as above)
\begin{align}
    (\residue _{s=1}\zeta _K(s))^t 
    \prod _{\idp }
    \left( 
        \frac{\totient (\idp )}{\Nrm (\idp )} 
        \right)^t
        \sum _{x\in \Omega \cap \Z^d} 
        \left( \prod _{i=1}^t \Lambda ^{\ida }_{\OKS ,Q} (\psi _i (x))
        \right) 
        .
\end{align}
The main term we get is (where $\idp $ runs over the same range again)
\begin{equation}
    \vol (\Omega )\prod _{ \idp 
        }
        \left( \frac{\totient (\idp )}{\Nrm (\idp )} \right) ^t \beta _{\idp } ,
\end{equation}
and if we take $D:= N^{0.4}$, the error term is as asserted.
\end{proof}

\bibliographystyle{amsplain}
\bibliography{WataruBibtex.bib}

\end{document}